\documentclass[oneside,12pt,a4paper,openany]{report}
\usepackage[Rejne]{fncychap}

\usepackage{fancyhdr} 
\usepackage[printonlyused,withpage]{acronym}

\usepackage{arabtex}
\usepackage{utf8}

\usepackage[LAE,T1]{fontenc}

\usepackage{multirow}
\usepackage{float}
 
 \usepackage[figuresright]{rotating}
 
\usepackage{lscape} 
\usepackage{pgfplots}
\pgfplotsset{compat=1.8}
\usetikzlibrary{pgfplots.statistics}
\usepgfplotslibrary{patchplots,colormaps}
\pgfplotsset{compat=1.9} \usepackage{pgfplots}
\usepackage{setspace}
\usepackage{hyperref}
\hypersetup{
	colorlinks=true,
	linkcolor= blue,
	filecolor=black,
	urlcolor=black,
	citecolor=black,
}
\usepackage[flushleft]{threeparttable}
\usepackage{caption}
\usepackage{pgfplots}
\usepgfplotslibrary{patchplots,colormaps}
\pgfplotsset{compat=1.9} 
\usepackage{amsfonts}
\usepackage{lscape} 
\usepackage{commath}
\usepackage{pdflscape}
\usepackage{graphicx}
\usepackage{pdflscape}
\usepackage{amsmath}

\usepackage{xcolor}

\usepackage{mdframed}
\usepackage{float}
\usepackage{rotating}
\usepackage[linesnumbered,ruled,vlined]{algorithm2e}

\usepackage{tabularx}

\usepackage{amsmath}
\usepackage{amsthm}
\usepackage{afterpage}

\usepackage{bbding}
\usepackage{pifont}
\usepackage{wasysym}

\theoremstyle{definition}
\newtheorem{definition}{Definition}[section]

\mdfdefinestyle{MyFrame}{%
	linecolor=black,
	outerlinewidth=2pt,
	innertopmargin=4pt,
	innerbottommargin=4pt,
	innerrightmargin=4pt,
	innerleftmargin=4pt,
	leftmargin =4pt,
	rightmargin = 4pt
}

\usepackage{setspace}

\usepackage{hyperref}
\hypersetup{
	colorlinks=true,
	linkcolor= black,
	filecolor=black,
	urlcolor=black,
	citecolor=black,
}

\usepackage{geometry}
\usepackage[numbers]{natbib}

\DeclareMathOperator*{\argminA}{arg\,min} 

\newcommand*\rot{\rotatebox{45}}

\begin{document}

\pagenumbering{gobble}

	\begin{titlepage}
		\vspace*{-3cm}
\begin{figure}[!htb]
		\begin{minipage}{0.2\textwidth}
		\centering
		\includegraphics[width=1.0\linewidth]{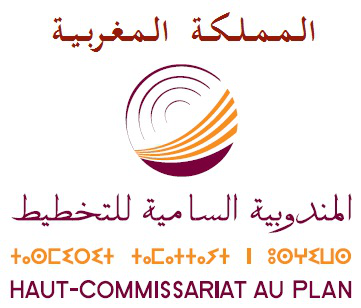}
	\end{minipage}\hfill
	\begin{minipage}{0.2\textwidth}
		\centering
		\includegraphics[width=1.0\linewidth]{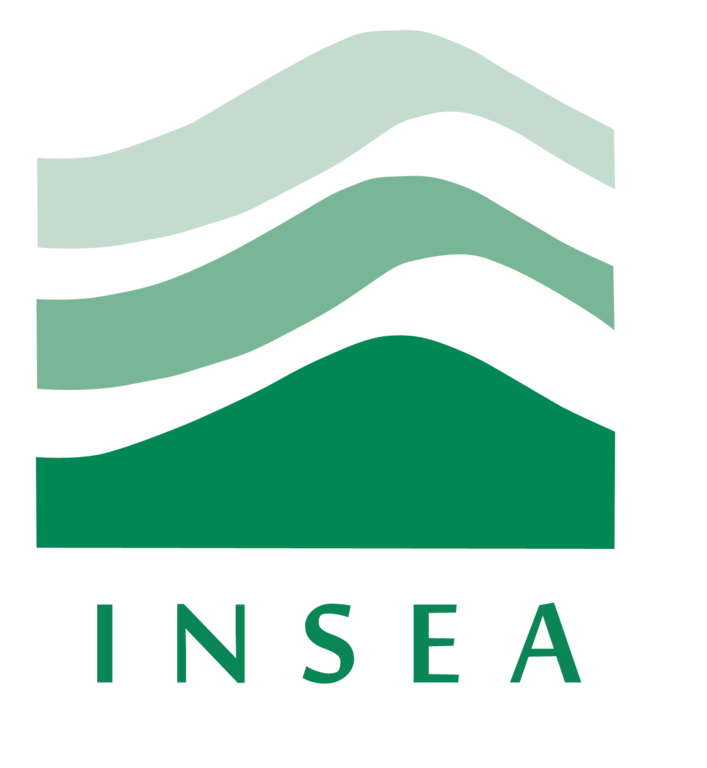}
	\end{minipage}\hfill
   
		\vspace*{-0.9cm}
\end{figure}

		\begin{center}
			\begin{minipage}[t]{1\textwidth}
				\begin{flushleft}
					\begin{spacing}{2}
						\begin{center}
							
						{\large \textbf{ CEDoc SIDD}\\}
						{\large \textbf{Research Laboratory : SI2M}}\\[-0.4cm]

							{\Huge\textbf{Doctoral thesis} \\}
							{\large \textbf {Presented to obtain the degree of doctor \\ in Applied Mathematics } }\\
					
						{Presented by:}\\[-0.4cm]
							{\Large \textbf{Mohammed Bazirha}\\}

						\end{center}
					\end{spacing}
				\end{flushleft}
			\end{minipage}\\[-1cm]

			\vspace{1cm}
			\rule{\linewidth}{0.5mm}
			\begin{center}
			\huge  	\textbf{Mathematical models and heuristics \\for the home health care routing \\ and scheduling problem} \\
			\end{center}
			\rule{\linewidth}{0.5mm}\\
			[0.7cm]
			
			 	{Supervised by:}\\[0.3cm]
			 
			 \begin{tabular}{ll}
			 	\textbf{Pr.}& \textbf{Abdeslam Kadrani}\\[0.2cm]
			 	
			 	\textbf{Pr.} & \textbf{Rachid  Benmansour}\\

			 \end{tabular}
		 
			\vspace{6mm}
			\rule{\linewidth}{0.1mm}\\
			\vspace{3mm}
		  Defended on  \textbf{December 29, 2022} at \textbf{INSEA-Rabat}\\[0.5cm]

			\begin{flushleft}
\textbf{Jury:}\\
[0.5cm]

				\begin{tabular}{lllll}

		     	Pr. Adil & Kabbaj         & PES at INSEA-Rabat & President\\
		     Pr. Abdeslam &Kadrani   &  PES at INSEA-Rabat & Thesis Supervisor\\
		     
		     
		     Pr. Rachid& Benmansour & PH\hspace{2mm} at INSEA-Rabat & Thesis Co-Supervisor\\
		     
		     Pr. Rachid& Ellaia & PES at EMI-Rabat & Reviewer\\
		     
		     Pr. Souad& El Bernoussi  & PES at Faculty of Sciences Rabat & Reviewer\\
		     
		     Pr. Youssef& Benadada  & PES at ENSIAS-Rabat & Reviewer\\
		    	
			\end{tabular}

			\end{flushleft}
		
		\rule{\linewidth}{0.1mm}\\
		\end{center}	
		
	\end{titlepage}

\begin{titlepage}
	\vspace*{-3cm}
	\begin{figure}[!htb]
		\begin{minipage}{0.2\textwidth}
			\centering
			\includegraphics[width=1.0\linewidth]{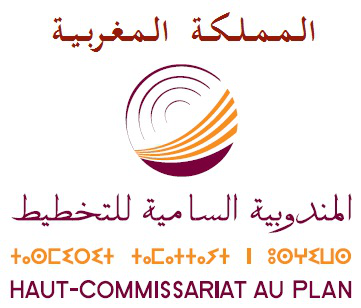}
		\end{minipage}\hfill
		\begin{minipage}{0.2\textwidth}
			\centering
			\includegraphics[width=1.0\linewidth]{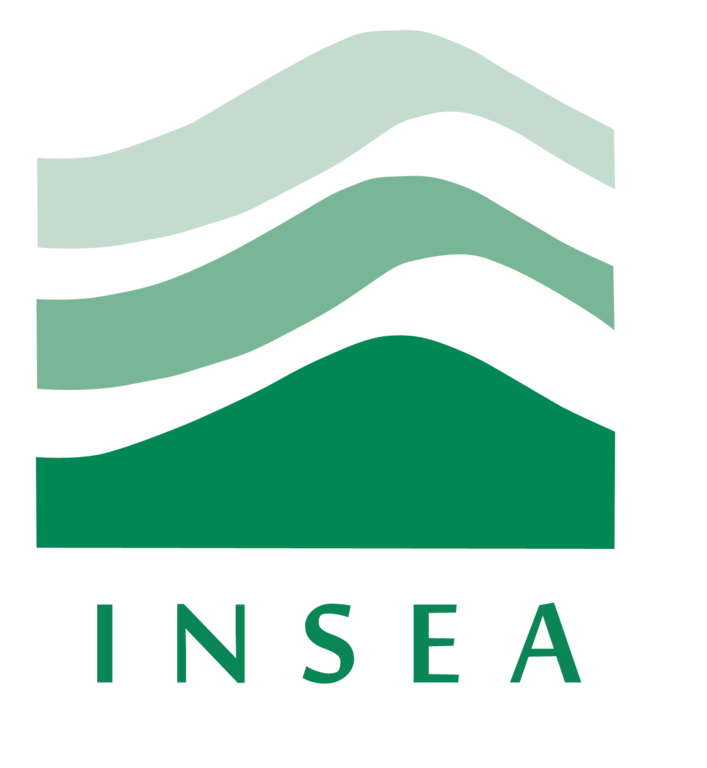}
		\end{minipage}\hfill
		
		\vspace*{-0.9cm}
	\end{figure}
	
	\begin{center}
		\begin{minipage}[t]{1\textwidth}
			\begin{flushleft}
				\begin{spacing}{1.8}
					\begin{center}
							\vspace*{0.9cm}
						
						{\large \textbf{Doctoral Studies Center "SCIENCE, ENGINEERING AND SUSTAINABLE DEVELOPMENT" (CEDoc SIDD)}\\}

						{\large \textbf{Research Laboratory : SI2M}\\}
						
						{\large \textbf {Doctoral thesis in Applied Mathematics  } \\}

						{\Large \textbf{Under the theme: }\\}

					\end{center}
				\end{spacing}
			\end{flushleft}
		\end{minipage}

		\rule{\linewidth}{0.5mm}
		\begin{center}
		\huge  	\textbf{Mathematical models and heuristics \\for the home health care routing \\ and scheduling problem} \\
		\end{center}
\rule{\linewidth}{0.5mm}\\
[0.5cm]

	 	{Prepared by:}\\[0.2cm]
      \textbf{Mohammed Bazirha}\\[0.5cm]
 
 	 	{Supervised by:}\\[0.2cm]
 
 \begin{tabular}{ll}
 	\textbf{Pr.}& \textbf{Abdeslam Kadrani}\\[0.1cm]
 	
 	\textbf{Pr.}& \textbf{Rachid  Benmansour}\\

 \end{tabular}
 
 \vspace{6mm}
 \rule{\linewidth}{0.1mm}\\
 \vspace{3mm}
 Defended on  \textbf{December 29, 2022} at \textbf{INSEA-Rabat}\\[0.5cm]

 \begin{flushleft}
 	\textbf{Jury:}\\
 	[0.3cm]

 	\begin{tabular}{lllll}

 		Pr. Adil & Kabbaj         & PES at INSEA-Rabat & President\\
 		Pr. Abdeslam &Kadrani   &  PES at INSEA-Rabat & Thesis Supervisor\\
 		
 		
 		Pr. Rachid& Benmansour & PH\hspace{2mm} at INSEA-Rabat & Thesis Co-Supervisor\\
 		
 		Pr. Rachid& Ellaia & PES at EMI-Rabat & Reviewer\\
 		
 		Pr. Souad& El Bernoussi  & PES at Faculty of Sciences Rabat & Reviewer\\
 		
 		Pr. Youssef& Benadada  & PES at ENSIAS-Rabat & Reviewer\\
 		
 	\end{tabular}

 \end{flushleft}

	\begin{center}
			\vspace{2mm}
		{ \large \textit{\textbf{Academic year: 2021-2022.}}}\\
	\end{center}
	\end{center}	
	 	\vspace{0.5mm}
	\rule{\linewidth}{0.1mm}\\
\end{titlepage}

 {\fontfamily{ptm}\selectfont
 	\setstretch{1.20}
 	
 	  \newpage
 
 \topskip0pt
 \vspace*{\fill}
 
 \begin{center}
 	\textbf{\Large Acknowledgment}
 \end{center}

 { \large
 
 	\vspace{10 mm}
 
 Words cannot express my gratitude to my thesis supervisor Pr.  \textbf{Abdeslam Kadrani} and my thesis co-supervisor  Pr. \textbf{Rachid Benmansour} for their supervision as well as for their support, their  relevant remarks and their encouragement.
 
 \vspace{10 mm}

 Also, I have the honor to express my sincere thanks to the President Pr. \textbf{Adil Kabbaj}.  I am also grateful to reviewers Pr. \textbf{Rachid Ellaia}, Pr. \textbf{Souad El Bernoussi} and Pr. \textbf{Youssef Benadada}  for agreeing to be members of the jury, for the time they devoted to reading the report and for their remarks and comments that helped to improve the quality of the report. 
 
 \vspace{10 mm}
 
 Lastly, I would be remiss in not mentioning my family, especially my \textbf{parents}. Their belief in me has kept my spirits and motivation high during this process. I would also like to thank  my \textbf{friends} and everyone who has encouraged me.
 
 } 
  \vspace*{\fill}
  
  \newpage

  \vspace*{\fill}
  
  \begin{center}
  	\textbf{\Large Dedication}
  \end{center}
  
  { \large
  	 \centering
  	 	\vspace{10 mm}
  
  To my parents
  
   \vspace{6 mm}
  No words can express my love, respect and gratitude for your tenderness and the sacrifices you have made for my education and well-being.
   
   \vspace{16 mm}
   
  To all my family members
  
   \vspace{6 mm}
  For the love and affection they have given me throughout my life, for their sacrifice, their trust and their patience.
  
   \vspace{16 mm}
   
  To anyone 
  
   \vspace{6mm}
  Who has taught me to appreciate the essence of life.
  
   \vspace{16 mm}
  I dedicate this  work.

  } 
  \vspace*{\fill}

	\clearpage

    \vspace*{\fill}
    
      \begin{center}
\textbf{\Large  {List of Abbreviations}}
  \end{center}
\vspace{10mm}
   
   {\small
  \begin{acronym}
  	{ 
  		 
  	 \acro {HHC} {Home Health Care}
  	 \acro{HHCRSP} {Home Health Care Routing and Scheduling Problem}
  	
  	\acro{GA} {Genetic Algorithm}

  	\acro{GVNS} {General Variable Neighborhood Search}
  	
  	\acro{NSGA-II} 	{Non-Dominated Sorting Genetic Algorithm II }
  	
  	\acro{MOEA/D} {Multi-Objective Evolutionary Algorithm based on Decomposition}

     \acro{VRPTW} {Vehicle Routing Problem with Time Windows }
     
     \acro {SPR} {Stochastic Programming with Recourse}
     
      \acro{VRP}  {Vehicle Routing Problem}

    	\acro{VNS} {Variable Neighborhood Search}
     
     \acro {TS} {Tabu Search}
     \acro {SA} {Simulated Annealing}
     
      \acro{MOEAs} {Multi-Objective Evolutionary Algorithms}
  }
  \end{acronym}
    
}

  \vspace*{\fill}
  
\begin{abstract}
 This thesis addressed the \ac{HHCRSP}, which is a class of workforce scheduling problems. The \ac{HHCRSP} is an extension of the \ac{VRPTW} to which the constraints related to the \ac{HHC} context are added. It aims to provide care services to patients at their homes instead of going to the hospital. 
 We dealt with three different problems from the optimization viewpoint. In the first one, we considered a deterministic model to tackle the \ac{HHCRSP} with multiple time windows, multiple services, their synchronization if they are required to be simultaneous and skill requirements. We proposed a new mathematical to solve this problem along with a \ac{GVNS} based heuristic to solve large instances.
 In the second problem, we extended the deterministic model to cope with uncertainties in terms travel and service times.  We proposed two \ac{SPR} models. In the first \ac{SPR} model, we defined the recourse as a penalty cost for the tardiness of services and a remuneration for caregivers’ overtime. In the second \ac{SPR} model, we defined the recourse as skipping patients if their time windows should be violated. We embedded Monte Carlo simulation, which is used to estimate the expected value of recourse, into a  \ac{GA} based heuristic to solve SPR models. 
 In the last problem, we kept the multi-objective aspect of the deterministic model without aggregating its objective functions, and we used algorithms based on Pareto dominance to find the non-dominated solutions and then involve the decision-maker to select which one he prefers.  Two approaches, Pareto and decomposition based, with multi-objective evolutionary algorithms are adopted to solve the \ac{HHCRSP}. Three algorithms are implemented:  \ac{NSGA-II}, \ac{MOEA/D} and a hybrid \ac{NSGA-II} with \ac{MOEA/D} (hybrid) algorithm.
 Computational results highlighted the efficiency of \ac{GVNS} to solve the deterministic model and the adequacy of the \ac{GA} to be used with the simulation to solve the \ac{SPR} model. For the multi-objective  \ac{HHCRSP}, computational results and performance measures inferred that the hybrid algorithm found solutions that better approximate the Pareto front while the \ac{MOEA/D} algorithm solved instances faster in terms of CPU running times. 
\end{abstract}

\vspace{1cm}

\begin{keywords}
Mathematical modeling; Stochastic optimization; Multi-objective optimization;  Home health care; Routing and scheduling; Multiple time windows; Synchronization; Meta-heuristics.
\end{keywords}
\newpage
\begin{arabtext}\\
	\setcode{utf8}
	\begin{center}
		\textbf{\large ملخص}
	\end{center}
	
	تتناول هذه الأطروحة مشكل تحديد المسارات والجدولة الزمنية للرعاية الصحية المنزلية  \LR{(HHCRSP)}، وهي فئة من مشاكل الجدولة الزمنية للقوى العاملة. يعد \LR{HHCRSP}  امتدادًا لمشكل تحديد مسارات السيارات مع النوافذ الزمنية \LR{(VRPTW)} والتي تُضاف إليها قيود تتعلق بسياق الرعاية المنزلية. يهدف \LR{HHCRSP}  إلى تقديم خدمات الرعاية للمرضى في منازلهم وتجنب التنقل إلى المستشفيات.
	تعاملنا مع ثلاث مشاكل مختلفة من منظور  البحث عن الحل الأمثل. في الأول، اعتبرنا نموذجًا حتميا لمعالجة \LR{HHCRSP} بنوافذ زمنية متعددة، المهارات المطلوبة، وخدمات متعددة، ومزامنتها إذا كان مطلوبًا أن تكون الخدمات في أن واحد. لقد اقترحنا نموذجًا رياضيًا جديدًا و خوارزمية إرشادية تعتمد على البحث العام متغير الجوار \LR{(GVNS)} لحل الحالات الكبيرة للمشكل.
	في المشكلة الثانية، قمنا بتوسيع النموذج الحتمي للتعامل مع عشوائية أوقات السفر والخدمة. اقترحنا نموذجين للبرمجة العشوائية مع الرجوع  \LR{(SPR)}. في نموذج \LR{SPR} الأول، حددنا الرجوع على أنه تكلفة جزائية لتأخر الخدمات ومكافأة عن العمل الإضافي لمقدمي الرعاية. في نموذج \LR{SPR} الثاني، قمنا بتعريف اللجوء على أنه إلغاء زيارات المرضى الذين انتهكت نوافذهم الزمنية. قمنا بتضمين محاكاة مونت كارلو، والتي تُستخدم لتقدير القيمة المتوقعة للرجوع، في  خوارزمية إرشادية تعتمد على الخوارزميات الجينية\LR{(GA)}  لحل النموذجين\LR{SPR}.   
	في المشكلة الأخيرة، احتفظنا بطبيعة النموذج الحتمي متعدد الأهداف دون تجميع دواله الهدف، واستخدمنا الخوارزميات القائمة على سيادة باريتو لإيجاد حلول غير مسودة وإشراك صانع القرار في اختيار الحل الذي يفضل. تم اعتماد نهجين، قائمين على باريتو وعلى التقسيم، مع خوارزميات تطورية متعددة الأهداف لحل \LR{HHCRSP}. تم تنفيذ ثلاث خوارزميات: التصنيف الجيني الغير المسود\LR{(NSGA-II)}، الخوارزمية التطورية متعددة الأهداف القائمة على التقسيم \LR{(MOEA/D)} وخوارزمية هجينة بين \LR{NSGA-II} و \LR{MOEA/D}.
	أبرزت النتائج الحسابية كفاءة \LR{GVNS} في حل النموذج الحتمي ومدى ملاءمة \LR{GA} لاستخدامها مع المحاكاة لحل نموذج \LR{SPR}. بالنسبة إلى \LR{HHCRSP} متعدد الاهداف، ابرزت النتائج الحسابية ومقاييس الأداء أن الخوارزمية الهجينة وجدت حلولًا تقارب بشكل أفضل حدود باريتو بينما قامت خوارزمية \LR{MOEA/D}بحل الحالات بشكل أسرع من حيث مدة تشغيل وحدة المعالجة المركزية.

\vspace{1cm}

	\textbf{الكلمات الرئيسية: }
	النمذجة الرياضية؛
	الأمثلة العشوائية؛
   الأمثلة متعددة الأهداف؛
	الرعاية الصحية المنزلية ؛
	تحديد المسارات والجدولة الزمنية  ؛
	نوافذ زمنية متعددة؛
	التزامن؛
	طرق إرشادية.

\end{arabtext}

\newpage

\begin{center}
\textbf{\large Résumé}
\end{center}

Cette thèse traite le problème de routage et de planification des soins de santé à domicile (HHCRSP), qui est une classe de problèmes de planification de la main-d'œuvre. Le HHCRSP est une extension du problème de tournées de véhicules avec fenêtres de temps (VRPTW) auquel sont ajoutées des contraintes liées au contexte des soins à domicile. Il vise à fournir des services de soins aux patients à leur domicile et à éviter les déplacements vers les hôpitaux.
	Nous avons traité trois problèmes différents du point de vue de l'optimisation. Dans le premier, nous avons considéré un modèle déterministe pour traiter le HHCRSP avec des fenêtres de temps multiples, des services multiples, leur synchronisation s'ils doivent être simultanés et des compétences requises. Nous avons proposé un nouveau modèle mathématique et une heuristique basée sur la recherche générale à voisinage variable (GVNS) pour résoudre les grandes instances. 
Dans le deuxième problème, nous avons étendu le modèle déterministe pour faire face aux incertitudes en termes de temps de déplacement et de service.  Nous avons proposé deux modèles de programmation stochastique avec recours (SPR). Dans le premier modèle SPR, nous avons défini le recours comme un coût de pénalité pour le retard des services et une rémunération pour les heures supplémentaires des soignants. Dans le second modèle SPR, nous avons défini le recours comme le fait de sauter des patients si leurs fenêtres de temps ne sont pas respectées. Nous avons intégré la simulation de Monte-Carlo, qui est utilisée pour estimer l’espérance du recours, dans l'heuristique basée sur l'algorithme génétique (AG) pour résoudre les modèles SPR. 
 	Dans le dernier problème, nous avons gardé la nature du modèle multi-objectif, sans agréger ses fonctions objectifs, et nous avons utilisé les algorithmes basés sur la dominance de Pareto pour trouver les solutions non dominées et impliquer par la suite le décideur pour choisir celle qu'il préfère.  Deux approches, basées sur Pareto et sur la décomposition, avec des algorithmes évolutionnaires multi objectifs sont adoptés pour résoudre le HHCRSP. Trois algorithmes sont mis en œuvre : l’algorithme génétique de tri non dominé II (NSGA-II), l'algorithme évolutionnaire multi objectifs basé sur la décomposition (MOEA/D) et un algorithme hybride NSGA-II avec MOEA/D (hybride).
 Les résultats de calcul ont mis en évidence l'efficacité du GVNS pour résoudre le modèle déterministe et l'adéquation de l'AG à être utilisé avec la simulation pour résoudre le modèle SPR. Pour le modèle multi objectifs, les résultats de calcul et les mesures de performance ont montré que l'algorithme hybride a trouvé des solutions qui s'approchent mieux de la frontière de Pareto, tandis que l'algorithme MOEA/D a résolu les instances plus rapidement en termes de temps d’exécution du CPU.

\vspace{1cm}

\textbf{\textit{Mots clés :}}
	Modélisation mathématique; Optimisation stochastique; Optimisation multi-objectifs; Soin à domicile; Routage et planification; Fenêtres de temps multiples; Synchronisation; Méta-heuristiques.

\newpage
\begin{center}
	\textbf{\large Résumé substantiel}
\end{center}

Le domaine des soins de santé à domicile (HHC) est un domaine  de recherche qui a pris ces dernières années un énorme essor et a fait l’objet d’une attention croissante en raison de l'augmentation de l'espérance de vie et du faible taux de natalité. Les HHC visent à fournir un large éventail de services aux patients à leur domicile, dans un environnement personnel \cite{cisse2017or}. Ces services permettront aux patients de vivre de manière autonome le plus longtemps possible, même en cas de blessure, de maladie, de vieillissement ou de maladie invalidante.  Ils couvrent des soins médicaux, comme les soins infirmiers, la physiothérapie et l'orthophonie. Les services de HHC peuvent même concerner des tâches dans d'autres activités de la vie quotidienne des personnes âgées telles que les aider à manger, à s'habiller et se laver et autres tâches ménagères.

Le HHC est l'une des alternatives aux traitements hospitaliers classiques qui ont vu le jour grâce aux progrès techniques, organisationnels et thérapeutiques. Il était initialement destiné à réduire l'engorgement des hôpitaux pour certains patients. Ce type de soins est en plein développement et est de plus en plus demandé par les patients puisqu'il leur permet d'être traités à domicile. Le HHC devrait générer au moins trois avantages qui sont : la diminution des admissions à l'hôpital, la diminution de la durée d'hospitalisation et la possibilité pour les patients de rester à leur domicile et de recevoir des soins et de l'assistance \cite{cisse2017or}.

 Dans le rapport mondial sur le vieillissement et la santé de 2015, il était mentionné que les gens devraient vivre en moyenne jusqu'à 77 ans, dont 15 seraient passés avec une forme de handicap. Le rapport mentionne également que dans l'un des plus grands hôpitaux français, 20 \% de tous les patients âgés de plus de 70 ans étaient significativement moins capables d'effectuer les tâches de base nécessaires à la vie quotidienne. En outre, la proportion de personnes âgées ne cesse d'augmenter dans les pays européens et devrait encore s'accroître dans les décennies à venir \cite{tarricone2008home}; de même, les membres de la famille vivent plus souvent séparés les uns des autres \cite{mankowska2014home}.

Ces dernières années, les questions liées à la santé ont attiré l'attention d'un grand nombre de chercheurs. Leurs études sont menées non seulement du point de vue des techniques sociales et médicales, mais aussi du point de vue de l'optimisation. L'optimisation a concerné un large éventail de caractéristiques et de contraintes des opérations de soins de santé à domicile. Ces caractéristiques concernent l'organisation des soins à domicile, les soignants et les patients. Le HHCRSP est considéré comme une version étendue du problème de tournés de véhicules avec fenêtres de temps (VRPTW) auquel sont ajoutées des contraintes liées au contexte des HHC. Outre la définition des itinéraires des soignants, il s'agit également de définir leur affectation aux patients et de planifier leurs visites pour fournir les services demandés.
Bien que les études précédentes aient traité un grand nombre de caractéristiques et de modèles, elles présentent encore certaines limites. La première limite dans les modèles déterministes qui sont bien étudiés par la communauté des chercheurs, est la non prise en compte des fenêtres des temps multiples. De plus, dans la plupart des études, les services multiples et leur synchronisation ne sont considérés que comme des services doubles par patient. De toute évidence, la combinaison de ces deux caractéristiques dans un même modèle n'a pas été envisagée jusqu'à présent. Par conséquent, les modèles précédents (eg. \cite{mankowska2014home,afifi2016heuristic,rasmussen2012home,frifita2020vns,lin2021matching}) deviennent impuissants de traiter le cas où les patients peuvent être disponibles dans différentes périodes et qu'ils demandent des services multiples.

La deuxième limite est l'incertitude des paramètres, tels que les temps de voyage et de service, qui sont rarement pris en compte dans le contexte des soins à domicile. Dans le monde réel, ces deux paramètres ne peuvent être déterministes et sont soumis à l'incertitude. Par conséquent, le calendrier prédéfini doit être révisé à chaque fois qu'un changement imprévu dans les situations pratiques est survenu. Autrement, si aucune action d'adaptation n'y est apportée, il est probable qu'il y aura des retards dans les services surtout pour les patients qui n’ont pas encore été visités. Cela qualifiera le service par de qualité médiocre, voire même un service à risque.  

Une autre reproche pour ces modèles concerne leurs approches de résolution, qui utilisent souvent les techniques d'agrégation pour traiter le problème d'optimisation multi-objectif. Le souci concerne essentiellement les pondérations des fonctions objectifs qui doivent être connues a priori. L'affectation des poids à des objectifs contradictoires est une tâche confuse, le décideur doit être impliqué a posteriori pour sélectionner une solution préférée parmi l'ensemble des solutions non dominées.

Cette thèse vise à proposer de nouveaux modèles et approaches pour surmonter les limitations susmentionnées dans le HHCRSP. Le premier objectif est de proposer un modèle déterministe général qui peut traiter simultanément un nombre arbitraire de services ainsi qu'un nombre arbitraire de fenêtres de temps pour chaque patient.  Le second objectif est lié à l'incertitude des paramètres. Il vise à proposer une solution qui prend en compte l'incertitude des temps de déplacement et de service pour établir une planification robuste. Le troisième objectif est de développer des algorithmes qui peuvent fournir l'ensemble des solutions non dominées, puis impliquer le décideur a posteriori pour sélectionner la solution qu'il préfère.

La première contribution traite une version déterministe du HHCRSP où des nouvelles contraintes sont considérées, telles que les fenêtres de temps multiples, les services multiples par patient, leur synchronisation s’ils doivent être fournis simultanément, et les compétences requises. Dans une première version, on étudie le problème uniquement avec des fenêtres de temps multiples considérées comme flexibles, et on montre par une comparaison les avantages de l'utilisation de fenêtres de temps multiples. Ensuite, dans la version étendue, on propose un nouveau modèle où chaque patient peut demander plusieurs services qui peuvent être simultanés, ce qui implique la synchronisation des visites des soignants. Compte tenu de la complexité du problème, qui est NP-difficile, nous proposons une heuristique basée sur la recherche générale à voisinage variable (GVNS) pour résoudre les grandes instances dans un temps de calcul raisonnable. 

\newpage

La deuxième contribution aborde l'incertitude des temps de déplacement et de services qui peuvent survenir à tout moment dans les tournées des soignants. Deux modèles de programmation stochastique avec recours (SPR) sont proposés pour établir une planification robuste qui prend en compte de l'incertitude des paramètres.  Dans le premier modèle, le recours est défini comme un coût de pénalité pour le retard des opérations de services et une rémunération pour les heures supplémentaires des soignants. Ce modèle tient compte de plusieurs services par patient et de leur synchronisation s'ils doivent être simultanés. Dans le deuxième modèle, le recours est défini comme le fait de sauter un patient si ses périodes de disponibilité ne sont pas respectées. Ce modèle suppose que les patients peuvent spécifier plusieurs fenêtres de temps, qui sont considérées comme fixes. La simulation de Monte Carlo est utilisée pour estimer la valeur de recours attendue pour les deux modèles, qui est intégrée dans l'algorithme génétique pour résoudre les deux modèles SPR.

La troisième contribution concerne le HHCRSP avec multi-objectifs résolu par les algorithmes conçus pour approcher la frontière de Pareto. On suppose que le décideur doit intervenir a posteriori pour sélectionner une solution parmi un ensemble de solutions non dominées. Trois objectifs sont optimisés simultanément. Il s'agit de minimiser les temps de déplacement et d’attente des soignants et d'équilibrer leur charge de travail.  Des approches basées sur Pareto et la décomposition avec des algorithmes évolutionnaires multi-objectifs sont adoptées pour résoudre le HHCRSP avec des fenêtres de temps multiples. Trois algorithmes sont mis en œuvre à savoir: l'algorithme génétique de tri non dominé II (NSGA-II), l'algorithme évolutionnaire multi-objectif basé sur la décomposition (MOEA/D) et un algorithme hybride.

\newpage

 \pagestyle{fancy}\fancyhf{}
\renewcommand{\chaptermark}[1]{%
	\markboth{\chaptername \ \thechapter.\ #1}{}}
\fancyhead[L]{\leftmark}
\fancyhead[R]{}
\rfoot{ \thepage}
\renewcommand{\footrulewidth}{0.4pt}

	\pagenumbering{arabic}
  {\fontfamily{ptm}\selectfont
 	\setstretch{1.2}
	\tableofcontents
	\newpage
 {\fontfamily{ptm}\selectfont
	\setstretch{1.20}
	
	\listoftables
	
	\listoffigures

\chapter*{Introduction} 
\markboth{Introduction}{}
\addcontentsline{toc}{chapter}{Introduction}
 
 	

\vspace{-10mm}

 The \ac{HHC} is an evolving research area that has received growing consideration in the recent years due to increasing life expectancy and low birth rates. It aims to provide a wide range of services to patients at their homes in a personal environment \cite{cisse2017or}. It will allow them to live autonomously for as long as possible even with an injury, illness, aging or disabling disease.  It covers a wide range of services that may involve medical care, such as nursing, physical therapy, and speech therapy. It may include helping elderly individuals with activities of daily living such as eating, dressing, and bathing. It can also involve assistance with cooking, cleaning, and other housekeeping. 
 
The \ac{HHC} is one of the alternatives to classical hospital treatments that have arisen thanks to technical, organizational, and therapeutic progress. It was initially intended to reduce hospital overcrowding for some patients. This type of care is in full development and is increasingly requested by patients since it allows them to be treated at home. The \ac{HHC} is expected to generate at least three advantages which are: a decrease in hospital admissions, a decrease in hospitalization duration and the ability for patients to remain in their homes and receive care and assistance \cite{cisse2017or}.

\vspace{2mm}
\textbf{ \underline{Background and motivation}}
\vspace{2mm}
 
In the world report on ageing and health of 2015, it was mentioned that people are expected to live on average to 77 years, of which 15 would be spent with some form of disability. The report also mentions that in one of France's largest hospitals, 20 percent of all patients over the age of 70 were significantly less able to perform the basic tasks necessary for daily living. In addition, the proportion of older people is increasing steadily in European countries and is predicted to rise still further in the coming decades \cite{tarricone2008home} as well as  family members live more often separated from each other \cite{mankowska2014home}.

\ac{HHC} issues have attracted the attention of a large number of researchers in recent years. Their studies are not only from the perspective of social and medical techniques, but also from the optimization viewpoint. The optimization has concerned a wide range of characteristics and constraints of \ac{HHC} operations. These features concern the \ac{HHC} organization, caregivers, and patients. The \ac{HHCRSP} is considered as an extended version of the \ac{VRPTW} to which constraints related to the \ac{HHC} context are added. In addition to defining caregivers’ routes, it also involves their assignment to patients and scheduling their visits to provide the requested services.

We choose to deal with the \ac{HHCRSP}  to contribute in improving health care systems and patients' daily living. Thus, allowing people requiring care  to be treated at home in order to avoid moving to hospitals and overcrowding them. The \ac{HHCRSP} is considered as  one of the constrained problems where limited resources must be optimized  to satisfy the maximum of constraints.  To achieve that, we aim to propose new models and heuristics to efficiently assign  the limited number of caregivers to patients by optimizing one or more criteria taking into account many real-life constraints.

\vspace{2mm}
\textbf{ \underline{Research context }}
\vspace{2mm}

Although previous studies have dealt with a large number of features and models, they still have some limitations. The first limitation concerns deterministic models, which are well studied by the research community. In these models, multiple time windows are not considered. According to \cite{cisse2017or}, it is interesting to  define all periods in which patients will be available to receive care services. Furthermore, multiple services and their synchronization in most studies considered them only as double services per patient.  In \cite{redjem2016operations}, the authors have conducted a survey of many French home care structures, which showed that patients need several care services per day. Therefore, the synchronization of services should be considered if a patient needs more than a caregiver to provide one or multiple simultaneous services such as  dressing, bathing, getting out of bed and dosing medicine \cite{rasmussen2012home}. Obviously, combining the two features in the same model has not been considered so far. These previous models are powerless to deal with cases where patients can be available in many periods as well as if he requested multiple services.
	The second limitation is the uncertainty of parameters, such as travel and service times, which are less considered in the \ac{HHC} context. In the real world, these two parameters cannot be deterministic and are exposed to uncertainty. Therefore, the predefined schedule should be revised whenever an unforeseen change in practical situations has occurred. Otherwise, if no adaptation action is taken, there will be delays in services, especially for patients who have not yet been visited. This could qualify the service as a poor quality, or even risky service. 
The last limitation for these models is related to their resolution approaches, which often use aggregation techniques in most studies to deal with the multi-objective optimization problem. The issue mainly concerns the weights of the objective functions that should be known a priori. Assigning weights to conflicting objectives is a difficult task, the decision maker needs to be involved a posteriori to select a preferred solution from the non-dominated solutions set.

\vspace{2mm}
\textbf{ \underline{Objectives, research questions and contributions }}
\vspace{2mm}

This thesis aims to propose new models and methods  to overcome the aforementioned limitations in the \ac{HHCRSP}. The first objective is to propose a general deterministic model that can simultaneously deal with  an arbitrary number of services as well as an arbitrary number of time windows for each patient.  The second one is related to the uncertainty of parameters. It aims to propose a solution that takes into account the uncertainty of travel and service times to establish a robust planning. The third objective is to develop algorithms than can provide the non-dominated solutions set, then involve the decision maker a posteriori  to select which solution he prefers.

In the first contribution, we dealt with a deterministic model where we considered new complicated constraints. In \cite{bazirha2019daily}, we developed a new mathematical model to deal with multiple soft time windows with a maximum of earliness ($E_{max}$) and a maximum  tardiness ($T_{max}$) of services. Patients could define all periods in which they are available to receive care services. A soft/flexible time windows would increase the chance of finding a feasible schedule as delays are accepted. However, any delay that occurs, in the case of hard/fixed time windows, the schedule becomes impractical and the company must resort to more caregivers (i.e. an additional cost) to keep a reasonable quality of service.  The soft/flexible time windows proposed with  $E_{max}$ and $T_{max}$   is a general case of hard/fixed time windows ($E_{max}=0$ , $T_{max}=0$)  and of  soft/flexible time windows($E_{max}=\infty$, $T_{max}=\infty$).  We developed a  \ac{GVNS} based heuristic to reduce the computational time needed to solve the model.  Computational results show that instances with multiple time windows are better optimized since the decision maker has many possibilities to schedule the visits. 
 Since patients need several care activities per day \cite{redjem2016operations}, we extended the model to allow them to request many services with a possible synchronization of several ones.  A strategy to select a time window for each patient and to ensure  synchronization of simultaneous services, is proposed and combined with the \ac{GVNS}-based heuristic. This strategy is mandatory since, for each patient, the selection of a time window could not be done independently of his availability periods. In fact, since a single period should be chosen to provide multiple services, this choice may be suitable for one caregiver but not for another.  \ac{GVNS} generates solutions (caregivers' routes and their assignment) while  this strategy  selects, for each patient, a time window and ensures the synchronization of simultaneous services.

In the second contribution \cite{bazirha2020scheduling,bazirha2021stochastic,bazirha2022optimization}, we dealt with two stochastic models where we considered the uncertainty of travel and service times parameters. We proposed two stochastic programming models with recourse to deal with the uncertainty. We estimated the expected value of recourse using Monte Carlo simulation since computing the expected real value by an explicit mathematical formula is very complex.  We proposed a new stopping condition for the simulation, instead of just running it for a maximum number of iterations (100 iterations as in \cite{shi2018modeling}), to find a good estimation of the expected real value.  We embedded it into  \ac{GA} and \ac{GVNS} based heuristics to solve the models with stochastic parameters. 
	In the first \ac{SPR} model \cite{bazirha2021stochastic}, we  dealt with the \ac{HHCRSP} with multiple synchronized services. The objective is to minimize the transportation cost and the expected value of recourse. The recourse is defined as a penalty cost for patients’ delayed services and a remuneration for caregivers’ overtime. We show through numerical tests that \ac{GVNS} is not suitable to be combined with the simulation to solve the \ac{SPR} model and the adequacy of the \ac{GA} since its parameters do not depend the problem parameters.  
In the second \ac{SPR} model \cite{bazirha2022optimization}, we assumed that time windows are hard/fixed, which must be respected without any earliness or tardiness. To ensure that the requested services are provided within patients' time windows, the recourse is defined as skipping patients when their time windows will be violated. To increase the chance of providing the maximum of   services, patients can specify multiple time windows. The objective is to minimize the transportation cost and the expected value of recourse, which expressed as the average number of unvisited patients.

Dealing with multi-objectives using aggregation techniques imply involving the decision maker a priori to assign weights, which is a difficult task.  To overcome this issue, the last contribution is related to the multi-objectives case where the decision maker is assumed to be involved a posteriori to select the solution he prefers from then non-dominated solutions set. We use Pareto and decomposition bases approaches  with multi-objective evolutionary algorithms to simultaneously optimize three objectives trying to approximate the Pareto set. 
Minimizing caregivers' total traveling times ($f_1$) is more important than minimizing their transportation cost since the traveling time depends on road conditions, which will be higher if such conditions are poor. Waiting time ($f_2$) is considered to be unproductive \cite{redjem2016operations} and must be minimized. It is assumed that caregivers are paid for their regular working time regardless of the amount of care they do \cite{braekers2016bi}. To ensure fairness among caregivers,   their workload must be balanced ($f_3$). 
Three algorithms are implemented:  \ac{NSGA-II} , \ac{MOEA/D} and a hybrid \ac{NSGA-II} with \ac{MOEA/D} (hybrid) algorithm.

\vspace{2mm}

\textbf{ \underline{Road map}}
\vspace{2mm}

This manuscript is divided into four chapters. \textbf{Chapter \ref{ChapterLR}} gives a literature review of the studies dealing with the routing and scheduling in home health care. We first focus on the travelling salesman problem and then the \ac{VRP} and its variants since it is  extended by the \ac{HHCRSP} to cover the additional constraints related to health care. A review of the \ac{HHCRSP} is done with a focus on objectives, constraints and resolution methods considered in the literature. The constraints are notably divided into three groups according to whether they concern the \ac{HHC} organization, patients, or caregivers. Finally, we show some limitations in models and methods proposed by researchers and tackle it. 

\textbf{Chapter \ref{ChapterDM}} deals with a deterministic version of the \ac{HHCRSP} where new constraints are considered,  such as multiple time windows, multiple services per patients, and the synchronization. In \textbf{section \ref{SectionDMTW}}, we study the problem only with multiple time widows, which are considered as soft/flexible, and we show through a comparison the advantages of using multiple time windows. Then  in \textbf{section \ref{SectionDMSMTW}}, an extended version of the model is proposed where each patient can request many services that can be simultaneous, which involve the synchronization of caregivers’ visits. Given the complexity of the problem, which is NP-hard, we propose a  \ac{GVNS}  based heuristic to solve large instances in a reasonable computation time. 

\textbf{Chapter \ref{ChapterSM}} addresses the uncertainty of travel and services times that may arise anytime in caregivers’ routes. Two \ac{SPR} models are proposed to establish a robust planning that takes into consideration the uncertainty of parameters.  In \textbf{section \ref{SectionSMS}}, the recourse is defined as a penalty cost for the tardiness of services operations and  a remuneration for caregivers’ overtime. This model takes into consideration multiple services per patient and their synchronization if they are required to be simultaneous. In \textbf{section \ref{SectionSMTW}}, the recourse is defined as skipping a patient if his availability periods will not be respected. This model assumes that patients can specify multiple time windows, which are considered hard/fixed. Monte Carlo simulation is used to estimate the expected recourse value for both models, which is embedded into the \ac{GA} to solve the two \ac{SPR} models.

\textbf{Chapter  \ref{ChapterMOMTW}} deals with the multi-objective \ac{HHCRSP} using algorithms designed to approach the Pareto front. It is assumed that the decision maker is to be involved a posteriori to select a solution from not dominated solutions set. Three objectives are simultaneously optimized which are: minimizing caregivers’ travel and waiting times as well as balancing their workload.  Approaches based on Pareto and decomposition with multi-objective evolutionary algorithms are adopted to solve the \ac{HHCRSP} with  multiple time windows. Three algorithms are implemented: \ac{NSGA-II}, \ac{MOEA/D} and a hybrid \ac{NSGA-II} with \ac{MOEA/D} (hybrid) algorithm. 


This thesis is concluded by a general conclusion  in which we summarize the work accomplished and the results obtained. Research perspectives are also proposed for the pursuit of our research works.

\vspace{2mm}
\textbf{ \underline{Publication and research activities}}
\vspace{2mm}

The following is a summary of the papers produced and published in the course of this thesis, listed in chronological order.


\begin{itemize}
	
	\item M. Bazirha, A. Kadrani, and R. Benmansour. Optimization of the stochastic home health
	care routing and scheduling problem with multiple hard time windows.\textbf{\textit{ International
	Journal of Supply and Operations Management}}, 9(2):235–250, 2022.  
		
	\item M. Bazirha, A. Kadrani, and R. Benmansour. Stochastic home health care routing and
	scheduling problem with multiple synchronized services.\textbf{\textit{ Annals of Operations Research}}, pages 1–29, 2021.

	\item M. Bazirha, A. Kadrani, and R. Benmansour. Scheduling optimization of the home health
	care problem with stochastic travel and care times. I\textbf{\textit{n 2020 5th International Conference
	on Logistics Operations Management (GOL)}}, pages 1–8. IEEE, 2020.
	
		\item M. Bazirha, A. Kadrani, and R. Benmansour. Daily scheduling and routing of home
	health care with multiple availability periods of patients. \textbf{\textit{In International Conference on
			Variable Neighborhood Search}}, pages 178–193. Springer, 2020.
		
	\item  Bazirha, M., Kadrani, A., Benmansour, R.. Home Health Care Routing and Scheduling Problem with Multiple Time Windows and Multiple Synchronized Services. (\textbf{Preprint: submitted})

	\item  Bazirha, M., Kadrani, A., Benmansour, R.. Pareto and Decomposition Based Approaches for the Multi-Objective Home Health Care Routing and Scheduling Problem. (\textbf{Preprint: submitted})
	
\end{itemize}

\chapter{ Literature review}\label{ChapterLR}

\section{Introduction}
Workforce scheduling problems consist of defining the schedules of hourly workers to meet the demands in a workplace such as a restaurant, a hospital, or a retail store. It is also a matter of defining schedules so that the transition from one shift to another is smooth and the work remains continuous. In many scenarios in which workers must carry out tasks at different locations they require some form of transportation. 

The workforce scheduling problem extends the travel salesman problem and \ac{VRP} to include also scheduling of workforce in the case that workers must travel to carry out tasks at different locations.  In \cite{castillo2016workforce}, the authors presented four types of workforce scheduling problems that require some form of workers’ transportation, which are: security personnel routing and rostering, scheduling technicians, manpower allocation and home health care. The  \ac{HHCRSP} is a type of workforce scheduling problem that is well-known and studied in the last years due to its importance for people in need of daily living services, our contributions are related to this problem.

\section{ Home health care routing and scheduling problem}

The \ac{HHCRSP} is class of workforce scheduling and routing problems (WSRPs) that extend the \ac{VRP} to also include the scheduling of caregivers' visits to provide patients' requested services at different locations. Classical models and methods studied in the \ac{VRP} are not suitable to be used in general for the WRSPs and in particular for the \ac{HHCRSP} because it contains additional characteristics related to workers such as skills requirements, dependencies among tasks, working overtime and workload balancing... etc.

In the following we will present characteristics and constraints considered in the literature for the \ac{HHCRSP} and then show some limitations related to the studies to tackle it.  We start with related problems to the \ac{HHCRSP}, then the most criteria considered in the literature,  and then the constraints and features considered by researchers. The constraints are notably divided into three groups according to whether they concern the \ac{HHC} organization, patients, or caregivers.

 \subsection{Related problems}
 Before going into greater detail about our contributions in the \ac{HHC} context, it is important to situate the \ac{HHC} problem in the global context of the vehicle routing problem.  The latter problem is extended by the \ac{HHCRSP} to cover the additional constraints related to health care context.
 \subsubsection{Travel salesman problem}

 The traveling salesman problem (also called TSP) is a combinatorial optimization problem, studied since the 19$^{th}$  century, is one of the most well-known in the field of operations research. It can be described as follows: Given a set of cities and the distance between every pair of cities, the problem is the challenge of finding the shortest and the most efficient route (Hamiltonian Cycle) for a person to take given a list of specific destinations (cities). 
 
 The TSP was mathematically formulated in the 1800s by the Irish mathematician W.R. Hamilton and by the British mathematician Thomas Kirkman. In 1954, Dantzig et al. \cite{dantzig1954solution} proposed a Mixed-integer linear program to solve the TSP. Garey et al. \cite{garey1979guide} have proved that the TSP is NP-hard. The TSP is divided into two categories. If the distance from one node to another is different from the inverse distance, the TSP is called asymmetric (ATSP). In symmetric TSP (STSP), this distance is the same. Sub-tours elimination constraints   are the hardest to satisfy since their number increases exponentially and equals $2^{n-1}-2$. Miller et al. \cite{miller1960integer}  have proposed a new formulation to reduce the complexity of sub-tours constraints from $O(2^n)$ to $O(n^2)$. 
 

 
 
 
 




 \subsubsection{Vehicle routing problem}

 The  \ac{VRP} is an extension of the TSP, where it deals with many vehicles. The problem is to define a tour for each vehicle to visit assigned customers. Each one has a start and an end location, which typically is the depot.   
 
 According to Braekers et al. \cite{braekers2016vehicle}, Dantzig et al. \cite{dantzig1954solution} were the first to introduce the \ac{VRP} called "Truck Dispatching Problem". They proposed a formulation of the problem to optimize the total traveled distance that involves a fleet of homogeneous vehicle to serve the demand for oil of a number of gas stations from a central hub. A few years later,  Clarke et al. \cite{clarke1964scheduling}
 generalized this problem to a linear optimization problem commonly encountered in the field of logistics and transportation. It is defined as follows: a set of customers geographically scattered around a central depot must be served by a fleet of vehicle with varying capacities. 
 
 The classical \ac{VRP}, also known as the capacitated \ac{VRP} (CVRP), it has been extended in many ways to cover additional real-life constraints or features, which has led to a large number of variants of \ac{VRP}. In the case of vehicles with varying capacities, the extended version is described in the literature as  the heterogeneous fleet \ac{VRP} (HFVRP), also known as the mixed fleet \ac{VRP}. Each vehicle has a different capacity, the problem is how to assign customers to vehicles in such way that the constraint of vehicle capacity is respected along its route.
 
 Another well-studied variant is the \ac{VRP} with time windows (\ac{VRPTW}). This variant assumes that deliveries to a given customer must occur within a certain time interval, which varies from customer to another. These time windows could be considered as hard/fixed or soft/flexible. In the first case, deliveries must respect the time interval, which the customer has specified without any earliness or tardiness. However, in the second case, time windows could not be respected, and delays can be accepted with a penalty cost. Others studies such as \cite{belhaiza2014hybrid,favaretto2007ant} have considered Multiple Time Windows for the \ac{VRP}, customers are supposed to be available in many periods, which they prefer.

 The \ac{VRP} with deliveries and pickups (VRPDP) \cite{savelsbergh1995general} extends the \ac{VRP} in such way to allow the transportation of goods from a depot to customers as well as from customers to the depot. It is defined in \cite{savelsbergh1995general}  as follows: A fleet of vehicles is available with a given capacities, a start location, and an end location. Each request is defined by the size of load to be transported, the origins and the destinations. The origins denote the locations where goods will be picked up and the destinations denote where the goods will be delivered. Each load has to be transported by one vehicle from its origin location to its destination location without any transshipment to others locations. VRPDP is an important logistics problem that has many applications, mainly in reverse logistics. An increasingly environment-friendly population implies more collection of recyclable goods. Also, companies are increasingly aware of the savings that can made by combining deliveries and collections.
 
 The VRPDP has also many variants such as: The \ac{VRP} with backhauling (VRPB), also known as the linehaul-backhaul problem. Linehaul (delivery) points are locations that will receive a quantity of goods. Backhaul (pickup) points are locations that send a quantity of goods. The critical assumption is that all deliveries must be made on each route before any pickups can be made. The \ac{VRP} with mixed deliveries and pickups (VRPMDP) does not make the above assumption and allows deliveries and pickups to occur in any order on a vehicle route.  In the \ac{VRP} with simultaneous deliveries and pickups (VRPSDP), the customers simultaneously have delivery and pickup demand and the vehicles makes a single stop to both deliver and pick up goods.

 The multiple depot \ac{VRP} (MDVRP) is an extension of the classic \ac{VRP}. In its original version, \ac{VRP} considers only one depot from which a fleet of vehicles departs and ends its tour. In the MDVRP, the fleet of vehicles is assigned to multiple depot locations. Each vehicle comes from a depot to serve assigned customers will return to the same depot.

 The periodic \ac{VRP} (PVRP) is a variation of the classic \ac{VRP} in which delivery routes are constructed for over a period of time (for example, multiple days). The PVRP was introduced by Beltrami et al. \cite{beltrami1974networks} and has evolved into a large body of works, with many interesting variants and applications emerging in recent years. The PVRP occur in a wide range of applications, such as elevator maintenance and repair, courier services, vending machine replenishment and the collection of waste \cite{francis2008period}.
 
 Despite a vast literature concerning the \ac{VRP} and its variants, it is however very difficult, if not impossible, to apply studied cited above on workforce scheduling problems in general and particularly for the home health care routing and scheduling problem. Indeed, workforce scheduling problems  is a specific case of application bringing new constraints, which requires an adapted modeling and methods of resolution.
 
 In the following, we will present the scientific literature specific to workforce scheduling problem that involve the transportation of workers in the \ac{HHC} context addressed in this work.
 
 \subsection{Criteria considered in the  home health care models}

Finding the optimal program with respect to a set of constraints implies the optimization of one or more criteria. Otherwise, the problem is reduced to constraint satisfaction, which seeks a feasible solution according to a set of constraints to be satisfied.
According to \cite{cisse2017or}, four main types of criteria are the most objective functions tackled in the literature, which are: minimize caregivers’ transportation cost, minimize the number of unassigned services, minimize the number of caregivers and maximize satisfaction. 

Approximately 90\% of papers \cite{cisse2017or} consider the transportation cost as a criteria to be optimized, which often is expressed as the sum of distances traveled or the travel time incurred by caregivers. Minimize the number of unassigned services arises when the \ac{HHC} company may not have sufficient caregivers to provide all services requested by patients. To overcome this issue, the company can subcontract to another \ac{HHC} company these unassigned services with a cost, which is always higher than the cost incurred by the provider if it is performed internally. Some studies such as \cite{hewitt2016planning,allaoua2013matheuristic,ikegami2007bounds} considered minimizing the number of caregivers as a criteria. The \ac{HHC} company must be flexible to adjust the number of caregivers on the routes each day \cite{cisse2017or}. 

Maximize the satisfaction involve both caregivers and patients. For patients, this goal can be achieved, for example, by respecting their availability periods, which could be as hard/fixed time windows \cite{akjiratikarl2007pso,bard2014sequential,liu2017mathematical,redjem2016operations}, or soft/flexible \cite{hiermann2015metaheuristics,mankowska2014home,nickel2012mid,trautsamwieser2011optimization}. In the first case, providing services to patients must be done without any earliness or tardiness. In the second case, a penalty cost is considered for any earliness or tardiness. Another example of considering patients’ satisfaction is related to their preferences regarding caregivers. A patient may prefer a caregiver over others, especially in the case of continuity of care where patient may be visited several times. Therefore, he can prefer visited by the same caregiver. Gender and language are some examples of preferences.

Caregivers' satisfaction is modeled as balancing their workload, which is the main criterion considered in the literature \cite{cisse2017or}. The decision maker has to ensure that caregivers have worked in a balanced way by avoiding gap between them in terms of workload. This goal could be achieved by balancing the number of services assigned to each caregiver or balancing the working time in such way that amount of time worked by each caregiver must be close to each other. Balancing the number of services assigned to each caregiver would only be effective if all services have very close processing times. The decision maker also should consider the unavailability of caregivers, especially in the case of a medium- or long-planning horizon \cite{nickel2012mid,wirnitzer2016patient,bard2013weekly}. Workday preferences, vacations, days off or sick leaves are some examples used in the literature that motivate using the unavailability of caregivers constraint in the studies \cite{cisse2017or}.

When many criteria are considered and must be simultaneously optimized, aggregation techniques are the most used in the literature. In \cite{duque2015home}, a lexicographical order is used to solve a bi-objective problem that aims to maximize the service level and to minimize caregivers’ traveling distance. This approach is based on a priori knowledge of some decision maker preferences to establish the lexicographic order, which is not a simple task. Few studied adopted Pareto approach to deal with the multi-criteria such as \cite{fathollahi2020bi,decerle2019memetic,braekers2016bi}. As far as we know approaches based on indicator and decomposition are not used  so far in the  \ac{HHC} context.

 \subsection{Constraints considered in the  home health care models}

HHC issues have attracted the attention of a large number of researchers in recent years. Their studies are not only from the perspective of social and medical techniques, but also from the optimization viewpoint \cite{shi2018modeling}. A significant diversity among features are considered in the existing \ac{HHCRSP} models to deal with the characteristics and constraints of \ac{HHC} operations. According to \cite{cisse2017or}, these features can be divided into three groups, depending on whether they are related to the \ac{HHC} service organization, patients, or caregivers. 

\subsubsection{Constraints related to the \ac{HHC} organization}

The constraints related to the \ac{HHC} organization could be divided into three groups according to \cite{cisse2017or}, which are: temporal, assignment, and geographic constraints.

In temporal constraints, two characteristics are to be considered when the decision maker defines caregivers' schedule. The first one is the planning horizon that could be considered as short, medium, or long term. The length of the planning horizon considered in the \ac{HHCRSP} models has been often one day \cite{mankowska2014home,rasmussen2012home,akjiratikarl2007pso,trautsamwieser2011optimization} or one week \cite{bard2014sequential,bard2013weekly}. In \cite{wirnitzer2016patient}, a planning horizon of four weeks
has been considered. The second one is the frequency of routing decision, which refers to how many times the decision maker repeat the routing decision within the planning horizon. The long-term horizon is affected by new patients entering the system continuously and therefore the conditions of previously admitted patients and the availability of caregivers change. Therefore, either the routing decisions can be changed at fixed time intervals or the update can be done when certain conditions are met \cite{cisse2017or}.

In assignment constraints, which are largely related to the continuity of care, also called patient–nurse loyalty in \cite{nickel2012mid} or employee regularity \cite{gamst2012branch}, patients must be visited by the same caregiver and during approximately the same time. Care services can be performed for patients under full, partial or no continuity of care, which avoid losses of information among caregivers \cite{cisse2017or} and builds a relationship of confidence with patients \cite{freeman2010continuity}. In the first case, a patient assigned to one and only one caregiver who is responsible for his/her care. The second case usually occurs when a patient requires more than one type of care. A reference caregiver can be affected, as in the case of full continuity of care, to coordinate providing the requested services by other caregivers. In the last case, the decision maker does not need to consider the previous information related to the caregiver–patient assignment.

 The dispersion of patients in a geographic area may let the decision maker to cluster the existing teams into districts based on relevant criteria such as territory, caregiver skills, and patient needs. Such a district-based approach allows to reduce travel times in caregivers’ routes within an assigned district and lets to form smaller care teams that can be easily managed. In the single-district case, the decision maker does not split caregivers into smaller clusters, which is the simplest way of managing teams and is considered in \cite{mankowska2014home}. However, in the case of multiple districts \cite{duque2015home,eveborn2006laps,lanzarone2012operations}, a caregiver is assigned to a geographical area or a subset of areas and may only serve patients within districts in which he is assigned.

Other types of \ac{HHC} services rather than classical ones have been considered in the literature, which are related to geographic constraints. These services include collecting biological samples and/or delivery of drugs or materials. Problems including these type of services such as \cite{kergosien2014routing,braysy2009optimization,liu2013heuristic} are similar to \ac{VRP} with pickup and delivery.  Locations where pickup and delivery services are performed must be included int the parameters considered in the \ac{HHCRSP} models, where routes are optimized considering these locations. The \ac{HHC} center, medical laboratories, pharmacies and chemotherapy/radiotherapy centers are some examples of these locations \cite{cisse2017or} that deliver these services.

\subsubsection{Constraints related to patients}

HHC companies deal with a wide range  of service types, each one define  the nature and frequency of its visit. In the literature, regarding constraints on the frequency of visits, some authors assume that a patient needs to be visited once a day \cite{hiermann2015metaheuristics} or once a week \cite{bard2014sequential}. However, other studies consider a patient can be visited several times a day \cite{rasmussen2012home} or a week \cite{nickel2012mid}. These later cases arise when patients usually require multiple services.   The services may have a temporal dependency (e.g. \cite{redjem2016operations,mankowska2014home,rasmussen2012home}) if they must be provided on the same day.  Some studies consider synchronized services, which arise when a patient need more than a caregiver to perform his requested services. For instance, overweight patient requires two caregivers to lift him/her. Other studies consider also disjunction among services, which must not performed simultaneously. For example, physical therapy cannot be delivered while a blood sample is collected.

Most studies consider a single availability period for each patient in which he will be available to receive caregiver(s). As described in the \ac{VRPTW}, studies in \ac{HHCRSP} also consider time windows as hard/fixed time \cite{akjiratikarl2007pso,bard2014sequential,liu2017mathematical,redjem2016operations} or soft/flexible \cite{hiermann2015metaheuristics,mankowska2014home,nickel2012mid,trautsamwieser2011optimization}. In the first case, patients’ availability periods must be respected without any violation while in the second case, these time windows can be violated with a penalty cost.  In \cite{bertels2006hybrid}, the authors assign to each patient two time windows, the first one is soft/flexible that is included the second one, which is hard/fixed. Some studies such as \cite{wirnitzer2016patient,trautsamwieser2011securing} consider patients' preferences where a caregiver cannot be assigned to a patient for gender incompatibility or personal reasons.

In terms of geographic constraints, most studies consider Euclidean distances between patients’ locations for experiments \cite{cisse2017or}. Some studies compute travel times based on a geographic information system \cite{begur1997integrated}.  Studies such as \cite{hiermann2015metaheuristics,rendl2012hybrid} have considered a multi-modal transportation network, where caregivers can take three main modes of transportation between patients’ homes, which are: car, walking and public transportation. Travel times between patients can change significantly according to the time of day, e.g., rush hour or off-peak times. Therefore, Rest and Hirsch \cite{rest2016daily} have considered a model that takes into account time-dependent travel times.

\subsubsection{Constraints related to caregivers}

 Most \ac{HHCRSP} studies consider a predefined working time for caregivers during which they will visit patients to provide their requested care services in their homes. Caregivers’ availability is either considered  as hard/fixed \cite{redjem2016operations,liu2017mathematical,braekers2016bi} or soft/flexible  \cite{trautsamwieser2011securing}. In the first case, caregivers’ overtime is not allowed and exceeding the predefined working time of a caregiver is considered as an infeasible solution. However, in the soft/flexible case, caregivers’ extra working times is allowed with a penalty cost in the objective function.

Breaks are also considered in the literature for caregivers, for instance, to take lunch and some rest. They should be programmed into the caregivers' routes by defining the start and end of the break as well as its duration. The break can be considered by forcing the caregiver to take a break either within a specific time interval \cite{trautsamwieser2011optimization,liu2017mathematical}, or when a certain route length is reached \cite{thomsen2006optimization}.

Caregivers’ assignment must respect skill (also called qualification) requirements \cite{mankowska2014home,braekers2016bi,lin2021matching,hiermann2015metaheuristics,duque2015home} while providing a specific service.  Three formulations are used in the literature: The first one is to assign a single qualification for each caregiver.  The second approach is to assign several skills for each caregiver, which is more general. And the third case that is based on level of qualification, which represents the minimum level for a caregiver to provide a requested service.

Caregivers are paid for their regular working time regardless of the amount of care they do. To ensure equity among them, the decision maker may consider caregivers’ workload balancing while assigning them to patients. Ensure a perfect balance of the workload among caregivers is difficult if it not impossible task, thus this constraint is considered as a soft and  generally embedded in the objective function.

Like the MDVRP, the \ac{HHCRSP} studies also consider multiple \ac{HHC} centers, but most of them assume a single center from which caregivers begin and end their tours. In \cite{trautsamwieser2011securing}, the authors proposed three types of starting and ending locations for caregivers which are: 1) caregivers start and finish at the \ac{HHC} center; 2) caregivers start and end their tours at home; 3) and  caregivers start and end their tours at home, but working time do not include  travel times between caregivers’ homes and patients. In \cite{bard2014sequential},  the authors dissociated the starting and ending depots for each caregiver, which is more generic and cover both single and multiple depots.

\subsection{Resolution methods}

In the following, we will present the different solution methods used in the literature to solve the HHCRSP.	

\subsubsection{Exact methods}

 Exact methods guarantee to reach the optimal solution for a given problem. However, the solving algorithms have an exponential complexity, which increases with the problem size. Therefore, they can only solve instances of limited size, depending on the nature of the problem. In addition, they are used to benchmark approximate methods and analyze their behavior on optimal solutions. In the literature, studies that used exacts methods are:  branch-and-price \cite{rasmussen2012home}, column generation \cite{yuan2015branch}  and branch-and-cut-and-price \cite{trautsamwieser2014branch}. 
 
 \subsubsection{Approximate methods}
 
Since exacts methods are powerless to reach the optimal solution or even a feasible one for large instances, approximate methods have been used widely by the community of researchers to overcome exact methods limitation.  The approximate methods are more flexible as they allow to control the stopping criteria: maximum number of iterations, maximum CPU running time or maximum number of no-improvement over the best solution. However, they do not guarantee to reach the optimal solution, but at least they are able to find a feasible one in a reasonable computational time.

 In the literature, there are a wide variety of approximate solution methods. Heuristics that are problem-dependent and meta-heuristics are problem-independent techniques that can be applied to a broad range of problems. 
 
 Regarding heuristics used in the literature, Redjem et al. \cite{redjem2016operations} defined a heuristic with two-phases. The first one is to search the optimal tours by calculating for each caregiver the shortest travel duration. The second one is to introduce the precedence and the synchronization constraints. Cire et al. \cite{cire2012heuristic} proposed a heuristic adaptation of logic-based Benders decomposition to solve the \ac{HHC} problem. A greedy heuristic enhanced by constraint programming (CP) is used to solve the Benders master problem.   Bowers et al. \cite{bowers2015continuity} modified the Clarke-Wright algorithm \cite{clarke1964scheduling} to include additional components in the “saving”, which are the continuity of  care and the satisfaction of patients’ preferences

 Concerning meta-heuristics, a  very large number of them are applied to the \ac{HHC} problem. These methods can broadly divided into three groups: 1) single solution methods such as \ac{VNS}, \ac{TS} and \ac{SA} \cite{mankowska2014home,afifi2016heuristic,frifita2020vns,trautsamwieser2011securing,liu2013heuristic}; 2) population-based methods such as \ac{GA}, particle swarm optimization (PSO) \cite{akjiratikarl2007pso, liu2013heuristic}  and 3) hybrid methods such as \cite{liu2013heuristic,bredstrom2008combined,bertels2006hybrid,braekers2016bi}. The first group  is very good in exploitation, which ensures the searching of optimal solutions within the given region while the second group is known to be good in exploration, which allows the algorithm to reach different promising regions of the search space. The third group try to combine several methods to benefice from their advantages such as combining exacts methods with heuristics or combine single solutions methods with population-based methods to have both a good exploration and exploitation.

\subsection{Positioning our contributions with respect to the state of art}

We have summarized above the criteria and constraints considered in the \ac{HHCRSP} studies. These features are further detailed in the following literature reviews \cite{cisse2017or,grieco2020operational,fikar2017home}. In \cite{paraskevopoulos2017resource}, the authors presented a review of the general problem of resource constrained routing and scheduling, which has a wide range of applications, such as: forest management, installations maintenance and repairs, airport operations and  home health care.
 
In addition, the \ac{HHCRSP} extends the \ac{VRPTW} to cover the additional constraints related to \ac{HHC}. Drexl et al. \cite{drexl2012synchronization} reviewed the different type of synchronization constraints in the \ac{VRP}. Vidal et al. \cite{vidal2013heuristics} presented a survey of heuristics used in the literature  for multi-attribute vehicle routing problems. In the following, we focus on studies that have considered time windows, multiple services, their synchronization if they are to be performed simultaneously. Then, we propose three types of models: deterministic, stochastic, and multi-objectives to address the \ac{HHCRSP}.

\subsubsection{A deterministic version of the \ac{HHCRSP}}

\begin{table}[ !htb]
	\caption{ Constraints and assumptions considered in home health care problems.}
	\centering
	\def\arraystretch{1.3}
	\begin{threeparttable}
		
		\begin{tabular*}{\textwidth}{@{\extracolsep{\fill}} ll ccc ccc  }
			\hline
			
			Reference & App & $\#$ objs&  Skills & TW  &MTW  & Syn & Horizon  \\   
			
			\hline
			
			Duque et al. \cite{duque2015home}       & LO      &  2	 & --	& --  &  -- &  --&Long\\

			Breakers et al. \cite{braekers2016bi}         & PA       &  2	 & \checkmark	&  \checkmark & -- &  --&  Short\\ 
			Fathollahi et al.\cite{fathollahi2020bi}     & PA       &  2	 & --	&  \checkmark & -- &  --&Long\\

			Liu et al. \cite{liu2017mathematical}     &LA &  2	 	&-- & \checkmark  & --  &  --&Short\\  
			
			Nickel et al. \cite{nickel2012mid}        &LA &  4	 & \checkmark	&  \checkmark & --  &  --&Long\\ 
		Trautsamwieser et al. 	\cite{trautsamwieser2011optimization}     &LA &  7	 & \checkmark	&  \checkmark & --&  --  & Short\\ 
		Hiermann et al. 	\cite{hiermann2015metaheuristics}&LA &  13	 & \checkmark	&  \checkmark & -- &  --&  Short\\ 
	Shi et al. 		\cite{shi2019robust}&LA &  2	 & \checkmark	&  \checkmark & -- &  --&  Short\\ 
	Frifita et al. 		\cite{frifita2020vns}& -- &  1	 & --	&  \checkmark & -- &  \checkmark& Short \\ 
	Lin et al. 		\cite{lin2021matching}& LO &  5	 & \checkmark	&  \checkmark & -- &  \checkmark& Short \\ 
	Decerle et al. 		\cite{decerle2018memetic}&LA& 2	 & --	&  \checkmark & -- &  \checkmark& Short\\
	Bredstrom et al. 		\cite{bredstrom2008combined}&LA& 2	 & --	&  \checkmark & -- &  \checkmark& Short\\
	Redjem et al. 		\cite{redjem2016operations}     &LA &  2	 &-- 	&  \checkmark &   --&  \checkmark& Short\\  
	Mankowska et al. 		\cite{mankowska2014home} &LA &  3  & \checkmark	&  \checkmark & --  &  \checkmark& Short\\   
	Rasmussen et al. 		\cite{rasmussen2012home} &LA &  3  & \checkmark	&  \checkmark & --  &  \checkmark& Short\\

		Afifi et al. \cite{afifi2016heuristic}  &LA &  3  & \checkmark	&  \checkmark & --  &  \checkmark& Short\\ 
		Parragh et al. \cite{parragh2018solving} &LA &  2  & \checkmark	&  \checkmark & --  &  \checkmark& Short\\ 
		H{\`a} et al. \cite{ha2020new} &-- &  1  & 	--&  \checkmark & --  &  \checkmark& Short\\ 
		Dohn et al.  \cite{dohn2011vehicle} &-- &  1  &-- 	&  \checkmark & --  &  \checkmark& Short\\ 		
			
			\hline
			Our model &LA& 2   & \checkmark  &  \checkmark & \checkmark&\checkmark & Short \\
			\hline

		\end{tabular*}
		\begin{tablenotes}
			\small
			\item \textbf{Abbreviations:} \textbf{App} (Approaches), \textbf{LA} (Linear aggregation), \textbf{PA} (Pareto based), \textbf{LO} (Lexicographical order), "\textbf{$\#$ objs}" is the number of objectives considered. Columns: \textbf{Skills}, \textbf{TW} (time windows), \textbf{MTW} (multiple time windows) and \textbf{sync} (synchronization of services) express if these constraints are considered in studies (\checkmark)  or not(-). 
		\end{tablenotes}
	\end{threeparttable}
	
	\label{table:DLR}

\end{table}

Requested services must be provided within patients’ time windows to respect their availability periods. The decision maker should respect as much as possible these periods when scheduling  patients' visits. Two types of time windows are considered in the literature: hard/fixed time windows  \cite{akjiratikarl2007pso,bard2014sequential,liu2017mathematical,redjem2016operations}, or soft/flexible \cite{hiermann2015metaheuristics,mankowska2014home,nickel2012mid,trautsamwieser2011optimization}. In the first case, the decision maker has to schedule the visit within the time window without any earliness or tardiness of services. In the second case, time windows could not be respected, and delays can be accepted with a penalty cost. Some services are bound to specific times of the day (e.g. serving meals or medicine intakes) and are considered as mandatory requests with hard/fixed time windows. However, some services, such as cleaning or bathing, could be considered soft or flexible services, thus allowing the patient to specify multiple time windows in which he will be available to receive them. But, most studies in \ac{HHC} context, consider only a single time window per patient and it is interesting to consider them as multiple  defining all periods in which patients will be available to receive care services \cite{cisse2017or}.

Redjem et al. \cite{redjem2016operations} conducted a survey of many French home care structures, which showed that patients need several care services per day. Therefore, the synchronization of services should be considered if a patient needs more than a caregiver to provide one or multiple simultaneous services such as  dressing, bathing, getting out of bed and dosing medicine \cite{rasmussen2012home}. For example, an overweight patient might need two people to lift her or him. Moreover, the arrival of caregivers at different times  will disrupt patients and increase the time needed to provide requested services. However, starting at the same time could  provide the requested services within an optimal time, which is obviously equal to the time devoted to process the longest of these services. The synchronization becomes more complicated with the limited number of caregivers. Starting time of visits cannot be independently scheduled from the other routes \cite{cisse2017or}, which increases the complexity of the problem. 

 For the  deterministic version, we focus on studies dealing with multiple services, their synchronization and multiple time windows. Table \ref{table:DLR} shows that only few studies that consider synchronization of services \cite{redjem2016operations,mankowska2014home,rasmussen2012home,bredstrom2008combined,decerle2019memetic, frifita2020vns,lin2021matching}. However, most of these studies consider the synchronization only for double services. Redjem et al. \cite{redjem2016operations} considered the synchronization between multiple  services, but patients are supposed to be assigned to caregivers. Moreover, most studies assign to each patient a single time window in which he would be available to receive care services.
 
 In the first step, we consider the \ac{HHCRSP} with multiple times windows. We propose a new mathematical model as well as a \ac{GVNS} based heuristic to solve large instances. We show through a comparison the advantage of using multiple time windows for patients then we generalize the model to deal with multiple services and their synchronization if they are required to be simultaneous.

 \subsubsection{A stochastic  version of the \ac{HHCRSP}}
 
 Several models and methods have been proposed in the literature  to deal with the \ac{HHCRSP} but most of them are deterministic and generally less robust. The predefined schedule should be adapted for any change in practical situations. Otherwise, there will probably be delays in the services for patients who have not yet been visited, which will cause their dissatisfaction. Travel and service times are critical elements in the planning, any change could affect the overall planning and service quality would be poor or even risky. The uncertainty in travel times may be due to common factors such as varying road conditions, driving skills and weather conditions \cite{shi2019robust}. However, the service time is not always fixed as estimated due to practical reasons, such as diagnosing time and parking situations \cite{shi2019robust}.
 
 Most previous efforts have been focused on studying
 the \ac{VRP} with stochastic parameters such as demands, travel and service times \cite{laporte1992vehicle,li2010vehicle,tacs2014time,tacs2014vehicle,luo2016adaptive,errico2016priori,marinaki2016glowworm,mendoza2016hybrid}. But there are only a few  works that have dealt with the \ac{HHCRSP} with stochastic parameters (e.g., \cite{shi2018modeling, shi2019robust,yuan2015branch,bazirha2020scheduling}) (see Table \ref{table:SLR}). 
 These studies only consider a single availability period as well as a single service operation per  patient and no study, as far as we know, has used multiple time windows or multiple services for patients for the stochastic \ac{HHCRSP}. In the following we propose two \ac{SPR} models to deal with the uncertainty in travel and service times.

 In the first \ac{SPR} model, we assume that patients can request multiple services, which involve assignment of several caregivers to perform these services. In addition, if a patient requested simultaneous services, caregivers must start at the same time, which involve the synchronization. It will be hard to respect patients’ time windows, especially if simultaneously services are requested.  Therefore, the recourse is defined as penalty cost for patients’ delayed services and a remuneration for caregivers' overtime. 
 
 In the second \ac{SPR} model, we assume that time windows are hard/fixed, which must be respected without any earliness or tardiness. To ensure that that requested services are provided within patients' time windows, the recourse is fixed as skipping patients when their time windows will be violated. To increase the chance of providing the maximum of services, patients are allowed to specify multiple time windows.

 \begin{table}[h]
 	\centering
 	\caption{ Stochastic parameters considered in \ac{VRP} and \ac{HHCRSP} problems   }
 	
 	\def\arraystretch{1.3}
 	
 	\begin{tabular*}{\textwidth}{@{\extracolsep{\fill}}lll }
 		\hline
 		
 		Reference  & Stochastic parameters& Problem    \\   
 		
 		\hline
 		Laporte et al. \cite{laporte1992vehicle}     &  Travel times& \ac{VRP}  \\
 		Li et al. \cite{li2010vehicle}   &  Travel and service times & \ac{VRP}  \\
 		Ta{\c{s}} et al. \cite{tacs2014time,tacs2014vehicle}   &  Travel times & \ac{VRP}  \\
 		Errico et al.\cite{errico2016priori}     &  Service times & \ac{VRP}  \\		
 		Marinaki et al.\cite{marinaki2016glowworm}    &  Demands & \ac{VRP}  \\	
 		Luo et al.\cite{luo2016adaptive}     &  Demands & \ac{VRP}  \\
 		Mendoza et al.\cite{mendoza2016hybrid}    &  Demands & \ac{VRP}  \\
 		Yuan et al. \cite{yuan2015branch}     & Service times & \ac{HHCRSP} \\
 		Shi et al.\cite{shi2018modeling,shi2019robust}   & Travel and service times & \ac{HHCRSP} \\
 		\hline
 		
 	\end{tabular*}

 	\label{table:SLR}

 \end{table}
 \subsubsection{A multi-objective version of the HHCRSP}
 
 Most studies transform multi-objective problems into a mono-objective case using aggregation techniques (see Table \ref{table:DLR}). Assigning weights to conflicting  objectives is a confusing task and requires the decision maker' experience and knowledge of the problem. In addition, a sensitivity analysis  of weights should be conducted, which makes the process of solving multi-objective problems more complicated. Duque et al. \cite{duque2015home} used a lexicographical order to solve a bi-objective problem  that aims to maximize the service level and to minimize caregivers' traveling distance. This approach is based on a priori knowledge of some decision-maker preferences to establish the lexicographic order, which is not a simple task.
 \par
 Few studies \cite{braekers2016bi,decerle2019memetic,fathollahi2020bi} have dealt with the multi-objectives case in the \ac{HHC} context by using methods based on the concept of Pareto dominance to approximate the Pareto front. Decerle et al. \cite{decerle2019memetic} proposed  a memetic algorithm for multi-objective optimization of the \ac{HHC} problem, which is a hybridization of \ac{NSGA-II} with the multi-directional local search (MDLS) \cite{tricoire2012multi}. Three objectives were proposed to be optimized, which are minimizing  caregivers' total working time , maximizing the quality of service and minimizing the maximal working time difference among caregivers. Braekers et al.  \cite{braekers2016bi} proposed a heuristic embedding a large neighborhood search heuristic in a MDLS framework to solve a bi-objective problem. The two objectives considered are the costs and clients' inconvenience. Fathollahi et al. \cite{fathollahi2020bi} worked on the bi-objective \ac{HHCRSP} considering patients’ satisfaction for a fuzzy environment. Jimenez’s approach is used to deal with the fuzzy parameters in the objectives and constraints. The authors implemented  NSGA-II, social engineering optimizer (SEO) and simulated annealing algorithms to solve the bi-objective HHCRSP.
 
  As presented above, there is so far no work that has been proposed to solve multi-objective problems in the \ac{HHC} context using the \ac{MOEA/D} algorithm. In this work, we adopt the Pareto and decomposition based approaches to deal with the multi-objective HHCRSP.  According to \cite{zhang2007moea}, it is very time-consuming, if not impossible, to obtain the entire Pareto front for most multi-objective problems. In addition, the decision maker may not be interested in having  all Pareto optimal  solutions. \ac{MOEA/D} algorithm solves multi-objective problems faster since ranking solutions based on Pareto dominance is not used. This study deals  with three objectives so \ac{NSGA-II} algorithm  remains efficient for this kind of optimization problems. According to \cite{trivedi2016survey}, efficiently combining dominance and decomposition based approaches can result in high performance many objective optimizers. MOEA/D, \ac{NSGA-II} and the hybrid \ac{NSGA-II} with \ac{MOEA/D} algorithms are used to approximate the Pareto front, in other words, find the non dominated solutions set.
  
  
    \vspace{1em}
  
  \section{Conclusion}

In this chapter, we have presented a review of the extensive literature related to the \ac{VRP} as well as the home health care routing and scheduling problem.

The literature review shows that studies have been not considered the multiple time windows for patients as well as most of them have dealt with the synchronization   only of double services. Most models proposed in studies are deterministic and powerless to deal with the uncertainty that can arise any time in parameters such as demand, travel and service times. In addition, aggregation techniques are used by most of researchers when multi-objectives need to be simultaneously optimized.

Finally, we presented our contributions in relation to the existing literature. To give more flexibility to decision makers, we dealt with three different problems from the optimization viewpoint: deterministic, stochastic, and multi-objectives. The next three chapters of this thesis are devoted to present our contributions.

  
  \chapter{Deterministic home health care routing and scheduling problem} \label{ChapterDM}
 
 \section {Introduction}
 
 In this chapter, we study the  \ac{HHCRSP} with particular interest in time windows, skills requirement and synchronization constraints. As presented in the previous chapter, the  \ac{HHCRSP} is subject to many variants and extensions. Among the most frequently considered constraints are time window and synchronization constraints. 
 
 In the literature, most  studies assign a single time window  for each patient in which he is available to be visited.  Some services are bound to specific times of the day (e.g. serving meals or medicine intakes) and are considered as mandatory requests with hard/fixed time windows. 
 However, some services, such as cleaning or bathing, could be considered soft or flexible services, thus allowing the patient to specify multiple time windows in which he will be available to receive them. But, most studies in \ac{HHC} context, consider only a single time window per patient and it is interesting to consider them as multiple defining all periods in which patients will be available to receive care services.
 
 In addition, Redjem et al. \cite{redjem2016operations} conducted a survey of many French home care structures, which showed that patients need several care services per day. Therefore, the synchronization of services should be considered if a patient needs more than caregiver to provide one or multiple simultaneous services such as  dressing, bathing, getting out of bed and dosing medicine \cite{rasmussen2012home}. 
 
 We propose a new mathematical model to simultaneously deal with multiple services, their synchronization if they are required to be simultaneous, multiple time windows  per patient and  skills requirements. 
 Although exact methods provide the optimal solution, the computation time increases monotonically with the size of the problem. Due to the weakness of local search methods that fall into a local optimum and cannot escape it, several heuristics extending the local search methods have been proposed  such as a variable neighborhood search \cite{mladenovic1997variable}. 
 We apply the  \ac{GVNS}  based heuristic to solve large instances of the problem since its performance is proven on the \ac{HHCRSP} \cite{mankowska2014home,trautsamwieser2011optimization,shi2019robust,bazirha2019daily}. Furthermore, the proposed solution coding takes into account  both caregivers‘ routes and their assignment. VNS is more suitable since several neighborhood structures could be used to better explore the search space. 

 In the first step, we consider the \ac{HHCRSP} with multiple soft/flexible times windows. We propose a new mathematical model as well as a  \ac{GVNS}  based heuristic to solve large instances. We aim to minimize earliness and tardiness of providing services as well as to minimize caregivers’ waiting times. Then, we show through a comparison the advantage of using multiple time windows for patients.  

In the second step, we generalize the model to deal with multiple services and their synchronization if they are required to be simultaneous by considering two new objective functions and  hard/fixed time windows.  A strategy to select a time window for each patient and to ensure  synchronization of simultaneous services, is proposed and combined with  \ac{GVNS}  based heuristic. This strategy is mandatory since, for each patient, the selection of a time window could not be done independently of his availability periods. In fact, since a single period should be chosen to provide multiple services, this choice may be suitable for one caregiver but not for another.   \ac{GVNS}  generates solutions (caregivers' routes and their assignment) while  this strategy  selects, for each patient, a time window and ensures the synchronization of simultaneous services.
 
 This chapter is organized as follows: The definition of the problems studied are described in sections \ref{PMTW} and \ref{PMSMTW}. In sections \ref{MFMTW} and \ref{MFMSMTW}, the problems are formulated as  mixed integer linear programming models  with a description of the different parameters, variables and constraints taken into account. Sections \ref{VNSMTW} and \ref{VNSMSMTW} describe the \ac{GVNS}. The test instances, the experimental settings and the performance of the  \ac{GVNS}  are presented in sections \ref{NEMTW} and \ref{NEMSMTW}. Finally, the chapter ends with some remarks and conclusions in section \ref{DMConclusion}.

  \section{Home health care routing and scheduling problem with multiple soft time windows} \label{SectionDMTW}

 \subsection{Problem statement} \label{PMTW}

Given a set of patients $ N= \{1, 2,\ldots, n\} $ where $n$ is the number of patients to visit on a day. 
They request services among S the set of services, $ S= \{1, 2,\ldots, q\} $ where $q$ is the number of services  that a \ac{HHC} company could provide. Each service $s$ requested by a patient $i$ has a service duration $t_{is}$.  Each patient $i$ possesses one or multiple time windows $\{ [a_{il}, b_{il}], l\in L= \{1, 2,\ldots,p\} \}$ in which he is available to receive caregivers. Parameters $a_{il}$ and $b_{il}$ are respectively, the earliest and latest service times of the availability period $l$ for the patient $i$ and $p$ is the number of time windows in the day. However, only one period should be selected to provide the patient with the requested services. Travel times between patients $i$ and $j$ are denoted by $T_{ij}$. Patients’ requested services are expressed by a matrix of binary parameters $ \delta_{is} $ where each patient is assumed requesting a single service operation.

\par
The set of available caregivers is denoted $K=\{1,2,\cdots, c\}$ with $c$ is their total number;  $S$ denotes also the set of different skills of caregivers. Each caregiver has a duty length $ [d_k, e_k ]$,  where $d_k$ and $e_k$ are, respectively, the earliest and latest service times. All caregivers start and finish their tours at the center of the \ac{HHC} company, which is represented by artificial nodes $0$ and $n+1$. They use the same mode of transport and are not allowed to work overtime. Their skills are expressed by a matrix of binary parameters  $ \Delta_{ks}$.

 \par 
 
The goal is to define a daily schedule that  minimizes caregivers’ waiting times and respect as much as possible  the selected  time window for each patient. Caregivers leaving the \ac{HHC} center must return within their duty length without exceeding the maximum of earliness and tardiness for providing services as well as the maximum of the waiting time fixed by the decision maker. Their assignment to patients must respect skills requirement.

 \par
 The main hypotheses of this problem are:
 \begin{enumerate} \label{HDMTW}

 	\item   The \ac{HHC} company provides a set of services;
 
 	\item	Caregivers start and finish tours at the \ac{HHC} center. They use the same mode of transportation and the overtime is not allowed;  
 
 	\item	Each caregiver has a time window and a subset of skills;

 	\item   Processing times of service operations are known and cannot be preempted;
 	
 	\item   Travel times between patients are known;
 	
 	\item   Each patient  has multiple time widows, but only one period will be selected;
 	 
 	\item	Each patient requests a single service;
 
 	\item	Caregivers depart as they are available from the \ac{HHC} center, i.e. waiting at \ac{HHC} center is not allowed;

 	\item	Patients' time windows are soft/flexible, it can be violated with a penalty cost;
 
 	\item	 Each service must be provided without exceeding the maximum of earliness and tardiness;
 
 	\item	Each caregiver's waiting time must not exceed the maximum value.

 \end{enumerate}}
 
 \subsection{Mathematical formulation}\label{MFMTW}
 In the following, we first present the notations used in the sequel and then the proposed mathematical formulation for the \ac{HHCRSP} with multiple time windows.
 \subsubsection{Sets}
 \begin{itemize}
 	\item $ N =\{1,2,\ldots,n\}$: set of patients; 
 	\item $ N^0 = N \cup \{0\} $ and $N^{n+1}=N\cup\{n+1\} $:   set of patients including the \ac{HHC} center, which is represented by artificial nodes $0$ and $n+1$;
 	\item $K =\{1,2,\ldots,c\}$:  set of caregivers;
 	\item $L =\{1,2,\ldots,p\}$:  set of patients’ time windows (availability periods).  Each patient $i$ has $p$ periods of availability: Actually, there are $L_i$ valid periods, the others (i.e. $p-L_i$) are null;
 	
 	\item 	$S =\{1,2,\ldots,q\}$:  set of services (skills).
 \end{itemize}
 \subsubsection{Parameters}

 \begin{itemize}
 	
 	\item	$ \alpha , \beta, \gamma $: the weights respectively, of total earliness and total tardiness of service operations, and caregivers’ total waiting time where $ \alpha + \beta + \gamma =1$;
 	\item 	$M$: large number;
 	\item $ [a_{il} , b_{il}] $: the $l^{th}$ availability period of the patient i;
 	\item $ [d_k, e_k ]$: caregivers’ time windows;
 	\item $T_{max}$: maximal tardiness of a service operation;
 	\item	$E_{max}$: maximal earliness of a service operation;
 	\item $W_{max}$: maximal waiting time for each caregiver;
 	\item $ T_{ij}$: travel time from the patient $i$ to the patient $j$;
 	\item $ t_{is}$: processing time of the service operation $s$ at the patient $i \in N$;
 	\item $ \delta_{is}$:equals to 1 if a patient $i \in N$ requires service the operation $s \in S$;
 	\item $ \Delta_{ks} $:equals to 1 if the caregiver $k \in K$ is qualified to provide the service operation $s \in S$.
 	
 \end{itemize}
 \subsubsection{Decision variables}
 
 \begin{itemize}
 	\item $ x_{ijk}$: binary, 1 if the caregiver $k$ visits the patient $j$ after the patient i, 0 otherwise;
 	\item $ y_{iks}$: binary, 1 if the service operation $s$ is provided by the caregiver k to the patient i, 0 otherwise;
 	\item $ z_{il}$: binary, 1 if the $l^{th}$ availability period will be chosen for the patient $i$, 0 otherwise;
 	\item $ u_i$: earliness of a service operation at the patient $i$;
 	\item $ v_i$: tardiness of a service operation at the patient $i$;
 	\item $ A_{ik}$: arrival time of the caregiver $k$ to the patient $i$;
 	\item	$ S_{ik}$: start time of a service operation at the patient $i$ provided by the caregiver $k$;
 	\item $ W_k$: total waiting time of the caregiver $k$.

 \end{itemize}

 \subsubsection{Mathematical model}\label{MDMTW}
 
 The MILP formulation of the problem statement, is an extension of \ac{VRPTW}  \cite{Solomon1988} adapted and augmented by constraints that are specific to the \ac{HHC} context,  is defined as follows:

 \begin{flalign} 
 \min \quad Z = \sum_{i=1}^n (\alpha  u_i + \beta  v_i) + \sum_{k=1}^c \gamma  (S_{ik} - A_{ik})  && \nonumber
 \end{flalign}

 \quad s.t.
 
 \setlength{\belowdisplayskip}{0pt} \setlength{\belowdisplayshortskip}{0pt}
 \setlength{\abovedisplayskip}{0pt} \setlength{\abovedisplayshortskip}{0pt}
 
 \begin{flalign} \label{eqn:ptVisits1STW}
 	\sum _{i=0}^n \sum _{k=1}^c x_{ijk}= 1  &&  \forall j\in N 
 \end{flalign}
 
 \begin{flalign}\label{eqn:ptVisits2STW}
 	\sum _{j=1}^{n+1} \sum _{k=1}^c x_{ijk}= 1 && \forall i\in N 
 \end{flalign}

 \begin{flalign}\label{eqn:leftSTW}
 	\sum _{i=0}^n  x_{i(n+1)k}= 1         &&  \forall k\in K 
 \end{flalign}
 
 \begin{flalign}\label{eqn:retrunSTW}
 	\sum _{j=1} ^{n+1}  x_{0jk}= 1         && \forall k\in K
 \end{flalign}
 
 \begin{flalign}\label{eqn:fluxConservationSTW}
 	\sum _{i=0} ^{n}  x_{imk}= \sum _{j=1}^{n+1}x_{mjk} &&  \forall m\in N,  k\in K
 \end{flalign}
 
 \begin{flalign}\label{eqn:startingTimeSTW}
 	S_{ik}+\sum _{s=1}^{q}t_{is}  y_{iks}  + T_{ij}  \leq S_{jk}+(1-x_{ijk})M &&\forall
 	i\in N^0, j\in N^{n+1},  k\in K
 \end{flalign}
 
 \begin{flalign} \label{eqn:arrivalTime1STW}
 	S_{ik}+\sum _{s=1}^{q}t_{is}  y_{iks}  + T_{ij}  \leq A_{jk}+(1-x_{ijk})M && \forall  i\in N^0, j\in N^{n+1},  k\in K 
 \end{flalign}
 
 \begin{flalign} \label{eqn:arrivalTime2STW}
 	S_{ik}+\sum _{s=1}^{q}t_{is}  y_{iks}  + T_{ij}  \geq A_{jk} + (x_{ijk} - 1 )M && \forall   i\in N^0, j\in N^{n+1},  k\in K 
 \end{flalign}

 
\begin{flalign} \label{eqn:waitingTimeToZeroSTW}
	A_{ik} \leq \sum_{s=1}^q y_{iks} M   && \forall  i\in N, k\in K
\end{flalign} 

\begin{flalign}\label{eqn:StartTimeToZeroSTW}
	S_{ik} \leq \sum_{s=1}^q y_{iks}M   && \forall i\in N, k\in K
\end{flalign}


\begin{flalign} \label{eqn:yDefinitionSTW}
	\sum _{j=1} ^{n+1}  x_{ijk}= \sum_{s=1}^q y_{iks}  && \forall   i\in N,  k\in K
\end{flalign}

	\begin{flalign}\label{eqn:skillsRequirementSTW}
		y_{iks} \leq \delta_{is}\Delta_{ks}   && \forall i\in N, s \in S, k\in K
\end{flalign}

 \begin{flalign}  \label{eqn:twCaregiver1STW}
 S_{0k} =  d_{k}   && \forall  k\in K  
 \end{flalign}
 
 \begin{flalign} \label{eqn:twCaregiver2STW}
 A_{(n+1)k} \leq  e_{k}    && \forall  k\in K
 \end{flalign}

	\begin{flalign} \label{eqn:twPatient1STW}
		(\sum_{l=1}^pz_{il} + \sum_{s=1}^q y_{iks} -2) M +  \sum_{l=1}^p a_{il}z_{il} - u_i \leq S_{ik}  && \forall 
		i\in N, k\in K
	\end{flalign}

	\begin{flalign} \label{eqn:twPatient2STW}
		S_{ik} + \sum_{s=1}^q t_{is} y_{iks}   \leq \sum_{l=1}^p b_{il}z_{il} + v_i + (2 - \sum_{l=1}^pz_{il} -\sum_{s=1}^q y_{iks})M   && \forall
		i\in N, k\in K
\end{flalign}

\begin{flalign}\label{eqn:zDefinitionSTW}
	\sum_{l=1}^p z_{il} =1   && \forall  i\in N
\end{flalign}

 \begin{flalign}\label{eqn:Emax}
 u_i \leq E_{max} && \forall  i\in N
 \end{flalign}

 \begin{flalign} \label{eqn:Tmax}
 v_i \leq T_{max}  && \forall  i\in N
 \end{flalign}

 \begin{flalign}\label{eqn:TWT}
 W_k = \sum_{i=1}^n (S_{ik} - A_{ik}) && \forall  k\in K  
 \end{flalign}

 \begin{flalign}\label{eqn:Wmax}
 W_k \leq W_{max} && \forall  k\in K  
 \end{flalign}

 \begin{flalign} \label{eqn:xZeroDef}
 x_{iik} =0         && \forall  i\in N,  k\in K  
 \end{flalign}
 
 
 \begin{flalign}
 S_{ik} \geq 0  && \forall  i\in N,  k\in K  
 \end{flalign}
 
 \begin{flalign}
 A_{ik} \geq 0  && \forall   i\in N,  k\in K 
 \end{flalign}
 
 \begin{flalign}
 u_i \geq 0  		&& \forall   i\in N
 \end{flalign}
 
 \begin{flalign}
 v_i \geq 0  	&& \forall   i\in N
 \end{flalign}

 \begin{flalign}
 x_{ijk} \in \{0,1\} && \forall   i\in N, j\in N, k\in K  
 \end{flalign}

 \begin{flalign}
 y_{iks} \in \{0,1\} && \forall   i\in N, k\in K, s\in S 
 \end{flalign}
 
 \begin{flalign} \label{eqn:zDef}
 z_{il} \in \{0,1\}      && \forall   i\in N, l\in L  
 \end{flalign}
 \\
 
 The objective function aims to minimize the total penalized earliness and tardiness of services operations as well as caregivers’ total waiting times. 
 Constraints (\ref{eqn:ptVisits1STW}) and (\ref{eqn:ptVisits2STW}) state that each patient will be visited exactly by one caregiver. Constraints (\ref{eqn:leftSTW}) and (\ref{eqn:retrunSTW})  state that each caregiver who has left the center to visit assigned patients must return to that center. Constraints (\ref{eqn:fluxConservationSTW}) express the flux conservation. Constraints (\ref{eqn:startingTimeSTW}) determine the service operations’ starting time of the patient $j$ with respect to service operations’ completion time of the patient $i$. These constraints enforce that the starting time of services along the route of a caregiver are strictly increasing. In doing so, they also eliminate sub-tours because a return to an already visited patient would violate the start time of the previous visit \cite{mankowska2014home}. Constraints (\ref{eqn:arrivalTime1STW}) and (\ref{eqn:arrivalTime2STW}) define the arrival time of a caregiver $k$ to the patient $j$. Constraints (\ref{eqn:waitingTimeToZeroSTW}) and (\ref{eqn:StartTimeToZeroSTW}) initialize the waiting and starting times to zero if the caregiver $k$ will not be affected to the patient $i$. Constraints (\ref{eqn:yDefinitionSTW}) define the variable $y_{iks}$. Constraints (\ref{eqn:skillsRequirementSTW}) ensure that a qualified caregiver $k$ performs a requested service operation $s$ to patient $i$. Constraints (\ref{eqn:twCaregiver1STW}) and (\ref{eqn:twCaregiver2STW}) enforce the respecting of caregivers’ time windows. Constraints (\ref{eqn:twPatient1STW}) and (\ref{eqn:twPatient2STW}) ensure the respecting of patients’ time windows. Constraints (\ref{eqn:zDefinitionSTW}) guarantee that only one period is selected  from patient’s time windows. Constraints (\ref{eqn:Emax}) and (\ref{eqn:Tmax}) guarantee not to exceed the maximal earliness and tardiness of a service operation. Constraints (\ref{eqn:TWT}) define the total waiting time for each caregiver. Constraints (\ref{eqn:Wmax}) ensure not to exceed the maximal waiting time for each caregiver. Constraints (\ref{eqn:xZeroDef} to \ref{eqn:zDef}) set the domains of the decision variables.
 
 \subsection{Variable neighborhood search} \label{VNSMTW}
 
 The drawback of local search strategies is known as falling into a local optimum with a poor value. Several heuristics, extend and improve the local search strategies, have been proposed to avoid being trapped in a local optima such as  \ac{TS} \cite{glover1986future},  \ac{SA} \cite{kirkpatrick1983optimization} and \ac{VNS} \cite{mladenovic1997variable}. \ac{VNS}    is based on the idea of systematic changes of neighborhoods in the search to find a better solution. When an initial solution is defined, the\ac{VNS}    proceeds by a descent method exploring the predefined neighborhoods of a solution to find a local minimum. Each time the descent method is trapped in a local optimum, the shaking phase is randomly applied to generate a new solution and start over the search. Many versions of\ac{VNS}    are used in the literature such as: 1) variable neighborhood descent (VND) is a deterministic version of VNS, the defined neighborhoods are applied to the initial solution in a predefined order, the searching  restart from the first neighborhood when a new local minimum is found. 2) reduced variable neighborhood search is a pure stochastic search method where only the shaking phase is applied to solutions; 3)  \ac{GVNS}  is a\ac{VNS}    where the local descend method is replaced by the VND.

 	\subsubsection{Encoding}
 A solution will be represented by a matrix where the number of columns equals to the number of patients $n$. Two lines are used, the first one will contain patients and requested services operations (included in parenthesis), and the second will contain assigned caregivers. 
 
 Example: we assume that we have 6 patients and 2 caregivers skilled to provide 3 types of services operations. A solution will be encoded as described in Table \ref{table:EncodingSS}.
 
 \begin{table}[!htb]
 	\caption{Example of solution encoding}
 	
 	\centering
 	\begin{tabular}{ | c | c | c | c | c | c | c| }
 		\hline
 		Patients    & 1 \color{black}(3)  & 3 \color{black}(1) & 4 \color{black}(3) & 2 \color{black}(1) &  6 \color{black} (2) & 5 \color{black}(2) \\ 
 		\hline
 		Caregivers  & 1    & 2    & 2    & 1    & 2     &   1\\   
 		\hline
 	\end{tabular}
 	\vspace{3mm}
 	
 	\label{table:EncodingSS}
 \end{table}
 
 The caregiver 1 will visit  patients 1, 2 and 5  to provide respectively  the services 3, 1 and 2.
 
 \subsubsection{Decoding}
 
 \begin{algorithm}[!htb]

 	\SetAlgoLined
 	
 	\eIf{$(\alpha \leq \gamma)$}{
 		
 		\uIf{$ (ET_{il} \leq E_{max} )$}{
 			\textbf{set} $ u_i \longleftarrow ET_{il}$  \;
 		}
 		\uElseIf{$ (W_k + ET_{il} - E_{max} \leq W_{max} )$}{
 			\textbf{set} $ u_i \longleftarrow E_{max}$  and $ W_k \longleftarrow W_k + ET_{il} - E_{max}$  \;
 			
 		}
 		\Else{
 			the availability period $l$ is infeasible 
 		}
 	}
 	{
 		\uIf{$( ET_{il} + W_k \leq W_{max}) $}{
 			
 			\textbf{set} $ W_k \longleftarrow W_k + ET_{il} $  \;
 		}
 		\uElseIf{$ ( W_k + ET_{il} - W_{max} \leq E_{max} )$}{
 			\textbf{set}  $ u_i \longleftarrow W_k + ET_{il} - W_{max}$  and $ W_k \longleftarrow W_{max}$    \;
 			
 		}
 		\Else{
 			the availability period $l$ is infeasible 
 		}
 		
 	}
 	
 	\caption{Caregivers' arriving early to patients  }
 	
 	\label{CaregiversArriveTimes}
 \end{algorithm}
 
 Given a solution encoded as proposed above. For each sub set of patients assigned to a caregiver, the starting, the arrival and the waiting times will be calculated iteratively in the same order as they appear at the matrix. For each patient $i$, the period $ l\in L$ , constraints (1) to (31) are taken into account, that minimize the waiting time, the earliness and the tardiness of the service operation is selected (see equation \ref{selectedPeriod}).

 

 \begin{equation}
 \label{selectedPeriod}
\argminA_{l\in L} \{\alpha  u_i + \beta  v_i + \gamma  (S_{ik}-A_{ik}) \; | \; s.t. \; Constraints \; (1) \; to \; (31)\} 
 \end{equation}
 
 \vspace{3mm}
 For each period $l \in L $, three possible cases of caregivers’ arrival times are to distinguish:
 
 \begin{enumerate}
 	\item The caregiver arrives to a patient and finish providing requested service operation within the availability period: 
 	$$   u_i=0,\quad v_i=0 \quad and \quad  S_{ik} - A_{ik}=0$$
 	
 	\item The caregiver arrives to a patient within the availability period  and finish providing requested service operation with a tardiness time:
 	$$   u_i=0,\quad v_i=A_{ik} + t_{is}-b_{il}\quad and \quad S_{ik} - A_{ik}=0$$
 	\item The caregiver arrives to a patient  before the availability period . In this case, many possibilities arise to calculate the waiting time and the earliness of the service operation. This problem could be formulated as a MIP problem to determine the optimal combination. However, the algorithm \ref{CaregiversArriveTimes} is used. The early time of the caregiver $k$  at patient $i$ for the availability period $l$ ($ET_{il}$) is defined by the following formula: $ET_{il} = A_{ik}  - a_{il} $. 
 	
 \end{enumerate}

 \subsubsection{Neighborhoods }
 The neighborhood of a solution is defined as a transformation function applied to this solution to get a set of solutions where one can move some amount in any direction away from that solution without leaving the set. Four neighborhoods structure are proposed, two neighborhoods are used to intensify patients’ assignment to caregivers (switch and inter-swap) and the two others are used to intensify the order visiting (intra-shift and intra-swap). 
 
 \begin{enumerate}
 	
 	\item 	Switch (i.e. Patients’ reassignment to caregivers): the neighborhood of a solution  is defined as a reassignment of another caregiver $k$ to a patient $i$. The size of possible neighborhoods will be less than $n  (c-1)$ depending on caregivers’ qualifications. The equality  could be hold if all caregivers are skilled to provide all services operations ( see Fig. \ref{fig:SwitchDMTW});

 	\item Inter-swap: this neighborhood aims to change patients’ assignment to caregivers. Given two patients, caregivers’ assignment are swapped. The size of possible neighborhoods equals to    $\frac{(n-1)\times n}{2}$ ( see Fig. \ref{fig:Inter-swapDMTW});
 	
 	\item Intra-shift: given a visiting order, the neighborhood is defined as shifting of a patient to another position. The size of possible neighborhoods equals to $(n-1) \times n$ ( see Fig. \ref{fig:Intra-shiftDMTW});
 	
 	\item Intra-swap: given a visiting order, the neighborhood is defined as two patients' position exchanging. The size of possible neighborhoods equals to $\frac{(n-1)\times n}{2}$ (see Fig. \ref{fig:Intra-swapDMTW}).

 	\subsubsection{Shaking}
 	
 	The shaking phase is the mechanism that allows to escape the local optima with a poor value \cite{mladenovic1997variable}. It consists of applying a series of moves to a solution to jump from a local optimum to a new solution, then restart improving the new solution using local search methods. The four proposed neighborhoods are used for the shaking phase ($ k_{max}$ = 4) as operators, which will be randomly applied $m_1$ times whenever the local search method is trapped in a local optimum.

 	\begin{figure}[!htb]
 		\centering
 		\begin{mdframed}[style=MyFrame,nobreak=true,align=center,userdefinedwidth=31em]

 			Original solution\quad :\hspace{1mm}
 			\centering
 			\begin{tabular}{ | c | c | c | c | c | c | c| }
 				\hline
 				Patients    &\textbf{\color{red} 1 (3)}& 3 (1) & 4 (3)& 2 (1)&6  (2)& 5 (2) \\ 
 				\hline
 				Caregivers  &\textbf{\color{red} 1 }& 2    & 2    & 1    & 2     &   1\\   
 				\hline
 				\hline
 			\end{tabular}
 			
 			Neighbor solution\quad :
 			\centering
 			\begin{tabular}{ | c | c | c | c | c | c | c| }
 				\hline
 				Patients    &\textbf{\color{red}1  (3)}& 3 (1) & 4 (3)& 2 (1)& 6 (2)& 5 (2) \\ 
 				\hline
 				Caregivers  &\textbf{\color{blue} 2 }& 2    & 2    & 1    & 2         &   1\\   
 				\hline
 			\end{tabular}
 		\end{mdframed}
 		\caption{Example of switch neighborhood moves}
 		\label{fig:SwitchDMTW}
 		
 	\end{figure}

 	\begin{figure}[!htb]
 		\centering
 		\begin{mdframed}[style=MyFrame,nobreak=true,align=center,userdefinedwidth=31em]

 			Original solution\quad :\hspace{1mm}
 			\centering
 			\begin{tabular}{ | c | c | c | c | c | c | c| }
 				\hline
 				Patients    &\textbf{\color{cyan}1  (3)}& 3 (1) & 4 (3)& 2 (1)& \textbf{\color{cyan}6 (2)}& 5 (2) \\ 
 				\hline
 				Caregivers  &\textbf{\color{red} 1 }  & 2    & 2    & 1    &\textbf{\color{blue} 2 } &   1\\   
 				\hline
 				\hline
 			\end{tabular}
 			
 			Neighbor solution\quad :
 			\centering
 			\begin{tabular}{ | c | c | c | c | c | c | c| }
 				\hline
 				Patients    &\textbf{\color{cyan}1  (3)}& 3 (1) & 4 (3)& 2 (1)& \textbf{\color{cyan}6 (2)}& 5 (2) \\ 
 				\hline
 				Caregivers  &\textbf{\color{blue} 2}& 2    & 2    & 1     &\textbf{ \color{red} 1}&1\\   
 				\hline
 			\end{tabular}
 		\end{mdframed}
 		\caption{Example of inter-swap neighborhood moves}
 		\label{fig:Inter-swapDMTW}
 		
 	\end{figure}

 	\begin{figure}[!htb]
 		\centering
 		\begin{mdframed}[style=MyFrame,nobreak=true,align=center,userdefinedwidth=31em]

 			Original solution\quad :\hspace{1mm}
 			\centering
 			\begin{tabular}{ | c | c | c | c | c | c | c| }
 				\hline
 				Patients    &\textbf{\color{blue} 1 (3)}& 3 (1) & 4 (3)& 2 (1)&6  (2)& 5 (2) \\ 
 				\hline
 				Caregivers  &\textbf{\color{blue} 1 }& 2    & 2    & 1    & 2     &   1\\   
 				\hline
 				\hline
 			\end{tabular}
 			
 			Neighbor solution\quad :
 			\centering
 			\begin{tabular}{ | c | c | c | c | c | c | c| }
 				\hline
 				Patients    & 3 (1) & 4 (3)& 2 (1)&\textbf{\color{blue} 1 (3)} &6 (2)& 5 (2) \\ 
 				\hline
 				Caregivers  & 2    & 2    & 1 & \textbf{\color{blue} 1 }   & 2         &   1\\   
 				\hline
 			\end{tabular}
 		\end{mdframed}
 		\caption{Example of intra-shift neighborhood moves}
 		\label{fig:Intra-shiftDMTW}
 		
 	\end{figure}

 	\begin{figure}[!htb]
 		\centering
 		\begin{mdframed}[style=MyFrame,nobreak=true,align=center,userdefinedwidth=31em]

 			Original solution\quad :\hspace{1mm}
 			\centering
 			\begin{tabular}{ | c | c | c | c | c | c | c| }
 				\hline
 				Patients    &\color{red} 1 (3)& 3 (1) & 4 (3)& 2 (1)&\color{blue}6  (2)& 5 (2) \\ 
 				\hline
 				Caregivers  &\color{red} 1    & 2    & 2    & 1    &\color{blue} 2     &   1\\   
 				\hline
 				\hline
 			\end{tabular}
 			
 			Neighbor solution\quad :
 			\centering
 			\begin{tabular}{ | c | c | c | c | c | c | c| }
 				\hline
 				Patients    &\color{blue}6  (2)& 3 (1) & 4 (3)& 2 (1)&\color{red} 1 (3)& 5 (2) \\ 
 				\hline
 				Caregivers  &\color{blue} 2& 2    & 2    & 1    & \color{red} 1         &   1\\   
 				\hline
 			\end{tabular}
 			
 		\end{mdframed}
 		\caption{Example of intra-swap neighborhood moves}
 		\label{fig:Intra-swapDMTW}
 		
 	\end{figure}

 \end{enumerate}

 \subsubsection{Local search methods}
 
 Whenever shaking phase is applied, local search methods are called to improve the new generated solution until they fall in a local optimum. The four proposed neighborhoods are used as local search methods, each one (intra-swap, intra-shift, inter-swap, and switch) is applied to a solution to explore the search space until it falls in a local optimum.   Local search methods adopt two types of algorithm while exploring neighborhoods: i) The first improvement is to start over the search whenever a neighbor that improves the best-known solution is found; ii) The best improvement is to start over  the search when all neighbors are explored and the best one is selected \cite{taillard1990some}. The the best improvement algorithm is adopted.

\subsubsection{Initial solution}

The initial solution is generated randomly as follow:

\begin{enumerate}
	\item  Patients are sorted by increasing end of their time windows for single period. for multiple periods the visiting order is generated randomly (see Table \ref{table:VODMTW});
	
	\begin{table}[H]
	 
		\caption{Example of visiting order}
		\centering
		\begin{tabular}{ | c | c | c | c | c | c | c| }
			\hline
			Patients    & 1 (3)  & 3 (1) & 4 (3) & 2 (1) &  6  (2) & 5 (2) \\ 
			\hline
			Caregivers  &    &    &    &   &      &   \\   
			\hline
		\end{tabular}
		
		\label{table:VODMTW}
	 
	\end{table}
	
	\item For each patient, assign a qualified caregiver selected randomly (see Table \ref{table:AssignDMTW});
	\begin{table}[H]
		\caption{example of caregivers’ assignment to patients}
		\centering
		\begin{tabular}{ | c | c | c | c | c | c | c| }
			\hline
			Patients    & 1 (3)  & 3 (1) & 4 (3) & 2 (1) &  6  (2) & 5 (2) \\ 
			\hline
			Caregivers  &1    &  2  &  2  &  1 &  2    & 1  \\   
			\hline
		\end{tabular}

		\label{table:AssignDMTW}
	 
	\end{table}
	
	\item Calculate caregivers’ waiting times, the earliness and tardiness of services operations using the decoding method proposed above;
	
	\item If the solution is infeasible repeat steps 1,2 and 3. Otherwise, go to the step 5;
	
	\item Calculate the objective function value.
	
\end{enumerate}

 \subsubsection{GVNS algorithm}
 
 \begin{algorithm}[!htb]
 	\SetAlgoLined
 	
 	\textbf{Initialization:} \;
 	\hspace{5mm}	 \textbf{- Set} $K_{max}=4$ and $l_{max}=4$  \;
 	\hspace{5mm}	 \textbf{- Generate} an initial solution $x$ \;
 	\While{(the stopping condition is not reached)}{
 		\For{$k\leftarrow 1$ \KwTo $k_{max}$}  {
 			
 			\textbf{Set} $ x' \longleftarrow$ \textbf{Shaking} ($x$, neighborhood$_k$, $m_2$) \;
 			\For{$l\leftarrow 1$ \KwTo $l_{max}$}{
 				
 				\textbf{Set} $ x'' \longleftarrow$ \textbf{Local search} ($x'$, neighborhood$_l$) \;
 				\eIf{$f(x'') < f(x')$}{
 					\textbf{Set} $ x' \longleftarrow x'' $ \;
 					\textbf{Set} $l \longleftarrow 1$ \;
 				}{
 					\textbf{Set} $l \longleftarrow l+1$ \;
 				}
 			}
 			\eIf{$f(x') < f(x)$}{
 				\textbf{Set} $ x \longleftarrow x' $ \;
 				\textbf{Set} $k \longleftarrow 1$ \;
 			}{
 				\textbf{Set} $k \longleftarrow k+1$ \;
 			}
 		}

 	}	\caption{ \ac{GVNS}  algorithm}
 	\label{GVNSDMTW}
 \end{algorithm}

 \ac{GVNS}  algorithm is executed until there is no improvement over the best-known solution for $m_1$ iterations. Therefore, the counter is set to zero if a better new solution is found.  \textbf{Line 6}  generates a new solution $x’$ by randomly applying $m_2$ times the neighborhood $k$ to the solution $x$.  \textbf{Lines 7-15} describe VND algorithm, which is applied to each solution $x'$ generated in the shaking phase. \textbf{Line 8 }improves solution $x’$ using neighborhood $l$ as a local search method and stores it in $x''$.  At each time, if the solution $x$ (resp. $x'$) is improved, the $k$ (resp. $l$) is initialized to 1 (see Algorithm \ref{GVNSDMTW}).

 \subsection{Numerical experiments } \label{NEMTW}
 The experiments run on the computer with Intel i7-7600U 2.80-GHz CPU and 16 GB of RAM under windows 10. The MIP model is implemented and tested using CPLEX IBM version 12.8. The metaheuristic based on \ac{GVNS}  is coded and tested using the language C++.

 \subsubsection{Test instances} \label{TIDMTW}
 
 Tests instances have been generated randomly using the benchmark instances from  Mankowska et al. \cite{mankowska2014home}. Patients and the \ac{HHC} office are placed at random locations in the area of $ 100 \times 100$ distance units. Travel times $T_{ij}$ are equal to the Euclidean distance between patients truncated to integer. Processing times of services operations $t_{is}$ are randomly chosen from the interval $[10, 20]$. Six types of services $S= \{1, . . ., 6\}$ are considered. Caregivers are grouped into two subsets with different skills. Each caregiver of the first group is qualified for providing at most three services, which are randomly selected from the subset $\{1, 2, 3\}$ of S. Accordingly, each caregiver of the second group is qualified for providing at most three services from subset $\{4, 5, 6\}$. Each patient requires a single service, which is randomly drawn from $S=\{1,…,6\}$. The time windows are of length 120 min (2 h) and are randomly placed within a daily planning period of 10 hours. Regarding instances which contains two availability periods, the first period will be placed in the first 5 hours, and the second period will be placed in the interval [5, 10]. No special preference for the three sub-goals and, therefore, weights are set to $\alpha=\beta=\gamma=1/3$. The maximum earliness and the maximum tardiness of service operations are set respectively to 0min and 15min, the maximum waiting time is set to 90min. 8 instances are generated, each one is used with 1-period and with 2-periods availability of patients. The instance Int1\_1 refers to the instance 1 with one availability period  and the instance Int1\_2 refers to the instance 1 with two availability periods. A series of tests was executed to find the best tuning parameters for the proposed metaheuristic. 
 
 \subsubsection{Computational results}
 Instances are generated as described above and solved. Instances with two  availability periods are generated feasible in the first time. However, instances with single time windows need many regenerations to get a feasible solution. Table \ref{table:NRDMTW} summarizes the results of CPLEX and \ac{GVNS}  according to the sizes and patients' availability periods of test instances. $LB$  is the lower bound  of the model given by the CPLEX solver, $Z$ is the objective function value, $GAP$ is calculated as $100\% \times$ (Z – lower bound)/Z and $CPU$ is the computing time elapsed of solved instances. For the \ac{GVNS}  algorithm, each instance has been executed 10 times and the best, the worst and the average solutions are considered. $GAP$ is calculated as $100\% \times$ (Average – lower bound)/Average. $CPU$ computing time corresponds to the sum of time elapsed to solve each instance, which is running 10 times.  Optimal solutions reached are in boldface.
 
 All instances are solved, using CPLEX IBM, to optimality except  instances Int7\_2 and Int8\_2. Instances Int6\_2, Inst7\_2 and Int8\_2 are hard to solve compared to others instances. The instance Int6\_2 required 47530s to be solved, the instance Int7\_2 and Int8\_2 was executed for 18 hours and CPLEX could not find respectively the optimal solution and a feasible solution, which shows the limit of the exact method used considering the problem is NP-hard. The proposed metaheuristic (GVNS), was able to find a good solution in a very short computational time, solution in bold are proven optimal( see Table \ref{table:NRDMTW}). The proposed \ac{GVNS}  could find  optimal solutions at least one time of the 10 executions for 11 instances. Due to the LB that equals zero, the gap of the instances Int4\_2,  Int5\_2, Int7\_2 and Int8\_2 is 100\%. Because when LB approaches zero, the gap tends to 100\%.

 	
 		
 	

 \begin{figure}[!htb]
 	\begin{center}
 		\begin{tikzpicture}
 		\begin{axis}[
 		width=\textwidth, 
 		height=10cm,
 		grid=major, 
 		grid style={dashed,gray!30},
 		xlabel= Instances, 
 		ylabel= Objective function values,
 		xtick=data,
 		legend style={at={(0.08,0.8)},anchor=west},
 		,xticklabels={Int\_1, Int\_2, Int\_3, Int\_4, Int\_5, Int\_6, Int\_7, Int\_8}
 		]
 		
 		\addplot+[] 
 		table[x=SS, y=TW1, col sep=tab] {src/Tw1VsTw2.txt};

 		\addplot+[] 
 		table[x=SS, y=TW2, col sep=tab] {src/Tw1VsTw2.txt};

 		\legend{Single TW,  Double TW}
 		\end{axis}
 		
 		\end{tikzpicture}
 	\caption{Objective function values according to the number of time windows (TW)}
 	\label{fig:figDMTW}

 	\end{center}
 \end{figure}
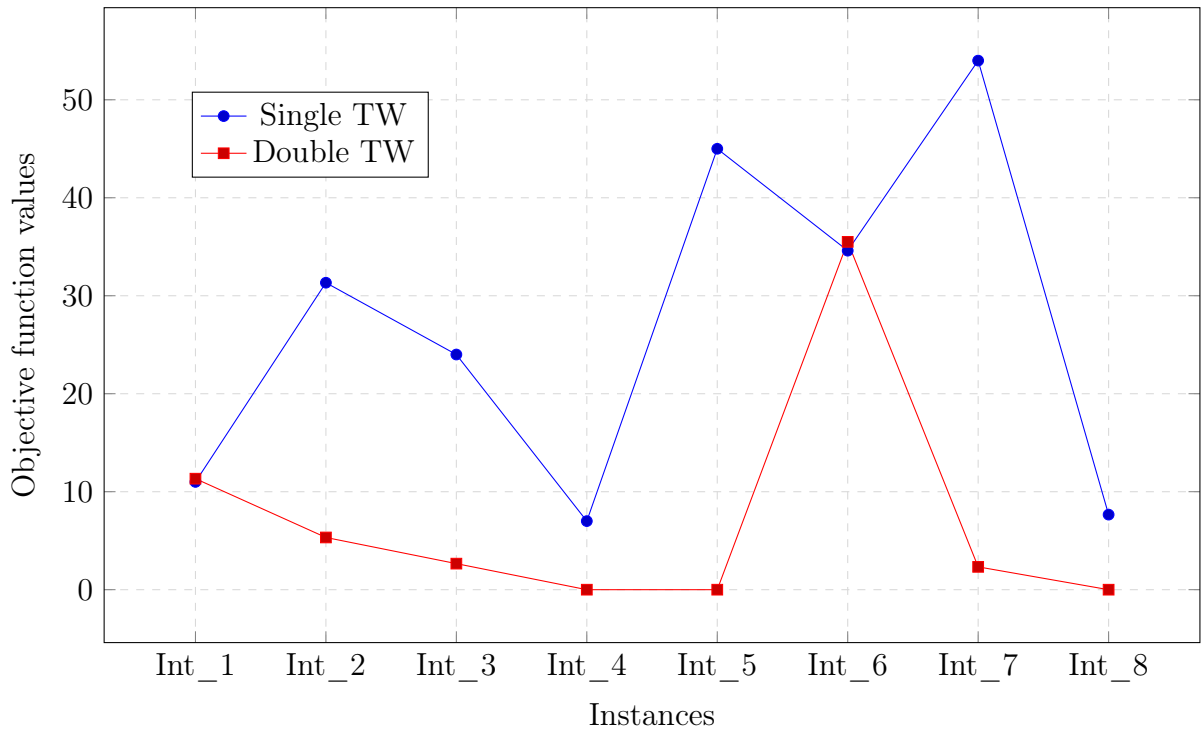

 \begin{table}[!htb]

 	\caption{Numerical results of tested instances}
 	\centering
 	\def\arraystretch{1.3}\setlength\tabcolsep{1.5pt}
 	\begin{tabular*}{\textwidth}{@{\extracolsep{\fill}}  c c  c c  | c c  c | c  c  c  c  c}
 		\hline
 		\multicolumn{4}{c|}{Instances} & \multicolumn{3}{c|}{\textbf{CPLEX} }&\multicolumn{5}{c}{ \textbf{\ac{GVNS}}} \\ 
 		\hline
 		\textbf{} & \textbf{ N} & \textbf{K}    &  \textbf{L}    &\textbf{LB} & \textbf{Z} & \textbf{CPU} & \textbf{Best} & \textbf{Worst}& A\textbf{verage} &\textbf{Gap} & \textbf{CPU}   \\

 		\hline


 		Int1\_1  & 7& 2 &  1 &  11 &\textbf{11} & 1.58 & \textbf{11} & \textbf{11}& \textbf{11} & 0 \% & $<1$   \\   
 		
 		\hline

 		Int2\_1  & 10& 3  &1 & 31.33 &\textbf{ 31.33} & 1.90 & \textbf{31.33} & 31.36& 31.53 & 1 \% & $<1$   \\

 		\hline
 		
 		Int3\_1  & 14& 4  & 1 &24 & \textbf{24} & 2.52 &\textbf{ 24}& 35.33& 27.76 & 13 \% & 1.79   \\    
 		
 		\hline
 		Int4\_1  & 20 & 5  &1 & 7 & \textbf{7} & 3.02 & \textbf{7}& 9.33& 7.43 & 6 \% & 8.82   \\   
 		\hline
 		Int5\_1  & 25 & 6  &1 & 45 &\textbf{ 45}& 5.16 & 59.86& 84& 65.1& 30 \% & 25   \\    
 		\hline
 		Int6\_1  & 30 & 6  &1 & 34.6 &\textbf{ 34.6}& 12.13 & 35.3& 44.6& 39.13& 17 \% & 25   \\    
 		\hline
 		Int7\_1  & 40 & 8 &1 & 54 & \textbf{54} & 255 & 72& 100& 78& 30 \% & 111   \\    
 		\hline
 		Int8\_1  & 50 & 10 &1 & 7.66 & \textbf{7.66} & 381 & 11& 18.66& 13.73& 44 \% & 210   \\    
 		\hline
 	
 		
 		Int1\_2  & 7& 2  &2 & 11.33 & \textbf{11.33} & 1.73 & \textbf{11.33} & \textbf{11.33} & \textbf{11.33} & 0 \% & $<1$   \\   
 		\hline
 		
 		Int2\_2  & 10& 3  &2 & 5.33 &\textbf{ 5.33} & 1.91 & \textbf{5.33} & 6.33& 5.43 & 2 \% & $<1$   \\   
 		\hline

 		Int3\_2  & 14& 4 &2 & 2.66 & \textbf{2.66} & 2.60 & \textbf{2.66}& 4& 3.16 & 16 \% & 1.74   \\    
 		\hline
 		
 		Int4\_2  & 20 & 5  &2 & 0 & \textbf{0} & 16.31 & \textbf{0}& 9& 2.7 & 100 \% & 9.77   \\    
 		\hline
 		
 		Int5\_2 & 25 & 6 & 2 &0 & \textbf{0} & 61 & \textbf{0}& 6.33& 4.366& 100 \% & 32  \\    
 		\hline
 		
 		Int6\_2 & 30 & 6 & 2 &35.3 & \textbf{35.3} & 47530 & \textbf{35.3}& 44.66& 41.96& 15 \% & 52.64  \\    
 		\hline
 		
 		Int7\_2 & 40 & 8 &2 & 0 & 2.33 &64800  & 5.66& 13.66& 7.86& 100 \% & 248  \\    
 		\hline
 		
 		Int8\_2  & 50 & 10 & 2 &0 & - & 64800& \textbf{0}& 1.33& 0.23& 100 \% & 561  \\    
 		\hline
 		
 	\end{tabular*}

 	\label{table:NRDMTW}

 \end{table}
 
 \par

 The objective function values of instances with two availability periods are less than those with single period, except the first and the sixth instances, their values are very close (see Fig. \ref{fig:figDMTW}). Therefore, patients will be more satisfied seeing that the earliness and tardiness of service operations are minimized. On the other hand, the waiting time is also minimized which will increase the productive time of caregivers since the waiting time is considered as unproductive~time.
 
 The soft/flexible time windows would increase the chance of finding a feasible solution since delays are accepted with a penalty cost. However, patients will be less satisfied because their availability periods are not respected. Furthermore, Patients needs several care activities per day and some of them are simultaneous by their nature. In the next section, we extend and adapt the model to deal with multiple hard/fixed time windows for patients and multiple synchronized services as well as to balance caregivers' workload.

 

 \clearpage

 \section{Home health care routing and scheduling problem with multiple hard time windows and multiple synchronized services } \label{SectionDMSMTW}

 \subsection{Problem statement}\label{PMSMTW}

We modify the problem statement defined the previous section \ref{PMTW} to allow patients to request multiple services, which can be provided as independent or simultaneous services. Also, we suppose that patients' time windows are hard/fixed and we consider two new objective functions.
 
 \par 
 
 In this problem, we aim to minimize waiting time and balance caregivers' workload. The waiting time ($f_1$) is considered to be unproductive \cite{redjem2016operations} and the caregivers are supposed to be paid for their regular working time regardless of the amount of care they provide \cite{braekers2016bi}. On the other hand, to ensure fairness among caregivers, their workload must be balanced ($f_2$). To achieve this goal, we consider the minimization of the sum of absolute difference between each caregiver’ working time and the average working time of all caregivers.  The total caregivers' working times is equal to the sum of service and travel times. 

\par

 Most studies consider the workload as the difference between the extreme values of working times. However, any change between the extremes has no effect on the value of workload \cite{matl2018workload}. In this study, the Mean Absolute Deviation (MAD) is used in our model to take into account all values but it is sensible to extreme values. To overcome this issue, caregivers' time windows are considered and the overtime time is not allowed, which ensures caregivers' working times do not exceed the duty length. 

\par
The goal is to define a daily schedule that minimizes waiting times for caregivers and balances their workload. Each patient should be visited within a time window selected from his availability periods. Services that need to be provided simultaneously should be synchronized, but if it is not necessary they would be performed independently. Caregivers leaving the \ac{HHC} center must return within their duty length and their assignment must respect skills requirement.

\par
The main hypotheses of this problem are:
\begin{enumerate}
	\item  The first 6 hypotheses defined in \ref{HDMTW};
	
	\item	Each patient can request one or more services, which should be  synchronized if they are required to be simultaneous;
	\item   Each requested service must be performed by one and only one caregiver. In others words, providing the same service  for the same patient by several caregivers is not allowed;
	\item   Each caregiver must provide one and only one service. In others  word, a caregiver $k$  is not allowed to carry out several services for the same patient;
	\item   Caregivers' departure are independent and each caregiver leaves  as soon as he  finished.
	
\end{enumerate}}

\subsection{Mathematical formulation} \label{MFMSMTW}
In the following, we first present the new notations, those not defined are the same from the previous section, used in the sequel and then the proposed mathematical formulation for the HHCRSP-MTW-MSS.

\vspace{-2mm}

\subsubsection{Parameters}

\begin{itemize}
	
	\item $\lambda_i$:  equals 1 if requested service operations by patient $i$ require to be simultaneous, and 0 otherwise.
	
\end{itemize}

\subsubsection{Decision variables}

\begin{itemize}
	\item $ u_{ik}$:  equals  1 if the arrival time of the caregiver $k$   is used as the starting time for the patient $i$, 0 otherwise;
	\item $ v_{ik}$:  equals 1 if the starting time $S_{ik}$ is greater than the arrival time $A_{ik}$ of the caregiver k in the case of independent services or  if  the synchronized starting time $SS_i$ is greater than the maximum of  arrival times of  the caregivers at patient $i$ in the case of simultaneous services; 0 otherwise;   
		\item $ D_k$:   absolute value of the difference between the working time of each caregiver $w_k$ and the average working time.
	
\end{itemize}

\subsubsection{Mathematical model}

The mathematical model proposed to deal with the new hypotheses is adapted from  the previous formulation \ref{MDMTW} to cover  additional assumptions, such as multiple services and their synchronization if they have to be simultaneous, minimizing caregivers' waiting times and balancing their workload.  The mathematical formulation of the new model is defined as follows: 

	\begin{flalign} 
		\min Z = \alpha \sum_{k=1}^c\sum_{i=0}^n ( S_{ik} - A_{ik})+ \beta \sum_{k=1}^c D_k  && \nonumber
\end{flalign}

\quad s.t.

\setlength{\belowdisplayskip}{0pt} \setlength{\belowdisplayshortskip}{0pt}
\setlength{\abovedisplayskip}{0pt} \setlength{\abovedisplayshortskip}{0pt}

\begin{flalign} \label{eqn:ptVisits1MSMTW}
	\sum _{i=0}^n \sum _{k=1}^c x_{ijk}=  \sum_{s=1}^q\delta_{js}  &&  \forall j\in N 
\end{flalign}

\begin{flalign}\label{eqn:ptVisits2MSMTW}
	\sum _{j=1}^{n+1} \sum _{k=1}^c x_{ijk}= \sum_{s=1}^q\delta_{is} && \forall i\in N 
\end{flalign}

\begin{flalign}\label{eqn:leftMSMTW}
	\sum _{i=0}^n  x_{i(n+1)k}= 1         &&  \forall k\in K 
\end{flalign}

\begin{flalign}\label{eqn:retrunMSMTW}
	\sum _{j=1} ^{n+1}  x_{0jk}= 1         && \forall k\in K
\end{flalign}

\begin{flalign}\label{eqn:fluxConservationMSMTW}
	\sum _{i=0} ^{n}  x_{imk}= \sum _{j=1}^{n+1}x_{mjk} &&  \forall m\in N,  k\in K
\end{flalign}

\begin{flalign}\label{eqn:startingTimeMSMTW}
	S_{ik}+\sum _{s=1}^{q}t_{is}  y_{iks}  + T_{ij}  \leq S_{jk}+(1-x_{ijk})M &&\forall
	i\in N^0, j\in N^{n+1},  k\in K
\end{flalign}

\begin{flalign} \label{eqn:arrivalTime1MSMTW}
	S_{ik}+\sum _{s=1}^{q}t_{is}  y_{iks}  + T_{ij}  \leq A_{jk}+(1-x_{ijk})M && \forall  i\in N^0, j\in N^{n+1},  k\in K 
\end{flalign}

\begin{flalign} \label{eqn:arrivalTime2MSMTW}
	S_{ik}+\sum _{s=1}^{q}t_{is}  y_{iks}  + T_{ij}  \geq A_{jk} + (x_{ijk} - 1)M && \forall   i\in N^0, j\in N^{n+1},  k\in K 
\end{flalign}


\begin{flalign} \label{eqn:waitingTimeToZeroMSMTW}
	A_{ik} \leq \sum_{s=1}^q y_{iks} M   && \forall  i\in N, k\in K
\end{flalign} 

\begin{flalign}\label{eqn:StartTimeToZeroMSMTW}
	S_{ik} \leq \sum_{s=1}^q y_{iks}M   && \forall i\in N, k\in K
\end{flalign}

\begin{flalign} \label{eqn:yDefinitionMSMTW}
	\sum _{j=1} ^{n+1}  x_{ijk}= \sum_{s=1}^q y_{iks}  && \forall   i\in N,  k\in K
\end{flalign}

\begin{flalign}\label{eqn:skillsRequirementMSMTW}
		y_{iks} \leq \delta_{is}\Delta_{ks}   && \forall i\in N, s \in S, k\in K
\end{flalign}

\begin{flalign} \label{eqn:assign1MSMTW}
	\sum_{s=1}^q y_{iks} \leq 1   && \forall  i\in N,  k\in K
\end{flalign}

\begin{flalign} \label{eqn:assign2MSMTW}
	\sum_{k=1}^c y_{iks} \leq 1   && \forall  i\in N,  s\in S
\end{flalign}

\begin{flalign} \label{eqn:twCaregiver1MSMTW}
	d_{k} =   A_{0k}   && \forall  k\in K
\end{flalign}

\begin{flalign} \label{eqn:twCaregiver2MSMTW}
	A_{(n+1)k} \leq  e_{k}   && \forall    k\in K
\end{flalign}

\begin{flalign} \label{eqn:twPatient1MSMTW}
		(\sum_{l=1}^pz_{il} + \sum_{s=1}^q y_{iks} -2) M +  \sum_{l=1}^p a_{il}z_{il}  \leq S_{ik}  && \forall 
		i\in N, k\in K
\end{flalign}

\begin{flalign} \label{eqn:twPatient2MSMTW}
		S_{ik} + \sum_{s=1}^q t_{is} y_{iks}   \leq \sum_{l=1}^p b_{il}z_{il} +(2 - \sum_{l=1}^pz_{il} -\sum_{s=1}^q y_{iks})M   && \forall
		i\in N, k\in K
\end{flalign}

\begin{flalign}\label{eqn:zDefinitionMSMTW}
	\sum_{l=1}^p z_{il} =1   && \forall  i\in N
\end{flalign}

	\begin{flalign} \label{eqn:synchronization1MSMTW}
		\sum_{k^{\prime}=1}^c S_{ik^{\prime}}  -  \sum_{s=1}^q\delta_{is} S_{ik} \leq (1 -\lambda_i \sum_{s=1}^q y_{iks})M  && \forall  i\in N, k\in K
	\end{flalign}
	\begin{flalign} \label{eqn:synchronization2MSMTW}
		\sum_{k^{\prime}=1}^c S_{ik^{\prime}}  -  \sum_{s=1}^q\delta_{is} S_{ik} \geq (\lambda_i\sum_{s=1}^q y_{iks} - 1)M  && \forall i\in N, k \in K
	\end{flalign}

	\begin{flalign} \label{eqn:YdefinitionMSMTW}
		w_k = \sum_{i=1}^n \sum_{s=1}^q t_{is}y_{iks} +  \sum _{i=0}^n \sum _{j=1}^{n+1} T_{ij}x_{ijk}     && \forall    k\in K
\end{flalign}


	\begin{flalign} \label{eqn:YAbsValue1MSMTW}
		D_k \geq  w_k - \frac{\sum_{k=1}^c {w_k}}{c}     && \forall\  k\in K
	\end{flalign}
 
	\begin{flalign} \label{eqn:YAbsValue2MSMTW}
		D_k \geq  \frac{\sum_{k=1}^c {w_k}}{c}   -  w_k     && \forall  k\in K
	\end{flalign}

	\begin{flalign} \label{eqn:earlyStarting1}
		\sum _{k^{\prime}=1}^{c}S_{ik^{\prime}}  \leq  \sum_{s=1}^q\delta_{is}A_{ik}+(1-\lambda_i u_{ik})M  && \forall  i\in N,  k\in K
	\end{flalign}

	\begin{flalign}\label{eqn:earlyStarting2}
		\sum _{k=1}^{c}S_{ik}  \leq \sum_{s=1}^q\delta_{is} \sum_{l=1}^p a_{il}z_{il}+(1-\lambda_i\sum _{k=1}^{c}v_{ik})M  && \forall  i\in N  
	\end{flalign}

\begin{flalign} \label{eqn:earlyStarting3}
	\sum_{k=1}^c (u_{ik} + v_{ik})\leq 1 +(1-\lambda_i)M  && \forall  i\in N
\end{flalign}

\begin{flalign} \label{eqn:earlyStarting4}
	\sum_{k=1}^c (u_{ik} + v_{ik})\geq 1 +(\lambda_i -1)M  && \forall  i\in N
\end{flalign}

\begin{flalign} \label{eqn:earlyStarting5}
	S_{ik}  \leq  A_{ik}+(1+\lambda_i  -u_{ik})M  && \forall  i\in N,  k\in K
\end{flalign}

\begin{flalign}\label{eqn:earlyStarting6}
	S_{ik}  \leq  \sum_{l=1}^p a_{il}z_{il}+(1 +\lambda_i - v_{ik})M  && \forall   i\in N,  k\in K
\end{flalign}

\begin{flalign} \label{eqn:earlyStarting7}
	u_{ik} + v_{ik} \leq 1 +\lambda_iM  && \forall  i\in N,  k\in K
\end{flalign}

\begin{flalign} \label{eqn:earlyStarting8}
	u_{ik} + v_{ik} \geq 1 -\lambda_iM  && \forall   i\in N,  k\in K
\end{flalign}

\begin{flalign} \label{eqn:earlyStarting9}
	u_{ik} \leq	\sum_{s=1}^q y_{iks}     && \forall  k \in K,i\in N
\end{flalign}

\begin{flalign} \label{eqn:earlyStarting10}
	v_{ik} \leq	\sum_{s=1}^q y_{iks}     && \forall  k \in K,i\in N
\end{flalign}

\begin{flalign} \label{eqn:xDomMSMTW}
	x_{iik} =0   && \forall  i\in N,  k\in K
\end{flalign}

\begin{flalign} \label{eqn:sDomMSMTW}
	S_{ik} \geq 0   && \forall  i\in N,  k\in K 
\end{flalign}

\begin{flalign}  \label{eqn:ADomMSMTW}
	A_{ik} \geq 0   && \forall  i\in N,  k\in K
\end{flalign}


\begin{flalign} \label{eqn:xDom1MSMTW}
	x_{ijk} \in \{0,1\}   && \forall  i\in N, j\in N, k\in K
\end{flalign}

\begin{flalign} \label{eqn:yDomMSMTW}
	y_{iks} \in \{0,1\} \  && \forall   i\in N, k\in K, s\in S
\end{flalign}

\begin{flalign} \label{eqn:zDomMSMTW}
	z_{il} \in \{0,1\}   && \forall   i\in N, l\in L
\end{flalign}

\begin{flalign} \label{eqn:uDomMSMTW}
	u_{ik} \in \{0,1\}  && \forall  i\in N, k\in K 
\end{flalign}

\begin{flalign} \label{eqn:vDomMSMTW}
	v_{ik} \in \{0,1\}  && \forall  i\in N, k\in K 
\end{flalign}

The objective function aims to minimize  caregivers' waiting times as well as to balance their workload.
Constraints (\ref{eqn:ptVisits1MSMTW}) and (\ref{eqn:ptVisits2MSMTW}) state that each patient will be visited  by a group of  caregivers, which their number depends on the services requested. Constraints (\ref{eqn:leftMSMTW}) and (\ref{eqn:retrunMSMTW}) state that each caregiver who has left the center to visit assigned patients must return to that center. Constraints (\ref{eqn:fluxConservationMSMTW}) impose route continuity (flux conservation) for the patients assigned to a caregiver $k$, while requiring that tours should be constructed rather than open paths. Constraints (\ref{eqn:startingTimeMSMTW}) determine the services operations’ starting time of the patient $j$ with respect to the completion time of service operations of the patient $i$. These constraints enforce that the starting time of services along the route of a caregiver are strictly increasing. In doing so, they also eliminate sub-tours because a return to an already visited patient would violate the start time of the previous visit  \cite{mankowska2014home}.
Constraints  (\ref{eqn:arrivalTime1MSMTW}) and  (\ref{eqn:arrivalTime2MSMTW}) define the arrival time of a caregiver $k$ to the patient $j$. Constraints (\ref{eqn:waitingTimeToZeroMSMTW}) and (\ref{eqn:StartTimeToZeroMSMTW})  initialize the starting and arrival times to zero if the caregiver $k$ will not be affected to the patient $i$. Constraints (\ref{eqn:yDefinitionMSMTW}) define the variables $y_{iks}$,  the caregiver~$k$ is assigned to the patient $i$ if he visits a patient or returns to the \ac{HHC} center after visiting the patient $i$. 
Constraints (\ref{eqn:skillsRequirementMSMTW}) ensure that a qualified caregiver $k$ performs a requested service operation~$s$ to the patient~$i$. Constraints (\ref{eqn:assign1MSMTW}) state that each caregiver provides one service to an assigned patient $i$. Constraints (\ref{eqn:assign2MSMTW}) ensure that for a patient $i$, each service is performed by one caregiver.
Constraints (\ref{eqn:twCaregiver1MSMTW}) and (\ref{eqn:twCaregiver2MSMTW}) ensure that each caregiver is available at the earliest service time $d_k$ at the \ac{HHC} center and must return before the latest service time $e_k$ to that center. Constraints (\ref{eqn:twPatient1MSMTW}) and (\ref{eqn:twPatient2MSMTW}) enforce the respecting of patients’ time windows. Constraints (\ref{eqn:zDefinitionMSMTW}) guarantee that only one time period  is selected from the patient’s availability periods.
Constraints (\ref{eqn:synchronization1MSMTW}) and (\ref{eqn:synchronization2MSMTW}) ensure the synchronization of caregivers' starting time if  requested service operations must be simultaneous ($ \lambda_i =1$) for patient $i$.  These constraints are the refined version of constraints  (\ref{eqn:synchronizationv11MSMTW}) and (\ref{eqn:synchronizationv12MSMTW}), which reduces the number of constraints from  $2nc^2$ to $2nc$. 
\begin{flalign} \label{eqn:synchronizationv11MSMTW}
	S_{ik'}  -   S_{ik} \leq (2 - \sum_{s=1}^q y_{ik's} - \sum_{s=1}^q y_{iks})M  && \forall  i \in N, k \in K, k' \in K
\end{flalign}

\begin{flalign}\label{eqn:synchronizationv12MSMTW}
	S_{ik'}  -   S_{ik} \geq (\sum_{s=1}^q y_{ik's}+ \sum_{s=1}^q y_{iks} -2)M  && \forall  i \in N, k \in K, k' \in K
\end{flalign}

Constraints (\ref{eqn:YdefinitionMSMTW}) define the total working time for each caregiver. Constraints (\ref{eqn:YAbsValue1MSMTW}) and (\ref{eqn:YAbsValue2MSMTW}) define the absolute value $D_k$ between the total working time of the caregiver $k$ and the average working time: $D_k = max \{w_k - \frac{\sum_{k=1}^c {w_k}}{c} , \frac{\sum_{k=1}^c {w_k}}{c}  - w_k \}$.
Constraints (\ref{eqn:earlyStarting1})-(\ref{eqn:earlyStarting8}) guarantee that caregivers start  working as soon as they are available at patients. If the requested services must be simultaneously performed  ($\lambda_i=1$) for a patient $i$, constraints (\ref{eqn:earlyStarting1})-(\ref{eqn:earlyStarting4}) must be verified. $SS_i$ denotes the synchronized starting time.  These constraints are equivalent to the following formula: $SS_{i} = max \{ a_{i\hat{l}},A_{ik_{1}}, A_{ik_{2}},\ldots \} $, which is the maximum between the earliest starting time of  the selected availability period  $\hat{l}$ for the patient $i$ and  caregivers' arrival times $\{ A_{ik_{1}}, A_{ik_{2}},\ldots \} $ at patient $i$. 
Otherwise, the requested service operations are not  needed to be simultaneous ($\lambda_i=0$), constraints (\ref{eqn:earlyStarting5})-(\ref{eqn:earlyStarting8}) ensure that the starting time for  each caregiver equals the maximum between the earliest starting time of the selected availability period  $\hat{l}$ for the patient~$i$ and caregiver $k$ arrival time $A_{ik}$ ($S_{ik} = max \{ a_{i\hat{l}},A_{ik}\} $). Constraints (\ref{eqn:earlyStarting9}) and (\ref{eqn:earlyStarting10}) set  $u_{ik}$ and $v_{ik}$ to zero if the caregiver $k$ will not be affected to the patient $i$. Constraints (\ref{eqn:xDomMSMTW})-(\ref{eqn:vDomMSMTW}) set the domains of the decision variables.

\par

Although the exact methods guarantee to find the optimal solutions, their computation time increases with the size of the problem. Heuristics are more flexible since their computational time could be controlled. A \ac{GVNS}  based heuristic is proposed in the next section to solve the model with large instances.

\subsection{Variable neighborhood search} \label{VNSMSMTW}
In the section we adapt the proposed \ac{GVNS}  in section \ref{VNSMTW} to deal with the new hypothesis. We mention that the same neighborhoods are used with a modification to include multiple services per patient. Therefore, we explain only the changes that has been done.
\subsubsection{Encoding}
We extend the encoding used in the previous chapter to deal with multiple services per patients. Therefor, we will duplicate a patient $i$ as many times as the number of services he requested (see Table~\ref{table:encodingDMSMTW}).
\par
\textbf{Example}: consider 5 patients and 2 caregivers skilled to provide 3 types of services operations, patient $4$ requests two services 2 and 3. A solution will be encoded as follows (see Table~\ref{table:encodingDMSMTW}).

\begin{table}[ !htb]
	
	\caption{Example of solution encoding}
	\centering
	
	\begin{tabular}{ | c | c | c | c | c | c | c| }
		\hline
		Patients   & 1 \color{black}(3)  & 3 \color{black}(1) & \textbf{ 4 \color{black}(3)} & 2 \color{black}(1) & \textbf{4 \color{black} (2)} & 5 \color{black}(2) \\ 
		\hline
		Caregivers  & 1    & 2    & 1    & 2    & 2     &   1\\   
		\hline
	\end{tabular}
	
	\label{table:encodingDMSMTW}
	
\end{table}

The caregiver 1 will visit patients 1, 4 and 5  to provide respectively the services 3, 3 and 2.

\subsubsection{Decoding}

 The decoding method used to compute the objective function is defined according to the following four cases:
 \begin{itemize}
 	\item \textbf{Case 1:  single service and single time window ($SSSTW$)}  for each patient. Since caregivers' routes are independents, their starting times will be iteratively calculated  in the order they appear in the matrix;
 	
 	\item \textbf{Case 2:  single service and multiple time windows ($SSMTW$) }for  each patient. Since caregivers' routes are independents, their starting times will be iteratively calculated  in the order they appear in the matrix. Since each patient has $p$ time windows, the earliest availability period $l\in L$ that minimizes the tardiness of providing the requested service operation $s$ will be chosen for each patient;
 	
 	\item \textbf{Case 3:  multiple services and single time window ($MSSTW$)}. If any patient requests simultaneous services,  caregivers' routes become dependents. Therefore, caregivers' starting times should be synchronized to ensure that are available at the same time to provide those simultaneous services. A synchronized starting time $SS_i$ per patient (in case $\lambda_i=1$) is considered and it is equal to the maximum between the earliest time window $a_{i1}$ for the patient $i$ and  caregivers' arrival times 	$  \{ A_{ik_{1}}, A_{ik_{2}},\ldots \} $. Algorithmically, arrival, starting and waiting times  will be iteratively calculated in the order they appear in the matrix. Each time a caregiver $k$  arrives to  a patient $i$, $SS_i$ is first initialized by $a_{i1}$ and then will take the value $A_{ik}$ if ($ A_{ik}  > SS_i $). These steps are  repeated until $SS_i$ will not change for each patient. In the end, caregivers' starting times at patients will be the synchronized starting time  $ S_{ik}=SS_i$;
 	
 	\item \textbf{Case 4:  multiple services and multiple time windows ($MSMTW$)}. This case generalizes Case~3 to cover multiple time windows per patient.  Since only a single period should be chosen to provide multiple services, this choice may be suitable for one caregiver but not for another.  Therefore, a two-phase strategy is used to select for each patient an availability period to receive requested care services and ensure the synchronization of services operations if are requested to be simultaneous. Phase $I$ aims to heuristically select for each patient a period to receive requested care services. It selects the earliest period that minimizes the tardiness in providing a service operation by the first caregiver arriving at an assigned patient. Then, the other caregivers must provide the requested service operations within that time window.  Phase $II$ will ensure the synchronization of service operations if they are requested to be simultaneous. Since for each patient, a time window is selected among her/his availability periods, Case 4 is reduced to Case 3 and the algorithm described in Case 3 will be applied to ensure the synchronization of multiple services.

 	\textbf{Remark}. The decoding method (cases 3 and 4) could fall into an infinite loop, so simultaneous services would never be synchronized and the solution will be consequently considered infeasible as illustrated by Table \ref{table:encodingDMSMTWInfeasible}. Suppose that patients 4 and 5 request two services 2 and 3. Therefore, services  requested by both patients 4 and 5 cannot be synchronized since caregiver 1 visits patient 4 then patient 5, while caregiver 2 visits patient 5 then patient 4. In the case of hard/fixed time windows, the infeasibility will be considered if after a number of iterations, the synchronized starting time value for a patient $i$ (in the case $\lambda_i=1$) will be greater than the latest service time $b_{il}$. However, in the case of soft/flexible time windows, the solution infeasibility is considered when the number of iterations exceeds the maximum number of iterations denoted by $MaxIterSyn$. 
 	
 	\begin{table}[!htb]
 		
 		\caption{Example of infeasible solution under the synchronization constraint}
 		\centering
 		
 		\begin{tabular}{ | c | c | c | c | c | c | c| c|}
 			\hline
 			Patients       & 1 (3)  & 3 (3) & \textbf{5 (3)}& \textbf{ 4 (3)} & 2 (1) & \textbf{ 4  (2)} & \textbf{5 (2)} \\ 
 			\hline
 			Caregivers   & 1    & 2   &\textbf{2} & \textbf{1 }   & 2    & \textbf{2 }    &  \textbf{ 1}\\   
 			\hline
 		\end{tabular}
 		
 		\label{table:encodingDMSMTWInfeasible}
 		
 	\end{table}

 \end{itemize}

 \subsubsection{Management of infeasible solutions}
 
 While exploring the search space, feasible and infeasible solutions are encountered. The following cases explain how infeasible solutions are managed:

 \begin{itemize}
 	\item \textbf{Patients' time windows and caregivers' duty length}: solutions not respecting constraints (\ref{eqn:twCaregiver2MSMTW}) and (\ref{eqn:twPatient2MSMTW}) are accepted with a penalty cost, which equals to the tardiness of services ($T_{i}$) and caregivers' overtime ($O_k$). The coefficient $\gamma$ is used to  penalize infeasible solutions and to ensure the convergence to feasible ones (see equation (\ref{tarOver})). The higher the penalty coefficient $\gamma$, the faster the convergence;
 	\begin{flalign} \label{tarOver}
 		Z_{GVNS} = Z + \gamma(\sum_{i=1}^n T_{i} + \sum_{k=1}^c O_k) 
 	\end{flalign} 
 	
 	\item \textbf{Skills requirements}: when applying switch and inter-swap neighborhoods, a caregiver could be assigned to  a service operation that he is not qualified  to provide (constraints (\ref{eqn:skillsRequirementMSMTW})). We ensure that only qualified ones are accepted;
 	
 	\item \textbf{Each caregiver must be assigned to a single service for a given patient:} when applying switch and inter-swap neighborhoods, a caregiver could be assigned to provide more than service operation for the same patient  or more than caregiver could be assigned to provide the same service operation for a given patient. This problem could be encountered only when a patient requires multiple services (constraints (\ref{eqn:assign1MSMTW}) and (\ref{eqn:assign2MSMTW})). To avoid this issue, a boolean table ($F_{nc}$) is used and first initialized by $false$. Rows contain patients ($n$ is the size) and columns contains caregivers ($c$ is the size). If patient~$i$ is assigned to caregiver $k$, the $F_{ik}$ is set to $true$.  Before assigning a caregiver to a patient, the table must contain $false$ value. This table is updated as the search space is explored to ensure that constraints  (\ref{eqn:assign1MSMTW}) and (\ref{eqn:assign2MSMTW}) are respected; 
 	
 	\item\textbf{ Constraints \ref{eqn:earlyStarting1}-\ref{eqn:earlyStarting9}:} These constraints are ensured since computing caregivers' starting  times are carried out at the earliest possible time.
 	
 \end{itemize}

 \subsubsection{Initial solution}
 
 The initial solution is generated as follows: 
 \begin{enumerate}
 	\item  Sort patients by the increasing start of their time windows;
 	
 	
 	
 	
 	

 	\item  For each patient in the visiting order found by step 1, assign a qualified caregiver who can arrive the earliest;
 	
 	
 	
 	
 	
 	
 	\item  Compute the objective function $Z$  using the decoding method.
 	
 \end{enumerate}
 Steps 1 and 2 are applied in case of a single time window per patient. In case of multiple time windows, the visiting order (step 1) and caregivers' assignment (step 2) are randomly done since infeasible solutions quickly converge  to feasible ones. However, with a single availability period per patient, penalization (tardiness of services and caregivers' overtime) takes time to converge to zero.
 
  \subsection{Numerical experiments }\label{NEMSMTW}
 \subsubsection{Test instances}

 The test instances have been randomly generated using benchmark instances from \cite{mankowska2014home}. The \ac{HHC} center is supposed to provide six types of services $S= \{1, . . ., 6\}$ to patients. Caregivers are divided into two groups with different skills. Each caregiver in the first group is qualified to provide at most three services, which are randomly selected from the subset $\{1, 2, 3\}$ of $S$. Accordingly, each caregiver in the second group is qualified to provide at most three services from subset $\{4, 5, 6\}$. Caregivers’ duty length is set to 10 hours (600 minutes).
 
 Patients and the \ac{HHC} center are randomly located in the area of $100 \times 100$ (small territory) or $200 \times 200$ (large territory) distance units. Travel times $T_{ij}$ are equal to the Euclidean distance between patients truncated to an integer. Processing times of services operations $t_{is}$ are randomly chosen from  $[10, 20]$ or $[20, 60]$. Patients’ time windows are of length 120 minutes \cite{mankowska2014home} and are placed at random within a daily planning period of 10 hours. For instances with two availability periods, the first period is planted within the first 5 hours. The second period is placed in the interval [5,10] with a deviation of at least 2 hours from the first period.  In the case of three time windows, we randomly place them in the following intervals: [0, 200], [200, 400] and [400, 600]. Each patient requires single or double services, and 30\% of patients are considered requesting double services \cite{mankowska2014home}, which is randomly drawn from $S=\{1,…,6\}$. Only 50\% of double services are considered to be simultaneous  and the other 50\% double services will be provided independently \cite{mankowska2014home}.

 \begin{table}[!htb]
 	\caption{ Tested instances details}
 	\centering
 	\def\arraystretch{1.3}
 	
 	\begin{tabular*}{\textwidth}{@{\extracolsep{\fill}}   lllllll }
 		\hline
 		
 		Set & Subset & Size& $ N$ &$S$  &$K$  &$L$   \\   
 		
 		\hline
 		& A  &	9	& 10	&	10 & [3, 4] & 1   \\   
 		
 		SSSTW   & B &	9	& 25	&	25 & [5, 7] & 1   \\   
 		
 		& C &	9	& 50	&	50 & [10,12] & 1   \\   
 		\hline
 		& D &	9	& 10	&	13 & [3, 4] & 1   \\   
 		
 		MSSTW 	& E &	9	& 25	&	33	 & [5, 7] & 1   \\   
 		
 		& F &	9	& 50	&	65 	& [10, 13] & 1   \\   
 		\hline
 		&G &	9	& 10	&	10 & 3 & 2   \\   
 		
 		SSMTW   & H &	9	& 25	&	25 &  [5, 7] & 2  \\   
 		
 		& I &	9	& 50	&	50 &  10 & 2 \\   
 		\hline
 		
 		& J &	9	& 10	&	13 &  [3, 4] & 2   \\   
 		
 		MSMTW  & K &	9	& 25	&	33 &  [5, 8] & 2  \\   
 		
 		& L &	9	& 50	&	65 &  [10, 14] & 2 \\   
 		\hline

 		Small  & M &	9	& 10	&	[12,16] &  [3, 5] &  [1, 3]\\

 		\hline
 		
 		Large  & N &	4	& [70, 200]	&	[100, 200] &  [20, 40] &  3\\

 		\hline
 		
 	\end{tabular*}

 	\label{table:TIDMSMTW}

 \end{table}

 Four sets of instances are generated, each one is  matched to  a case of the four cases ($SSSTW$, $SSMTW$, $MSSST$ and $MSMTW$) and contains three subsets as well as two other subsets $Small$ and $Large$. Table \ref{table:TIDMSMTW} summarizes the details of tested instances where ’Size’ is the number of instances used in each subset and {'$S$' is the number of total requested services by all patients (number of jobs). For each instance, we start with the minimum number of caregivers as considered in \cite{mankowska2014home}. When we consider large territory and service times in [20, 40], we increase the number of caregivers until a feasible solution is found (see $K$ in Table \ref{table:TIDMSMTW}).
 
 As described above, the same parameters used in \cite{mankowska2014home} are adopted. However, we considered large service time duration $[20, 60]$ and large territory $200 \times 200$, which is more realistic. For each subset (example: subset A), the first three instances (A1,A2 and A3) are used with service times in $[10, 20]$ and patients’ positions in $100 \times 100$; the second three instances (A4, A5 and A6) are used with service times in $[20, 40]$ and patients’ positions in $100 \times 100$; the last three instances (A7, A6 and A9) are used with service times in $[10, 20]$ and patients’ positions in $200 \times 200$. 
  In the set $small$, the proportion of patients requesting multiple services is varied from 20\% to 60\% and the number of time windows is also varied for each patient from 1 to 3. In the set $Large$, we assign to each patient 3 time windows and he can request up to 3 services. 35\% of patients are assumed requesting multiple services, see tables \ref{table:MSMTW1} and \ref{table:MSMTW2}. The parameter N$_k$ with $k \in \{1,2,3\}$, refers to the number of patients who requested $k$ services and therefore  $k$ caregivers should be assigned. Caregivers' qualifications are randomly assigned without sharing them into two groups as was done for the sets $SSSTW$, $SSMTW$, $MSSST$ and $MSMTW$.  
 The weights of the objective functions are set to 0.5 ($\alpha=\beta=0.5$). Several tests were performed to adjust the tuning parameters. Thus, the latter ones are defined as follows:
 \begin{itemize}
 	\item   Shaking order: switch, intra-swap,inter-swap and shift;
 	
 	\item   Local search methods order: inter-swap, intra-swap, shift and switch;
 	
 	\item   Stopping criterion: no improvement over the best solution for 100 iterations ($m_1=100$); 
 	\item  number of times a neighborhood is applied in shaking phase $c+1$ ($m_2=c+1$); 
 	\item  $\gamma=100$ and $MaxIterSyn = 2\times c$.
 	
 \end{itemize}

 \subsubsection{Computational results}

 
 \begin{table*}[!htb]
  
 	\caption{Numerical results for instances with single service and single time window.}
 	\centering
 	\def\arraystretch{1.3}
 	\begin{tabular}{  @{}  l  | llll |lllll}
 		\hline
 		\multicolumn{1}{c|}{\textbf{Set}} & \multicolumn{4 }{c|}{\textbf{CPLEX} }&\multicolumn{5}{c}{ \textbf{\ac{GVNS}} (10 runs)} \\ 
 		\hline
 		
 		\textbf{SSSTW}  &\textbf{LB} & \textbf{Z} & \textbf{GAP}&\textbf{CPU} & \textbf{Best} & \textbf{Worst}& A\textbf{verage} &\textbf{ Gap} & \textbf{CPU }   \\

 		\hline			
 		
 		A1	&	341.50	&	\textbf{341.50}	&	0.00\%	&	2.47	&	\textbf{341.50}	&   \textbf{341.50}	&	\textbf{341.50}	&	0.00\%	&	$<1$	\\
 		A2	&	211.50	&	\textbf{211.50}	&	0.00\%	&	1.98	&	\textbf{211.50}	&	\textbf{211.50}	&	\textbf{211.50}	&	0.00\%	&	$<1$	\\
 		A3	&	290.00	&	\textbf{290.00}	&	0.00\%	&	2.03	&	\textbf{290.00}	&	\textbf{290.00}	&	\textbf{290.00}	&	0.00\%	&	$<1$	\\
 		A4	&	206.00	&	\textbf{206.00}	&	0.00\%	&	2.09	&	\textbf{206.00}	&	\textbf{206.00}	&	\textbf{206.00}	&	0.00\%	&	$<1$	\\
 		A5	&	331.00	&	\textbf{331.00}	&	0.00\%	&	1.88	&	\textbf{331.00}	&	\textbf{331.00}	&	\textbf{331.00}	&	0.00\%	&	$<1$	\\
 		A6	&	331.00	&	\textbf{331.00}	&	0.00\%	&	2.07	&	\textbf{331.00}	&	\textbf{331.00}	&	\textbf{331.00}	&	0.00\%	&	$<1$	\\
 		A7	&	246.00	&	\textbf{246.00}	&	0.00\%	&	2.83	&	\textbf{246.00}	&	\textbf{246.00}	&	\textbf{246.00}	&	0.00\%	&	$<1$	\\
 		A8	&	345.00	&	\textbf{345.00}	&	0.00\%	&	2.71	&	\textbf{345.00}	&	\textbf{345.00}	&	\textbf{345.00}	&	0.00\%	&	$<1$	\\
 		A9	&	225.00	&	\textbf{225.00}	&	0.00\%	&	2.16	&	\textbf{225.00}	&	\textbf{225.00}	&	\textbf{225.00}	&	0.00\%	&	$<1$	\\
 		\hline	
 		B1	&	195.00	&	\textbf{195.00}	&	0.00\%	&	1787	&	\textbf{195.00}	&	207.00	&	201.10	&	3.03\%	&	5.65	\\
 		B2	&	211.50	&	\textbf{211.50}	&	0.00\%	&	149.16	&	\textbf{211.50}	&	217.50	&	214.20	&	1.26\%	&	5.85	\\
 		B3	&	308.50	&	\textbf{308.50}	&	0.00\%	&	697.50	&	\textbf{308.50}	&	317.50	&	311.90	&	1.09\%	&	5.21	\\
 		B4	&	340.50	&	\textbf{340.50}	&	0.00\%	&	24.41	&	\textbf{340.50}	&	\textbf{340.50}	&	\textbf{340.50}	&	0.00\%	&	4.04	\\
 		B5	&	277.50	&	\textbf{277.50}	&	0.00\%	&	12.73	&	\textbf{277.50}	&	\textbf{277.50}	&	\textbf{277.50}	&	0.00\%	&	5.75	\\
 		B6	&	345.50	&	\textbf{345.50}	&	0.00\%	&	466.13	&	\textbf{345.50}	&	367.00	&	353.25	&	2.19\%	&	6.23	\\
 		B7	&	437.00	&	\textbf{437.00}	&	0.00\%	&	1274	&	\textbf{437.00}	&	\textbf{437.00}	&	\textbf{437.00}	&	0.00\%	&	5.79	\\
 		B8	&	125.00	&	\textbf{125.00}	&	0.00\%	&	241.73	&	\textbf{125.00}	&	146.50	&	129.90	&	3.77\%	&	7.57	\\
 		B9	&	184.50	&	\textbf{184.50}	&	0.00\%	&	305.84	&	\textbf{184.50}	&	254.00	&	195.80	&	5.77\%	&	7.7	\\
 		
 		\hline	
 		C1	&	156.26	&	519.50	&	69.92\%	&	7200	&	461.50	&	471.45	&	483.00	&	67.65\%	&	107.49	\\
 		C2	&	58.78	&	672.00	&	91.25\%	&	7200	&	545.50	&	606.00	&	573.35	&	89.75\%	&	119.79	\\
 		C3	&	18.50	&	671.00	&	97.24\%	&	7200	&	575.00	&	602.00	&	584.70	&	96.84\%	&	116.56	\\
 		C4	&	46.50	&	665.50	&	93.01\%	&	7200	&	609.50	&	665.00	&	649.70	&	92.84\%	&	84.99	\\
 		C5	&	2.44	&	221.00	&	98.90\%	&	7200	&	221.00	&	280.50	&	235.90	&	98.97\%	&	91.27	\\
 		C6	&	97.57	&	369.00	&	73.56\%	&	7200	&	334.00	&	411.00	&	365.75	&	73.32\%	&	100.86	\\
 		C7	&	30.76	&	552.50	&	94.43\%	&	7200	&	514.00	&	574.00	&	551.65	&	94.42\%	&	112.14	\\
 		C8	&	0.00	&	179.00	&	100\%	&	7200	&	152.00	&	216.00	&	188.40	&	100\%	&	115.77	\\
 		C9	&	60.89	&	136.50	&	55.39\%	&	7200	&	132.00	&	239.00	&	196.65	&	69.04\%	&	117.36	\\

 		\hline

 	\end{tabular}

 	\label{table:SSSTW}

 \end{table*}

 \begin{table*}[!htb]

 	\caption{Numerical results for instances with multiple services and single time window.}
 	\centering
 	\def\arraystretch{1.3}
 	\begin{tabular}{@{}  l   | llll |lllll}
 		\hline
 		\multicolumn{1}{c|}{\textbf{Set}} & \multicolumn{4 }{c|}{\textbf{CPLEX} }&\multicolumn{5}{c}{ \textbf{\ac{GVNS}} (10 runs)} \\ 
 		\hline
 		
 		\textbf{MSSTW}  &\textbf{LB} & \textbf{Z} & \textbf{GAP}&\textbf{CPU} & \textbf{Best} & \textbf{Worst}& A\textbf{verage} &\textbf{ Gap} & \textbf{CPU }   \\   
 		\hline
 		
 		D1	&	159.00	&	\textbf{159.00}	&	0.00\%	&	2.32	&	\textbf{159.00}	&	\textbf{159.00}	&	\textbf{159.00}	&	0.00\%	&	$<1$	\\
 		D2	&	150.50	&	\textbf{150.50}	&	0.00\%	&	1.91	&	\textbf{150.50}	&	\textbf{150.50}	&	\textbf{150.50}	&	0.00\%	&	$<1$	\\
 		D3	&	208.00	&	\textbf{208.00}	&	0.00\%	&	3.16	&	\textbf{208.00}	&	214.00	&	209.55	&	0.74\%	&	$<1$\\
 		D4	&	115.00	&	\textbf{115.00}	&	0.00\%	&	3.39	&	\textbf{115.00}	&	\textbf{115.00}	&	\textbf{115.00}	&	0.00\%	&	$<1$\\
 		D5	&	275.50	&	\textbf{275.50}	&	0.00\%	&	2.44	&	\textbf{275.50}	&	\textbf{275.50}	&	\textbf{275.50}	&	0.00\%	&	$<1$\\
 		D6	&	218.50	&	\textbf{218.50}	&	0.00\%	&	2.96	&	\textbf{218.50}	&	\textbf{218.50}	&	\textbf{218.50}	&	0.00\%	&	$<1$	\\
 		D7	&	227.00	&	\textbf{227.00}	&	0.00\%	&	2.09	&	\textbf{227.00}	&	\textbf{227.00}	&	\textbf{227.00}	&	0.00\%	&	$<1$\\
 		D8	&	343.50	&	\textbf{343.50}	&	0.00\%	&	2.06	&	\textbf{343.50}	&	\textbf{343.50}	&	\textbf{343.50}	&	0.00\%	&	$<1$	\\
 		D9	&	113.50	&	\textbf{113.50}	&	0.00\%	&	2.1	&	\textbf{113.50}	&	\textbf{113.50}	&	\textbf{113.50}	&	0.00\%	&	$<1$	\\
 		
 		\hline
 		E1	&	42.50	&	\textbf{42.50}	&	0.00\%	&	2822	&\textbf{42.50}	&	58.00	&	48.75	&	12.82\%	&	37.58	\\
 		E2	&	42.75	&	67.00	&	36.19\%	&	7200	&	67.00	&	97.50	&	79.70	&	46.36\%	&	30.80	\\
 		E3	&	81.50	&	147.00	&	44.56\%	&	7200	&	147.00	&	163.00	&	155.00	&	47.42\%	&	29.81	\\
 		E4	&	186.25	&	249.00	&	25.20\%	&	7200	&	249.00	&	280.50	&	260.20	&	28.42\%	&	21.13	\\
 		E5	&	283.00	&	\textbf{283.00}	&	0.00\%	&	317.86	&	\textbf{283.00}	&	285.50	&	283.70	&	0.25\%	&	34.79	\\
 		E6	&	229.83	&	336.00	&	31.60\%	&	7200	&	336.00	&	843.50	&	643.35	&	64.28\%	&	32.17	\\
 		E7	&	82.00	&	\textbf{82.00}	&	0.00\%	&	167.25	&	\textbf{82.00}	&	155.00	&	103.00	&	20.39\%	&	26.86	\\
 		E8	&	195.50	&	\textbf{195.50}	&	0.00\%	&	367.14	&	\textbf{195.50}	&	373.00	&	248.45	&	21.31\%	&	45.25	\\
 		E9	&	80.50	&	\textbf{80.50}	&	0.00\%	&	551.05	&	\textbf{80.50}	&	152.00	&	123.35	&	34.74\%	&	35.12	\\
 		
 		\hline
 		F1	&	0.00	&		&		&	7200	&	202.50	&	237.00	&	220.40	&	100\%	&	537.73	\\
 		F2	&	0.00	&		&		&	7200	&	335.50	&	404.50	&	367.00	&	100\%	&	516.14	\\
 		F3	&	28.50	&		&		&	7200	&	530.50	&	580.00	&	546.70	&	94.79\%	&	475.02	\\
 		F4	&	0.00	&	495.50	&	100\%	&	7200	&	365.50	&	422.00	&	385.70	&	100\%	&	480.81	\\
 		F5	&	0.00	&	368.00	&	100\%	&	7200	&	295.50	&	325.00	&	310.50	&	100\%	&	504.58	\\
 		F6	&	22.00	&		&		&	7200	&	313.50	&	373.00	&	344.30	&	93.61\%	&	575.58	\\
 		F7	&	0.00	&		&		&	7200	&	242.00	&	289.00	&	266.90	&	100\%	&	815.32	\\
 		F8	&	24.00	&		&		&	7200	&	281.50	&	419.00	&	339.50	&	92.93\%	&	1015.05	\\
 		F9	&	13.50	&		&		&	7200	&	219.00	&	337.50	&	265.00	&	94.91\%	&	731.01	\\
 		
 		\hline

 	\end{tabular}

 	\label{table:MSSTW}

 \end{table*}

 \begin{table*}[!htb]

 	\caption{Numerical results for instances with single service and multiple time windows.}
 	\centering
 	\def\arraystretch{1.3}
 	\begin{tabular}{  @{}  l   | llll | lllll}
 		\hline
 		\multicolumn{1}{c|}{\textbf{Set}} & \multicolumn{4 }{c|}{\textbf{CPLEX} }&\multicolumn{5}{c}{ \textbf{\ac{GVNS}} (10 runs)} \\ 
 		\hline
 		
 		\textbf{SSMTW}  &\textbf{LB} & \textbf{Z} & \textbf{GAP}&\textbf{CPU} & \textbf{Best} & \textbf{Worst}& A\textbf{verage} &\textbf{ Gap} & \textbf{CPU }   \\   
 		\hline
 		
 		G1	&	159.50	&	\textbf{159.50}	&	0.00\%	&	3.00	&	\textbf{159.50}	&	164.50	&	160.00	&	0.31\%	&	$<1$\\
 		G2	&	51.50	&	\textbf{51.50}	&	0.00\%	&	2.40	&	\textbf{51.50}	&	61.00	&	53.40	&	3.56\%	&	$<1$	\\
 		G3	&	64.00	&	\textbf{64.00}	&	0.00\%	&	2.56	&	\textbf{64.00}	&	\textbf{64.00}	&	\textbf{64.00}	&	0.00\%	&	$<1$	\\
 		G4	&	181.50	&	\textbf{181.50}	&	0.00\%	&	2.86	&	\textbf{181.50}	&	188.50	&	184.30	&	1.52\%	&	$<1$	\\
 		G5	&	170.00	&	\textbf{170.00}	&	0.00\%	&	2.49	&	\textbf{170.00}	&	203.00	&	176.45	&	3.66\%	&	$<1$	\\
 		G6	&	129.50	&	\textbf{129.50}	&	0.00\%	&	2.42	&	\textbf{129.50}	&	154.50	&	132.00	&	1.89\%	&	$<1$	\\
 		G7	&	57.00	&	\textbf{57.00}	&	0.00\%	&	2.70	&	\textbf{57.00}	&	\textbf{57.00}	&	\textbf{57.00}	&	0.00\%	&	$<1$	\\
 		G8	&	163.00	&	\textbf{163.00}	&	0.00\%	&	2.33	&	\textbf{163.00}	&	\textbf{163.00}	&	\textbf{163.00}	&	0.00\%	&	$<1$	\\
 		G9	&	92.50	&	\textbf{92.50}	&	0.00\%	&	4.83	&	\textbf{92.50}	&	107.50	&	95.50	&	3.14\%	&	$<1$	\\
 		
 		\hline
 		H1	&	79.33	&	92.00	&	13.77\%	&	7200	&	92.00	&	132.00	&	116.40	&	31.85\%	&	13.79	\\
 		H2	&	39.07	&	121.00	&	67.71\%	&	7200	&	111.00	&	121.00	&	114.75	&	65.95\%	&	18.66	\\
 		H3	&	0.00	&	100	&	100\%	&	7200	&	88.00	&	134.00	&	109.10	&	100\%	&	14.22	\\
 		H4	&	239.50	&	\textbf{239.50}	&	0.00\%	&	6044	&	\textbf{239.50}	&	331.00	&	278.85	&	14.11\%	&	15.99	\\
 		H5	&	87.50	&	\textbf{87.50}	&	0.00\%	&	1157  &	\textbf{87.50}	&	144.50	&	98.55	&	11.21\%	&	16.56	\\
 		H6	&	101.50	&	205.50	&	50.61\%	&	7200	&	202.00	&	231.50	&	218.95	&	53.64\%	&	22.46	\\
 		H7	&	115.50	&	\textbf{115.50}	&	0.00\%	&	433.53	&	\textbf{115.50}	&	162.50	&	131.50	&	12.17\%	&	16.34	\\
 		H8	&	150.12	&	196.00	&	23.41\%	&	7200	&	184.00	&	220.50	&	200.40	&	25.09\%	&	22.95	\\
 		H9	&	6.50	&	46.50	&	86.02\%	&	7200	&	46.50	&	118.00	&	82.85	&	92.15\%	&	24.15	\\
 		\hline
 		
 		I1	&	0.00	&		&		&	7200	&	162.00	&	211.50	&	179.55	&	100\%	&	316.99	\\
 		I2	&	0.00	&		&		&	7200	&	126.50	&	155.50	&	140.15	&	100\%	&	200.36	\\
 		I3	&	0.00	&	258.50	&	100\%	&	7200	&	90.00	&	175.50	&	135.45	&	100\%	&	306.37	\\
 		I4	&	0.00	&		&		&	7200	&	336.00	&	437.50	&	369.35	&	100\%	&	329.10	\\
 		I5	&	0.00	&		&		&	7200	&	141.50	&	240.00	&	183.45	&	100\%	&	232.57	\\
 		I6	&	0.00	&		&		&	7200	&	239.50	&	330.50	&	274.20	&	100\%	&	194.80	\\
 		I7	&	0.00	&		&		&	7200	&	74.00	&	179.00	&	127.15	&	100\%	&	300.17	\\
 		I8	&	0.00	&		&		&	7200	&	108.50	&	1035.00	&	276.80	&	100\%	&	351.80	\\
 		I9	&	0.00	&		&		&	7200	&	64.00	&	214.00	&	136.10	&	100\%	&	254.03	\\

 		\hline	
 	\end{tabular}
 	
 	\label{table:SSMTW}
 \end{table*}

 \begin{table*}[!htb]
 	
 	\caption{Numerical results for instances with multiple services and multiple time windows.}
 	\centering
 	\def\arraystretch{1.3}
 	\begin{tabular}{ @{}   l   | llll| lllll}
 		\hline
 		\multicolumn{1}{c|}{\textbf{Set}} & \multicolumn{4 }{c|}{\textbf{CPLEX} }&\multicolumn{5}{c}{ \textbf{\ac{GVNS}} (10 runs)} \\ 
 		\hline
 		
 		\textbf{MSMTW}  &\textbf{LB} & \textbf{Z} & \textbf{GAP}&\textbf{CPU} & \textbf{Best} & \textbf{Worst}& A\textbf{verage} &\textbf{ Gap} & \textbf{CPU }   \\   
 		
 		\hline
 		
 		J1	&	161.00	&	\textbf{161.00}	&	0.00\%	&	3.71	&	\textbf{161.00}	&	242.00	&	199.90	&	19.46\%	&	1.19	\\
 		J2	&	130.50	&	\textbf{130.50}	&	0.00\%	&	12.89	&	\textbf{130.50}	&	239.00	&	168.00	&	22.32\%	&	1.05	\\
 		J3	&	102.50	&	\textbf{102.50}	&	0.00\%	&	6.01	&	\textbf{102.50}	&	192.00	&	152.90	&	32.96\%	&	1.13	\\
 		J4	&	311.00	&	\textbf{311.00}	&	0.00\%	&	6.00	&	\textbf{311.00}	&	361.00	&	341.10	&	8.82\%	&	1.30	\\
 		J5	&	160.50	&	\textbf{160.50}	&	0.00\%	&	13.60	&	\textbf{160.50}	&	192.00	&	168.40	&	4.69\%	&	1.65	\\
 		J6	&	160.00	&	\textbf{160.00}	&	0.00\%	&	10.40	&	\textbf{160.00}	&	213.00	&	172.55	&	7.27\%	&	1.84	\\
 		J7	&	246.50	&	\textbf{246.50}	&	0.00\%	&	5.65	&	\textbf{246.50}	&	278.50	&	254.90	&	3.30\%	&	1.31	\\
 		J8	&	57.50	&	\textbf{57.50}	&	0.00\%	&	6.61	&	\textbf{57.50}	&	145.50	&	66.30	&	13.27\%	&	1.94	\\
 		J9	&	165.00	&	\textbf{165.00}	&	0.00\%	&	32.49	&	\textbf{165.00}	&	227.00	&	192.55	&	14.31\%	&	1.75	\\

 		\hline

 		K1	&	0.00	&	177.00	&	100\%	&	7200	&	130.00	&	155.00	&	157.70	&	100\%	&	81.69	\\
 		K2	&	37.50	&	82.50	&	54.55\%	&	7200	&	121.50	&	321.50	&	194.65	&	80.73\%	&	76.40	\\
 		K3	&	0.00	&	118.50	&	100\%	&	7200	&	156.50	&	248.00	&	176.65	&	100\%	&	85.54	\\
 		K4	&	0.00	&	381.00	&	100\%	&	7200	&	414.50	&	549.00	&	491.60	&	100\%	&	91.87	\\
 		K5	&	0.00	&	219.50	&	100\%	&	7200	&	198.50	&	362.50	&	299.35	&	100\%	&	66.85	\\
 		K6	&	0.00	&	286.00	&	100\%	&	7200	&	245.00	&	420.50	&	314.45	&	100\%	&	82.84	\\
 		K7	&	0.00	&	294.00	&	100\%	&	7200	&	289.00	&	410.50	&	329.05	&	100\%	&	89.73	\\
 		K8	&	4.00	&	76.50	&	94.77\%	&	7200	&	66.00	&	240.00	&	161.20	&	97.52\%	&	97.75	\\
 		K9	&	0.00	&	201.50	&	100\%	&	7200	&	126.00	&	202.50	&	176.55	&	100\%	&	83.08	\\
 		
 		\hline
 		
 		L1	&	0.00	&		&		&	7200	&	205.50	&	328.50	&	240.30	&	100\%	&	1526.36	\\
 		L2	&	0.00	&		&		&	7200	&	108.50	&	227.50	&	155.55	&	100\%	&	1280.94	\\
 		L3	&	0.00	&		&		&	7200	&	127.50	&	219.00	&	177.75	&	100\%	&	2066.75	\\
 		L4	&	0.00	&		&		&	7200	&	409.00	&	566.00	&	462.05	&	100\%	&	1268.54	\\
 		L5	&	0.00	&		&		&	7200	&	237.50	&	424.50	&	313.60	&	100\%	&	1427.12	\\
 		L6	&	0.00	&		&		&	7200	&	490.00	&	803.50	&	661.60	&	100\%	&	1137.31	\\
 		L7	&	0.00	&		&		&	7200	&	256.50	&	717.00	&	425.50	&	100\%	&	1558.86	\\
 		L8	&	0.00	&		&		&	7200	&	214.50	&	384.00	&	284.50	&	100\%	&	1350.42	\\
 		L9	&	0.00	&		&		&	7200	&	251.50	&	488.00	&	326.50	&	100\%	&	1347.50	\\
 		
 		\hline	
 	\end{tabular}

 	\label{table:MSMTW}

 \end{table*}

 \begin{table*}[!htb]

 	\caption{Numerical results for small instances with varying the number of time windows from 1 to 3 and the percentage of patients requresting multiple services from 20\% to 60\%.}
 	\centering
 	\def\arraystretch{1.3}
 	 \hspace*{-1cm}
 	\begin{tabular}{  @{}  lllll   | lll | lllll}
 		\hline
 		\multicolumn{5}{c|}{\textbf{Set}} & \multicolumn{3}{c|}{\textbf{CPLEX} }&\multicolumn{5}{c}{ \textbf{\ac{GVNS}} (10 runs)} \\ 
 		\hline
 		
 		\textbf{Small}  &\textbf{N} &\textbf{K} &\textbf{S} &\textbf{L} &\textbf{LB} & \textbf{Z} & \textbf{CPU} & \textbf{Best} & \textbf{Worst}& A\textbf{verage} &\textbf{ Gap} & \textbf{CPU }   \\   
 		\hline

 		M1	&		&		&		&	1	&	223.00	&	\textbf{223.00}	&	2.82	&	\textbf{223.00}	&	\textbf{223.00}	&	\textbf{223.00}	&	0.00\%	&	$<1$	\\
 		M2	&	10	&	3	&	12	&	2	&	19.00	&	\textbf{19.00}	&	6.84	&	\textbf{19.00}	&	29.50	&	20.90	&	9.09\%	&	1.05	\\
 		M3	&		&		&		&	3	&	8.50	&	\textbf{8.50}	&	8.75	&	\textbf{8.50}	&	28.5	&	16.35	&	48.01\%	&	1.16	\\
 		
 		\hline
 		
 		M4	&	&		&		&	1	&	288.50	&	\textbf{288.50}	&	2.15	&	\textbf{288.50}	&	289.00	&	288.55	&	0.02\%	&	1.23	\\
 		M5	&	10	&	3	&	14	&	2	&	146.50	&	\textbf{146.50}	&	14.94	&	\textbf{146.50}	&	199.00	&	167.20	&	12.38\%	&	1.88	\\
 		M6	&	&		&		&	3	&	51.00	&	\textbf{51.00}	&	25.91	&	\textbf{51.00}	&	102.50	&	80.55	&	36.69\%	&	1.69	\\
 		
 		\hline
 		
 		M7	&		&		&		&	1	&	424.50	&	\textbf{424.50}	&	201.98	&	\textbf{424.50}	&	429.50	&	425.50	&	0.24\%	&	2.85	\\
 		M8	&	10	&	5	&	16	&	2	&	51.00	&	\textbf{51.00}	&	436.92	&	\textbf{51.00}	&	136.50	&	81.15	&	37.15\%	&	6.87	\\
 		M9	&		&		&		&	3	&	29.50	&	\textbf{29.50}	&	3202	&	\textbf{29.50}	&	51.00	&	35.95	&	17.94\%	&	7.74	\\

 		\hline	
 	\end{tabular}
 	
 	\label{table:MSMTW1}
 \end{table*}

 \begin{table*}[!htb]

 	\caption{Numerical results for large  instances with multiple services and multiple time windows.}
 	\centering
 	\def\arraystretch{1.3}
 	\begin{tabular}{ @{}   llllllll   | lll| lll}
 		\hline
 		\multicolumn{8}{c|}{\textbf{Set}} & \multicolumn{3 }{c|}{\textbf{CPLEX} }&\multicolumn{3}{c}{ \textbf{\ac{GVNS}} (1 run)} \\ 
 		\hline
 		
 		\textbf{Large}  &\textbf{S} & \textbf{K} & \textbf{N}&\textbf{N$_1$} & \textbf{N$_2$} & \textbf{N$_3$}& \textbf{L}&  \textbf{LB} &  \textbf{Z}& \textbf{CPU }  & \textbf{Z} &\textbf{ Gap} & \textbf{CPU }  \\   
 		
 		\hline
 		
 		N100s	&100	&20		&70		&46		&18	&6	&3	&0	&-	&7200	&179.00	&100\%	&1816.77\\
 		N100p 	&100	&20		&100	&100	&0	&0	&3	&0	&-	&7200	&95.50	&100\%	&342.40\\
 		N200s	&200	&40		&140	&92		&36	&12	&3	&0	&-	&7200	&488.00	&100\%	&8613.52\\
 		N200p	&200	&40		&200	&200	&0	&0	&3	&0	&-	&7200	&212.50	&100\%	&3314.38\\
 		
 		\hline
 	\end{tabular}

 	\label{table:MSMTW2}

 \end{table*}
 
 
 Instances are generated as described above  and solved within a time limit of 2 hours. 
 For the set $SSSTW$, the instances of the subsets $A$ and $B$ are optimally solved.  The CPLEX solver failed to reach optimal solutions for instances of the subset $C$, gaps of reached solutions are in the interval $[55.39\%, 100\%]$ (see Table \ref{table:SSSTW}). For the set $MSSTW$, the instances of the subset $D$ and the instances $E1$, $E5$,  $E7$, $E8$ and $E9$ are optimally solved.  The CPLEX solver failed to reach optimal solutions for the instances $E2, E3, E4, E6,  F4$ and $ F5$, gaps of reached solutions are in the interval $[25.20\%, 100\%]$.  A feasible solution has not been found for the instances $F1, F2, F3, F6, F7, F8$ and $ F9$  (see Table~\ref{table:MSSTW}).
 
 For the set $SSMTW$, the instances of the subset $G$ as well as the instances $H4, H5$ and $H7$ are optimally solved. The CPLEX solver was not able to reach optimal solutions for the remaining instances of the subset $H$ and the instance $I3$, the gaps of reached solutions are in the interval $[13.77\%, 100\%]$. A feasible solution was not found for the subset $I$ except the instance $I3$ (see Table \ref{table:SSMTW}).  For the set $MSMTW$, the instances of the subset $J$ are optimally solved. The CPLEX solver failed to reach optimal solutions  for the subset $K$, gaps of reached solutions are in the interval $[54.55\%, 100\%]$.  A feasible solution was not found for the subset $L$  (see Table \ref{table:MSMTW}).

 	For the set $Small$, all instances are optimally solved (see Table \ref{table:MSMTW1}) while the CPLEX solver failed to find feasible solutions for the instances of the set $Large$ (see Table \ref{table:MSMTW2}). Table \ref{table:MSMTW1} shows that when increasing the number of time windows, this leads to an increasing of the CPU running time as well as to a minimization of the objective function $Z$. As the percentage of patients requiring multiple services and the number of time windows per patient increases, the CPU running time increases. For the same number of patients (10), the CPU running time jumped from 2.82 seconds (instance $M1$ with 12 jobs and a single time window per patient) to 3202 seconds (instance $M9$ with 16 jobs and 3 time windows per patient).

 \par
 To sum up, for the  set $SSSTW,$ 18 instances are optimally solved while for the remaining instances, only feasible solutions are provided. 
 For the set $MSSTW,$ 14 instances are optimally solved, a feasible solution is found for 6 instances and the CPLEX solver failed to find a feasible solution for 7 instances.  
 For the set $SSMTW,$ 12 instances are optimally solved, a feasible solution is found for 7 instances and the CPLEX solver could not find a feasible solution for 8 instances. 
 For the set $MSMTW,$ 9 instances are optimally solved, a feasible solution is found for 9 instances and the CPLEX solver could not find a feasible solution for 9 instances
For the  subset $Small$, all instances are solved to optimality  and no feasible solution is found for the instances of the subset $Large$ (see Table \ref{table:slsCmp}, \textit{"CPLEX$_{Better}$"} presents the number of the CPLEX solutions that are better are  than those found by \ac{GVNS} ).
 
 Instances become more complicated when multiple services, synchronization and multiple time windows constraints are considered. Multiple care services increase the number of patients' visits, which in turn increases the number of caregivers involved in the scheduling. The synchronization of visits makes the problem more complex because caregivers should be available at the same time to provide simultaneous services. Multiple time windows constraints are hard to satisfy since their complexity is  exponential. Indeed, if there are $n$ patients with $p$ availability periods for each, we have $p^n$ possibilities to select for each patient an availability period to receive care services. 
 
 \begin{table}[ !htb]
 	
 	\caption{ Comparison between solutions found by the CPLEX solver and \ac{GVNS} }
 	\centering
 	\def\arraystretch{1.3}
 	\begin{threeparttable}
 		\begin{tabular*}{\textwidth}{@{\extracolsep{\fill}} l| lll| lll|ll}
 			
 			\hline
 			& \multicolumn{3}{c|}{\textbf{Cplex solutions}} &\multicolumn{3}{c|}{\textbf{\ac{GVNS} solutions}}&\multicolumn{2}{c}{\textbf{\ac{GVNS} vs Cplex}}\\  
 			Set & O  & F   & NS& O   & F & NS &   Cplex$_{Better}$ & \ac{GVNS} $_{Better}$\\ 
 			\hline
 			SSSTW & 18 & 9 & 0  &18  &9 & 0& 0 & 8  \\
 			\hline
 			MSSTW & 14 & 6 & 7  & 14&13 &0& 0 & 9 \\
 			\hline
 			SSMTW & 12 & 7 & 8 &12 &15 &0&0 & 13  \\
 			\hline
 			MSMTW & 9  &  9 &  9 &9 &18 &0& 3 &15 \\
 			\hline
 			Small & 9  &  0 &  0 &9 &0 &0& 0& 0\\
 			\hline
 			large & 0  &  0 &  4 &0 &4 &0& 0 &4 \\
 			\hline
 			Total & 62  &  31 &  28 &62 &59 &0&  3 &  49  \\
 			
 			\hline
 			
 		\end{tabular*}
 		\begin{tablenotes}
 			
 			\small
 			\item \textbf{Abbreviations:} \textbf{O} (Optimal), \textbf{F} (Feasible), and  \textbf{NS} (No solution).
 			
 		\end{tablenotes}
 	\end{threeparttable}
 	\label{table:slsCmp}
 	
 \end{table}

 \par 
 
 Although exact methods give the optimal solution, their computation time  increases monotonically with the size of the problem. The \ac{GVNS}  based heuristic combined with the proposed strategy, which aims to select a time window for each patient to receive visits and ensures the synchronization of simultaneous services, could solve all instances in a very short computational CPU running time.
 \ac{GVNS}  could find all optimal solutions in a short computing time at least once in 10 runs. In the case where the CPLEX solver fails to reach optimal solutions, \ac{GVNS}  showed the best performance in terms of the solution quality on the 49 instances. In Table \ref{table:slsCmp}, "\ac{GVNS}$_{Better}$" expresses the number of \ac{GVNS}  solutions that are better than those found by the CPLEX solver. 
 The worst CPU time of instances solved to optimality is on average  4.525 seconds (45.25 for 10 runs), which is elapsed to solve instance $E8$. Contrariwise, CPU times of  CPLEX optimal solutions are varying from 1 to 6044 seconds. \ac{GVNS}  was able to reach better solutions than most of the feasible solutions found by the CPLEX solver (see Tables \ref{table:SSSTW}, \ref{table:MSSTW}, \ref{table:SSMTW} and \ref{table:MSMTW}). The CPLEX solver failed to find feasible solutions for the instances of 50 patients with double services and/or double time windows except 3 from 27 instances (see the sets $F, I$ and $L$ in Tables \ref{table:MSSTW}, \ref{table:SSMTW} and \ref{table:MSMTW}). \ac{GVNS}  was able to solve up to 200 jobs with 3 times windows per patients and 3 services for some patients (see Table \ref{table:MSMTW2}). 
 
 \section {Conclusion} \label{DMConclusion}

 The \ac{HHC} companies aim to meet two requirements for the management of their medical service operations. They have to meet patient demands at the right time and coordinate caregivers’ activities in order to reduce operating costs.  From this perspective, this chapter addressed the  \ac{HHCRSP} with multiple time windows and multiple synchronized services (HHCRSP-MTW-MSS). 

First we dealt with the \ac{HHCRSP} with soft/flexible multiple time windows.  The objective aims to minimize the earliness and tardiness of providing services as well as to minimize caregivers’ waiting times. We showed through a comparison that instances with multiple time windows are better optimized than those with single availability period. Then, we generalized the model to allow patients to request multiple services with a possible synchronization. The objective is to optimize caregivers'  waiting time as well as to balance their workload in order to establish a daily planning. 

Exact methods are powerless to solve small instances  in a reasonable computation time due to the model that extends an NP-hard problem (VRP) as well as the consideration of multiple time windows and synchronization constraints. To overcome this limitation, a  \ac{GVNS}  based heuristic is proposed and combined with a strategy that aims to select a time window for each patient and  ensure the synchronization of simultaneous services in order to solve the model. Several tests were carried out on a randomly generated data to assess \ac{GVNS}  performance in terms of solution quality and efficiency. It is worth mentioning that unlike exact methods, the tests proved the high performance of the \ac{GVNS}  to deal with large instances in a very short computational time. 

The objectives aimed at improving the performance of \ac{HHC} companies are not only limited to those treated in this chapter. An interesting future work is to extend the proposed model to a multi-objective problem and develop approaches to handle the multi-criteria case. A second perspective is to extend the proposed model to deal with stochastic travel and service times since deterministic models ignore uncertainties that could happen. Also, in the event that a new risky service needs to be provided or a requested service is canceled, it would be interesting to propose solution approaches that avoid revising the overall planning.

   \chapter{Stochastic home health care routing and scheduling problem}\label{ChapterSM} 
 
 \section {Introduction}
 
In this chapter, we study the  \ac{HHCRSP} with particular interest in the uncertainty of parameters. Demand, travel and service times are the main parameters that are exposed to the uncertainty. We consider the uncertainty of travel and service times because these two parameters are critical elements in the planning, any change could affect the overall planning and service quality would be poor or even risky. The uncertainty in travel times may be due to common factors such as varying road conditions, driving skills and weather conditions \cite{shi2019robust}. However, the service time is not always fixed as estimated due to practical reasons, such as diagnosing time and parking situations \cite{shi2019robust}. Several models and methods have been proposed in the literature to deal with the \ac{HHCRSP} but most of them are deterministic and generally less robust. The predefined schedule should be revised for any change in practical situations. Otherwise, there will probably be delays in the services for patients who have not yet been visited, which will cause their dissatisfaction. 

Several models \cite{charnes1959chance,ben2009robust,bernard1955linear} have been proposed in the literature to deal with the uncertainty. The chance constrained model \cite{charnes1959chance} tries to find a solution for which the failure probability is less than some given threshold. In the failure case, the cost of corrective actions is not taken into consideration \cite{gendreau1996stochastic}. In the robust optimization \cite{ben2009robust}, the constructed  solution is feasible for any realization of the uncertainty in a given set \cite{bertsimas2011theory}. This model looks in a minimax fashion for the solution that provides the best “worst case” and solutions may be overly pessimistic. The uncertainty in this model is not stochastic, but rather deterministic and set-based \cite{bertsimas2011theory}. The \ac{SPR} model \cite{bernard1955linear} aims to find a first stage solution that minimizes the expected cost of the second stage solution, which is equal to the cost of the first stage solution plus the expected net cost of recourse \cite{gendreau1996stochastic}. A recourse is introduced in the second stage if a solution goes against the constraints. The objective function of the \ac{SPR} model is more meaningful than the chance constrained model \cite{gendreau1996stochastic}. The \ac{SPR} model is adopted to deal with stochastic parameters.

In this chapter, a \ac{SPR} model is proposed to deal with the  \ac{HHCRSP}  where uncertainties in terms of traveling and caring times are considered. The objective is to minimize the transportation cost and the expected value of recourse. We deal with stochastic travel and service times, the simulation is used since computing the expected real value by an explicit mathematical formula is very complex. Classical solvers such as CPLEX and GUROBI are not suitable to be used with the simulation since their computation time increases monotonically with the problem size and the simulation takes time to converge towards the expected real value. The \ac{GA} based heuristic is adopted and combined with Monte Carlo simulation to solve the \ac{SPR} model. Abundant experiments in the literature show that \ac{GA} based heuristic has a good ability for global searching \cite{shi2018modeling}. In addition, parameters of \ac{GA} are independent of problem parameters (e.g. number of services). In contrast, heuristics that extend and improve local search strategies such as \ac{TS} and \ac{VNS} suffer from this problem, the number of neighborhoods increases with the problem size. At each time a neighborhood is browsed, the simulation is carried out to estimate the expected value of recourse, it requires time to find a good estimation. Indeed, as much as the number of iterations tends to infinity, the estimated value tends towards the expected real value (Law of large numbers). The \ac{GA} is more suitable to solve the \ac{SPR} model as its parameters are not depending on the problem size. The only parameter that could increase the complexity (number of running the simulation) is the population size. However, it is controllable and could be fixed independently of the problem size. In doing so, the complexity faced in heuristics extending local search methods will be avoided. 

We deal with two stochastic problems. In the first one, we propose a  \ac{SPR} model to deal with the  \ac{HHCRSP}  where uncertainties in terms of traveling and caring times that may occur as well as synchronization of services are considered. The objective is to minimize the transportation cost and the expected value of recourse, which is estimated using Monte Carlo simulation. The recourse is defined as a penalty cost for patients’ delayed services and a remuneration for caregivers’ extra working time. We show through numerical tests that \ac{GVNS}  is not suitable to be combined with the simulation to solve the \ac{SPR} model and the adequacy of the \ac{GA} since its parameters do not depend the problem parameters as explained above.
 
In the second one, we assume that time windows are hard/fixed, which must be respected without any earliness or tardiness. To ensure that the requested services are provided within patients' time windows, the recourse is defined as skipping patients when their time windows will be violated. To increase the chance of providing the maximum of   services, patients can specify multiple time windows. The objective is to minimize the transportation cost and the expected value of recourse, which is estimated using Monte Carlo simulation.  The recourse expresses the average number of un-visited patients.

 This chapter is organized as follows: The definition of the problems studied are described in sections \ref{PSMSS} and \ref{PSMTW}. In sections \ref{PSMSS} and \ref{PSMTW}, the problems are formulated as \ac{SPR} models  with a description of the different parameters, variables and constraints taken into account. Sections \ref{GASMSS} and \ref{GASMTW} describe the GA. The test instances, the experimental settings and the performance of the \ac{GVNS}  and/or \ac{GA} are presented in sections \ref{NESMSS} and \ref{NESMTW}. Finally, the chapter ends with a conclusion in section \ref{SMConclusion}.


  \section{Stochastic home health care routing and scheduling problem with multiple synchronized services} \label{SectionSMS}

  \subsection{Problem statement}\label{PSMSS}
  
  in this section we consider the problem statement defined in \ref{PMSMTW} with stochastic travel and service times. However, each patient has only a single time windows in which he is available to be visited.   The problem is to define a daily planning in order to minimize the transportation cost and the expected value of recourse caused by patients’ delayed services and caregivers’ extra working time.  

  \par
  \par
  \par

  
  \subsection{Mathematical formulation} \label{MFSMSS}	
  
  The problem is formulated as a two-stage stochastic programming model with recourse. The first stage aims to compute the transportation cost while the second stage is to introduce the recourse, which is defined as  a penalty cost for patients' delayed services and a remuneration for caregivers' extra working time. The new notation of  sets, decision variables and parameters used in the model are defined as follows:
  
  	
  
  \subsubsection{Deterministic parameters}
  \begin{itemize}
  
  	\item $c_{ij}$ : transportation cost between patients $i$ and $j$.
  	
  \end{itemize}
  
  \subsubsection{Stochastic parameters}
  
  \begin{itemize}
  	
  	\item $ \widetilde{T}_{ij}$: travel time from the patient $i$ to the patient $j$;
  	
  	\item $ \widetilde{t}_{is}$: processing time of the service operation $s$ at the patient $i \in N$;
  	
  	\item $ \mathbb{E}(.)$: the expected value of the second stage of the model, which expresses the recourse value caused by extra working times and delayed  services operations.  
  \end{itemize}
  
  \subsubsection{Decision variables}
  
  \begin{itemize}
  	\item	$ \widetilde{S}_{ik}$: start time of a service operation at the patient $i$ provided by the caregiver $k$.
  \end{itemize}
  
  \subsubsection{Parameters for recourse model}
  
  \begin{itemize}
  	\item $v_i$: tardiness of a service operation at the patient $i$;
  	\item $o_k$: extra working time for the caregiver $k$;
  	
  	\item $\alpha$: unit penalty cost for a tardiness of a service operation;
  	\item $\gamma$: caregiver's remuneration  unit for an extra working time.
  \end{itemize}

  \subsubsection{Mathematical model}
  
  The stochastic programming model with recourse formulation  proposed to solve the \ac{HHCRSP} with stochastic travel and service times, which is adapted from  our previous work \cite{bazirha2019daily} and the  recourse model from  \cite{shi2018modeling}  by adding  other  constraints of the \ac{HHC} context (multiple services and synchronization), is defined as follows:

  \begin{flalign} 
  \min   Z =\sum _{k=1}^c\sum _{i=0}^n \sum _{j=1}^{n+1} c_{ij}x_{ijk} +\mathbb{ E}[\min (\alpha\sum _{i=1}^n  v_{i} +\gamma \sum _{k=1}^c  o_{k})] && \nonumber
  \end{flalign}
  
  \quad s.t.
  
  \setlength{\belowdisplayskip}{0pt} \setlength{\belowdisplayshortskip}{0pt}
  \setlength{\abovedisplayskip}{0pt} \setlength{\abovedisplayshortskip}{0pt}
  
  \begin{flalign} \label{eqn:ptVisits1SMTW}
  \sum _{i=0}^n \sum _{k=1}^c x_{ijk}= \sum_{s=1}^q\delta_{js} 	  && \forall  j\in N
  \end{flalign}
  
  \begin{flalign} \label{eqn:ptVisits2SMTW}
  \sum _{j=1}^{n+1} \sum _{k=1}^c x_{ijk}= \sum_{s=1}^q\delta_{is}   && \forall i\in N
  \end{flalign}

  \begin{flalign}\label{eqn:leftSMTW}
  \sum _{i=0}^n  x_{i, n+1, k}= 1        && \forall  k\in K 
  \end{flalign}

  \begin{flalign}\label{eqn:retrunSMTW}
  \sum _{j=1} ^{n+1}  x_{0, j, k}= 1        && \forall  k\in K
  \end{flalign}

  \begin{flalign}\label{eqn:fluxConservationSMTW}
  \sum _{i=0} ^{n}  x_{imk}= \sum _{j=1}^{n+1}x_{mjk} && \forall  m\in N,  k\in K
  \end{flalign}
  
  
  \begin{flalign}\label{eqn:startingTimeSMTW}
  \widetilde{S}_{ik} +\sum _{s=1}^{q}\widetilde{t}_{is}  y_{iks}  +\widetilde{T}_{ij}  \leq \widetilde{S}_{jk}+(1-x_{ijk})M && \forall 
  \quad  i\in N^0, j\in N^{n+1},  k\in K
  \end{flalign}

  \begin{flalign}\label{eqn:StartTimeToZeroSMTW}
  \widetilde{S}_{ik} \leq \sum_{s=1}^q y_{iks}M  && \forall i\in N, k\in K
  \end{flalign}

  \begin{flalign} \label{eqn:yDefinitionSMTW}
  \sum _{j=1} ^{n+1}  x_{ijk}= \sum_{s=1}^q y_{iks}  && \forall  i\in N,  k\in K 
  \end{flalign}

  \begin{flalign} \label{eqn:skillsRequirementSMTW}
  2 y_{iks} \leq \delta_{is}  +\Delta_{ks}   && \forall i\in N, s \in S, k\in K
  \end{flalign}
  
  \begin{flalign} \label{eqn:assign1SMTW}
  \sum_{s=1}^q y_{iks} \leq 1    && \forall i\in N,  k\in K
  \end{flalign}
  
  \begin{flalign} \label{eqn:assign2SMTW}
  \sum_{k=1}^c y_{iks} \leq 1   && \forall  i\in N,  s\in S
  \end{flalign}

  \begin{flalign}  \label{eqn:twCaregiver1SMTW}
  d_{k} \leq \widetilde{S}_{0k}       && \forall  k\in K 
  \end{flalign}
  
  \begin{flalign} \label{eqn:twCaregiver2SMTW}
  \widetilde{S}_{(n+1)k} \leq  e_{k}+o_k      && \forall  k\in K 
  \end{flalign}

  \begin{flalign} \label{eqn:twPatient1SMTW}
  ( \sum_{s=1}^q y_{iks} -1 ) M + a_{i} \leq 	\widetilde{S}_{ik}     && \forall i\in N, k\in K
  \end{flalign}
  
  \begin{flalign} \label{eqn:twPatient2SMTW}
  \widetilde{S}_{ik} + \sum_{s=1}^q \widetilde{t}_{is} y_{iks}   \leq b_{i} +v_i +(1-\sum_{s=1}^q y_{iks})M 
  && \forall i\in N, k\in K 
  \end{flalign}

  \begin{flalign} \label{eqn:synchronization1SMTW}
  \sum_{v=1}^c \widetilde{S}_{iv}  -  \sum_{s=1}^q\delta_{is} \widetilde{S}_{ik} \leq (2 - \lambda_i - \sum_{s=1}^q y_{iks})M && \forall  i\in N, k\in K
  \end{flalign}
  
  \begin{flalign} \label{eqn:synchronization2SMTW}
  \sum_{v=1}^c \widetilde{S}_{iv}  -  \sum_{s=1}^q\delta_{is} \widetilde{S}_{ik} \geq ( \sum_{s=1}^q y_{iks} + \lambda_i -2 )M && \forall i\in N, k \in K
  \end{flalign}

  \begin{flalign} \label{eqn:xDomSMTW}
  x_{iik} =0      && \forall i\in N,  k\in K
  \end{flalign}

  \begin{flalign} \label{eqn:sDomSMTW}
  \widetilde{S}_{ik} \geq 0   && \forall i\in N,  k\in K
  \end{flalign}

  \begin{flalign} \label{eqn:vDomSMTW}
  v_i \geq 0 	&& \forall  i\in N 
  \end{flalign}	
  
  \begin{flalign} \label{eqn:oDom}
  o_k \geq 0  		&& \forall  k\in K
  \end{flalign}
  
  \begin{flalign} \label{eqn:xDom1SMTW}
  x_{ijk} \in \{0,1\}  && \forall i\in N, j\in N, k\in K
  \end{flalign}

  \begin{flalign}\label{eqn:yDomSMTW}
  y_{iks} \in \{0,1\}  && \forall  i\in N, k\in K, s\in S
  \end{flalign}

  The objective function is to minimize the total transportation cost and the expected value of recourse caused by  patients' tardiness of services operations and caregivers' extra working times. Constraints (\ref{eqn:ptVisits1SMTW})-(\ref{eqn:yDomSMTW}) are already defined in \ref{MFMSMTW}.

  \subsubsection{The expected recourse estimation procedure}

  \begin{algorithm}[!htb]
  	
  	\SetAlgoLined
  	
  	\textbf{Initialization:} \;

  	\hspace{5mm}	 \textbf{- Set} $iter=0$ \;
  	\hspace{5mm}	 \textbf{- $V_k$:} set of patients assigned to the caregiver k   \;
  	
  	\While{(condition 1 or condition 2 )}{
  		\textbf{Generate} randomly $\widetilde{T}_{ij} $ \; 
  		\textbf{Generate} randomly $\widetilde{t}_{js}$   \;
  		
  		\While{($ SS_j$ is changing )}{
  			\textbf{Set} $sum_{T}=0$ and  $sum_{O}=0$\;	
  			
  			\For{$k\leftarrow 1$ \KwTo $K$}  {
  				
  				\For  {$  j\in V_{k} $}{
  					\textbf{Set}  $v_{j}=0$ \;
  					
  					\eIf{$ \lambda_j = 0$}{
  						\textbf{Calculate} $\widetilde{S}_{jk}$   \;
  						\If{$ \widetilde{S}_{jk} + \widetilde{t}_{js} > b_{j}$}{
  							
  							\textbf{Set} $ v_j \longleftarrow  \widetilde{S}_{jk} + \widetilde{t}_{js} - b_{j}$  \;
  						}
  					}{
  						\textbf{Calculate} $\widetilde{SS}_{j}$   \;
  						\If{$ \widetilde{SS}_{j} + \widetilde{t}_{js} > b_{j}$}{
  							
  							\textbf{Set} $ v_j \longleftarrow  \widetilde{SS}_{j} + \widetilde{t}_{js} - b_{j}$  \;
  						}
  					}

  					$sum_T \longleftarrow sum_T + v_j$ \;
  				}
  				
  				\textbf{Set}   $o_{k}=0$ \;
  				\If{$ \widetilde{S}_{(N+1)k}  > e_{k}$}{
  					\textbf{Set} $ o_k \longleftarrow  \widetilde{S}_{(N+1)k}  - e_{k}$ \;
  				}
  				
  				$sum_O \longleftarrow sum_O + o_k$ \;

  			}
  		}
  		$iter \longleftarrow iter + 1 $ \;
  		
  	}
  	\textbf{Set }  $ \widehat{\mathbb{E}(.)} \longleftarrow \frac {\alpha \times sum_T + \gamma \times sum_O}{iter}$ \;
  	\caption{The expected recourse estimation procedure}
  	\label{MCsimulation}
  \end{algorithm}

  The recourse model in stochastic programming is formulated by different ways depending on the nature of the problem. The recourse will be introduced if a solution goes against the constraints. On the \ac{VRP}  with capacity (CVRP) and stochastic demand, the recourse could be defined as returning to the depot if a vehicle is filled and assigned customers are not all visited. This recourse could be applied only to CVRP. In \cite{errico2016priori}, recourse policies are defined as:  skipping the service at the current customer and skipping the visit at the next customer. The main drawback of this recourse is to decrease patients' satisfaction since their requested services could be ignored. In our problem, the recourse model considered in \cite{shi2018modeling} is adopted, since we consider that all patients' requested services will be provided. Constraints related to travel and service times will be relaxed to soft ones, and  recourse will be introduced if a solution goes against these constraints. The recourse is defined as: a penalty cost for a tardiness of a service operation and a remuneration for an extra working time.
  
  \par 
  
  Monte Carlo simulation  \cite{von1951monte} is a method for estimating a numerical quantity that uses random numbers. Computing the real expected value ($ \mathbb{E}(.)$) is very complex since caregivers' routes are dependent and travel and service times are considered stochastic. This simulation is used as an alternative way to compute the estimated value of recourse ($ \widehat{\mathbb{E}(.)}$). $ \widehat{\mathbb{E}(.)}$ gives the average recourse that will be incurred for a given schedule. A robust planning will be performed since different scenarios that might occur have been simulated. The recourse of tardiness of services operations and caregivers' overtimes working  that might arise will be approximated by $ \widehat{\mathbb{E}(.)}$.
  
  \par 
  Algorithm \ref{MCsimulation} is performed to estimate the expected value of recourse ($ \widehat{\mathbb{E}(.)}$). $sum_T$ is the total tardiness of services operations occurred at patients. $sum_O$ is the total extra working times of caregivers and $SS_i$ is the synchronized starting time for the patient $i$.
  Given a solution, patients are assigned to caregivers and caregivers' routes (order of visiting) are defined. Firstly, caregivers' transportation cost are computed, which are  independent of travel and service times. Secondly, Monte Carlo simulation \ref{MCsimulation} is used to estimate the expected value of recourse caused by patients’ delayed services and caregivers’ extra working time. At each iteration, caregivers' travel times $\widetilde{T}_{ij} $ and patients' service times $\widetilde{t}_{is}$ are randomly generated, then the synchronization of services operations will be ensured. The estimation procedure is  running until either condition 1 or condition 2 is met. Condition 1 is fixed as a maximum number of iterations, denoted by $MaxIterMCS$, and condition 2 is defined as follows: 
  $ gap = \abs{\frac{E(.)_{(Iter -1)} -E(.)_{Iter}}{E(.)_{(Iter -1)}}  } < \epsilon $. This condition means that the gap must exactly hold at a maximum number of iterations, denoted by $ MaxIterGap $, and it expresses the gap between $E(.)_{(Iter -1)}$ and $E(.)_{Iter}$ with an error $\epsilon$. In other words, the estimated value is considered close to the real value when the formula is satisfied  ($gap < \epsilon$) for the fixed maximum number of iterations, which means that  $E(.)_{Iter}$ are close to each other for the last $MaxIterGap$ iterations. The counter is set to zero if the formula is not satisfied for an iteration. The more $\epsilon$ is reduced to 0 and  $MaxIterGap$ increases, the more  the estimated value tends towards the  expected real value. However, the computational time increases significantly given that for each individual the expected value is calculated. Condition 1, i.e. $MaxIterMCS$, is used to avoid falling into an infinite loop that could be caused by the condition 2, especially when $\epsilon$ tends towards 0 and $ MaxIterGap $ tends towards a large number. The condition used in \cite{shi2018modeling}, is a particular case of the proposed conditions ( $\epsilon =0$, $MaxIterGap = + \infty$ and $maxIterMCS=100$).

   
  \subsection{Genetic algorithm} \label{GASMSS}
  The main drawback of exact methods is the allowed running time, which is not always enough even to find a feasible solution, especially for NP-hard problems. Several heuristics  have been proposed and are able to yield near-optimal solutions to hard problems in a reasonable amount of time such as GA, which was invented  by John Holland in the 1960s \cite{holland1992genetic}. Abundant  experiments in the literature show that \ac{GA} based heuristics have a good ability for global searching \cite{shi2018modeling}, which involve a simulation of Darwinian "survival of the fittest". This simulation takes an initial population and through the mechanism of reproduction (selection, crossover and mutation), the produced offspring inherits parents’ characteristics and will be added to the next generation. This process keeps on iterating and the fittest individuals will be kept in the end.

  \subsubsection{Crossover operator} \label{crossoverSMSS}
  The crossover operation is the main genetic operator in \ac{GA} used to combine the genetic information of two parents to  stochastically generate new offspring. A large number of crossover operators have been proposed such as 1-point crossover, 2-point crossover and uniform order crossover (UOX). These 3 crossover operators  were implemented and tested on some instances to choose the best one, the UOX operator, was developed by David \cite{davis1991handbook}, gives better results. The UOX preserves the position of some genes and the relative ordering of the rest. After the selection of two parents is done, the two chromosomes are separated and the UOX operator is independently applied for each one (patients’ chromosome and caregivers’ chromosome) and then are recombined. A function is used to repair  infeasible offspring  according to assignment constraints. The UOX operator is applied as follows:
  
  \begin{enumerate}
  	
  	\item For each parent chromosome a binary string of the same length  is  randomly generated;
  	\item The intermediate offspring preserves genes of the first parent where the generated string contains \textbf{“1”};
  	\item Sort genes not preserved in the first parent in the same order as they appear in the second parent. For patients' chromosome, two genes are considered similar if patient's number and the requested service of each gene are equal. Replace genes not preserved in the first parent  by genes of second parent for caregivers' chromosome. 
  	
  \end{enumerate}
  
  Steps 1 and 2 are similar for both chromosomes (caregivers and patients). However, step 3 is adapted for each chromosome. Since the repetition of caregivers does not pose a problem, genes not preserved  in first parent are replaced by genes of the second parent (see Figures \ref{fig:StochPtOff} and \ref{fig:StochCrcross}). In contrast, each patient with the requested service operation must  appear only once in the solution, sorting  genes not preserved in the first parent in the same order as they appear in the second parent is mandatory to avoid  the duplication or the deletion of a patient from the solution (see Figures \ref{fig:StochPtOff} and \ref{fig:StochPtCross})

    \subsubsection{Mutation operator} \label{mutationSMSS}
  
  The role of the mutation operator, which is randomly performed with a small probability, is to avoid being trapped in local optima, to excavate the diversity of the individuals in the population and to diversify the search directions. Two mutation operators are proposed, patients' mutation operator is to exchange two patients' positions randomly generated, and caregivers' mutation operator is to switch an assigned caregiver to a patient randomly. Taking the example of the Table \ref{table:StochMutation1} and supposing that the probabilities of carried out mutation operators are satisfied for both chromosomes. The position for caregivers' chromosome is 1 and positions for patients' chromosome are 3 and 5. Table \ref{table:StochMutation2} shows the obtained solution after carry out the mutation operators.

  \begin{figure}[!htb]
  	
  	\begin{mdframed}[style=MyFrame,nobreak=true,align=center,userdefinedwidth=31em]
  		
  		\raggedleft
  		
  		Parent 1:
  		\begin{tabular}{ | c | c | c | c | c | c | c|c| }
  			\hline
  			\color{red} 5 (2)  & 2 (1) &6 (2)& \color{red} 3 (1) &\color{red} 1 (2) &  4 (3) & \color{red}1 (3) \\ 
  			\hline
  			
  			2    &\color{blue} 1 &  1 & 2    & \color{blue}1    & \color{blue}2     &   2\\ 
  			
  			\hline
  		\end{tabular}
  		
  		\vspace{1mm}
  		
  		Parent 2:
  		\begin{tabular}{ | c | c | c | c | c | c | c|c| }
  			\hline
  			
  			1 (3)  & 3 (1) & 1 (2) & 4 (3) & 2 (1) &  6 (2) & 5 (2) \\ 
  			\hline
  			1  & 1  & 2 & 2 & 1  & 1  &  2\\   
  			
  			\hline
  			
  		\end{tabular}

  		\vspace{1mm}
  		
  		Offspring:
  		\begin{tabular}{ | c | c | c | c | c | c | c|c| }
  			\hline
  			\color{red} 5 (2)  & 4 (3) & 2 (1)& \color{red} 3 (1) &\color{red} 1 (2) &  6 (2)& \color{red}1 (3) \\ 
  			\hline
  			1   & \color{blue}1 & 2 & 2    &\color{blue} 1    & \color{blue}2     &  2\\   
  			\hline
  		\end{tabular}

  	\end{mdframed}
  	\caption{Example of parents and an offspring }
  	\label{fig:StochPtOff}
  \end{figure}


  \begin{figure}[!htb]
  	
  	\begin{mdframed}[style=MyFrame,nobreak=true,align=center,userdefinedwidth=31em]
  		\raggedleft
  		
  		Binary string :
  		\begin{tabular}{  |c | c | c | c | c | c | c|c| }

  			\hline
  			 \hspace{1.95mm} 0  \hspace{1.95mm}   &  \hspace{1.95mm} 1  \hspace{1.95mm} &  \hspace{1.95mm}  0 \hspace{1.95mm}  &  \hspace{1.95mm} 0  \hspace{1.95mm}   &  \hspace{1.95mm} 1  \hspace{1.95mm} &  \hspace{1.95mm} 1  \hspace{1.95mm}  &  \hspace{1.95mm}   0  \hspace{1.95mm} \\   
  			\hline
  			
  		\end{tabular}
  		
  		\vspace{1mm}
  		
  		Parent 1:
  		\begin{tabular}{| c | c | c | c | c | c | c|c| }
  			
  			\hline
  			 \hspace{1.95mm} 2  \hspace{1.95mm}    & \hspace{1.95mm} \color{blue} 1 \hspace{1.95mm}  &  \hspace{1.95mm}  1  \hspace{1.95mm} &  \hspace{1.95mm} 2  \hspace{1.95mm}    &  \hspace{1.95mm} \color{blue}1   \hspace{1.95mm}   &  \hspace{1.95mm} \color{blue}2  \hspace{1.95mm}     &  \hspace{1.95mm}   2 \hspace{1.95mm} \\   
  			\hline
  		\end{tabular}
  		
  		\vspace{1mm}
  		
  		Parent 2:
  		\begin{tabular}{ | c | c | c | c | c | c | c|c| }
  			
  			\hline
  			 \hspace{1.95mm} 1  \hspace{1.95mm}   &  \hspace{1.95mm} 1  \hspace{1.95mm} &  \hspace{1.95mm} 2  \hspace{1.95mm} &  \hspace{1.95mm} 2  \hspace{1.95mm}    &  \hspace{1.95mm} 1  \hspace{1.95mm}    &  \hspace{1.95mm} 1  \hspace{1.95mm}    &  \hspace{1.95mm}   2 \hspace{1.95mm} \\  
  			\hline
  			
  		\end{tabular}
  		
  		\vspace{1mm}
  		
  		Step 2:
  		\begin{tabular}{ | c | c | c | c | c | c | c|c| }
  			
  			\hline
  			\hspace{8.7mm} & \color{blue}  \hspace{1.95mm} 1  \hspace{1.95mm}  &  \hspace{8.7mm} &   \hspace{8.7mm}    &\color{blue}  \hspace{1.95mm} 1   \hspace{1.95mm}   & \color{blue}  \hspace{1.95mm} 2  \hspace{1.95mm}   &  \hspace{8.7mm} \\   
  			\hline
  		\end{tabular}
  		
  		\vspace{1mm}
  		
  		Step 3:
  		\begin{tabular}{ | c | c | c | c | c | c | c|c| }
  			
  			\hline
  			 \hspace{1.95mm} 1  \hspace{1.95mm}   & \hspace{8.7mm}   &  \hspace{1.95mm} 2 \hspace{1.95mm}  &  \hspace{1.95mm} 2  \hspace{1.95mm}  &  \hspace{8.7mm}  & \hspace{8.7mm}      &  \hspace{1.95mm}  2  \hspace{1.95mm} \\   
  			\hline
  		\end{tabular}
  		
  		\vspace{1mm}
  		
  		offspring:
  		\begin{tabular}{ | c | c | c | c | c | c | c|c| }
  			
  			\hline
  		 \hspace{1.95mm} 	1   \hspace{1.95mm}  &  \hspace{1.95mm} \color{blue}1  \hspace{1.95mm}  & \hspace{1.95mm}  2  \hspace{1.95mm} &  \hspace{1.95mm} 2  \hspace{1.95mm}  & \hspace{1.95mm} \color{blue} 1 \hspace{1.95mm}     & \hspace{1.95mm}  \color{blue} 2  \hspace{1.95mm}     &  \hspace{1.95mm}  2 \hspace{1.95mm} \\   
  			\hline
  		\end{tabular}

  	\end{mdframed}
  	\caption{Example of  caregivers'  crossover  operator}
  	\label{fig:StochCrcross}
  \end{figure}


  \begin{figure}[!htb]
  	
  	\begin{mdframed}[style=MyFrame,nobreak=true,align=center,userdefinedwidth=31em]

  		\raggedleft
  		
  		Binary string :
  		\begin{tabular}{  |c | c | c | c | c | c | c|c| }
  			
  			\hline
  			  \hspace{1.95mm} 1   \hspace{1.95mm}    &   \hspace{1.95mm} 0  \hspace{1.95mm}  &   \hspace{1.95mm}  0   \hspace{1.95mm} &   \hspace{1.95mm} 1   \hspace{1.95mm}    &   \hspace{1.95mm} 1   \hspace{1.95mm}   &   \hspace{1.95mm} 0   \hspace{1.95mm}    &  \hspace{1.95mm} 1 \hspace{1.95mm} \\   
  			\hline			
  		\end{tabular}
  		
  		\vspace{1mm}
  		Parent 1:
  		\begin{tabular}{ | c | c | c | c | c | c | c|c| }
  			\hline
  			\color{red} 5 (2)  & 2 (1) &6 (2)& \color{red} 3 (1) &\color{red} 1 (2) &  4 (3) & \color{red}1 (3) \\ 
  			\hline
  			
  		\end{tabular}
  		
  		\vspace{1mm}
  		Parent 2:
  		\begin{tabular}{ | c | c | c | c | c | c | c|c| }
  			\hline
  			1 (3)  & 3 (1) & 1 (2) & 4 (3) & 2 (1) & 6 (2) & 5 (2) \\ 
  			\hline
  		\end{tabular}
  		
  		\vspace{1mm}
  		
  		Step 2:
  		\begin{tabular}{ | c | c | c | c | c | c | c|c| }
  			\hline
  			\color{red} 5 (2)  & \hspace{8.7mm} & \hspace{8.7mm} & \color{red} 3 (1) &\color{red} 1 (2) & \hspace{8.7mm} & \color{red}1 (3) \\ 
  			\hline
  			
  		\end{tabular}
  		
  		\vspace{1mm}
  		
  		Step 3:
  		\begin{tabular}{ | c | c | c | c | c | c | c|c| }
  			\hline
  			\hspace{8.7mm}  & 4 (3) & 2 (1)& \hspace{8.7mm} & \hspace{8.7mm}  &  6 (2)& \hspace{8.7mm} \\ 
  			\hline
  		\end{tabular}
  		
  		\vspace{1mm}
  		
  		Offspring:
  		\begin{tabular}{ | c | c | c | c | c | c | c|c| }
  			\hline
  			\color{red} 5 (2)  & 4 (3) & 2 (1)& \color{red} 3 (1) &\color{red} 1 (2) &  6 (2)& \color{red}1 (3) \\ 
  			\hline
  		\end{tabular}

  	\end{mdframed}
  	\caption{Example of  patients' crossover  operator}
  	\label{fig:StochPtCross}
  \end{figure}

  \begin{table}[ !htb]
  	
  	\begin{center}
  		
  		\caption{Example of a solution}
  		
  		\begin{tabular}{ |l| c | c | c | c | c | c | c|c| }
  			\hline
  			Patients  & 5 (2)  & 2 (1) &\color{red} \textbf{6  (2)}& 3  (1) & \color{red}\textbf{1 (2)} &  4 (3) & 1 (3) \\ 
  			\hline
  			
  			Caregivers & \color{red} \textbf{2}    & 1 &  1 & 2    &  1    &  2     &   2\\ 
  			
  			\hline
  		\end{tabular}
  		\label{table:StochMutation1}
  		
  	\end{center}
  \end{table}
  \begin{table}[!htb]
  	
  	\centering
  	
  	\caption{Example of mutation operators}
  	
  	\begin{tabular}{ |l| c | c | c | c | c | c | c|c| }
  		\hline
  		Patients  & 5  (2)  &   2  (1) & \textbf{\color{blue}1 	 (2)}& 3  (1) &\textbf{\color{blue}6  (2)} &4  (3)  & 1  (3) \\ 
  		\hline
  		
  		Caregivers & \textbf{\color{blue}1 }    & 1 &  1 & 2   &  1    &  2     &   2\\ 
  		
  		\hline
  		
  	\end{tabular}
  	\label{table:StochMutation2}
  	
  \end{table}

  \subsubsection{Fitness and selection }
  In selection phase, the fittest individuals are selected for reproduction, which let them pass their genes to the next generation. The most commonly known selection methods include roulette rank selection, wheel selection, and tournament selection \cite {goldberg1991comparative}. The first one suffers from the slow convergence and the sorting is done for the population to assign ranks, which increase the computational time. The second suffers from problem of premature convergence due to the possible presence of dominant individuals that always win the competition and are selected as a parent, which will be rewarded with a large number of offspring in the next generation while the population size is kept constant \cite{baker1985adaptive}.
  \par
  In this study, the tournament selection will be used. It consists of selecting some individuals from the population, then those individuals compete against each other. The one with the highest fitness wins and participates in reproduction. Two fitness functions are used to evaluate solutions, $F_{D}$ is used for the deterministic  model and $F_{SPR}$ is used for the SPR model and are defined as follows:
  
  \begin{flalign} \label{fitD}
  F_{D} =\sum _{k=1}^c\sum _{i=0}^n \sum _{j=1}^{n+1} c_{ij}x_{ijk}   + \beta( T_s + O_c)
  \end{flalign}
  
  \begin{flalign} \label{fitSPR}
  F_{SPR} =\sum _{k=1}^c\sum _{i=0}^n \sum _{j=1}^{n+1} c_{ij}x_{ijk} + \mathbb{E}(.)  
  \end{flalign}

  Infeasible solutions according to time windows constraints are accepted with a penalty cost, which is the sum of tardiness of services operations $T_s$ and caregivers' overtimes $O_c$. The coefficient $\beta$ is introduced to ensure the convergence of infeasible solutions to feasible ones for the deterministic model, the more $\beta$ is higher the more penalized solutions are to be eliminated from the population in selection step. The individual with small fitness will be selected. 
  
  \subsubsection{Genetic algorithm procedure}
  
  At each time crossover and mutation operations are applied to the selected parents, the generated offspring might not be feasible. A service operation could be assigned to a caregiver who is not qualified to provide. To avoid this issue, the repair function is called after applying crossover and mutation operations. This function checks each assigned service to a caregiver and verify if he is qualified to provide it, otherwise, a randomly qualified caregiver is selected.
  
  \begin{algorithm}[!htb]
  	
  	\SetAlgoLined
  
  	\textbf{Initialization:} \;
  	\hspace{5mm}	 \textbf{- Define} $P_{size}$: population size\;
  	\hspace{5mm}	 \textbf{- Define}  $P_c$: crossover probability   \;
  	\hspace{5mm}	 \textbf{- Define}  $P_S$: mutation probability   \;
  	\hspace{5mm}	 \textbf{- Define} $P$: population \;
  	\hspace{5mm}	 \textbf{- Define } $O_{ff}$: Offspring \;
  	
  	\While{(the stopping condition is not reached)}{
  		\For{$i\leftarrow 1$ \KwTo $P_{size}$}  {

  			\textbf{Selection}: $P_1$ and $P_2$   \;
  			\textbf{Generate randomly} $p_1$   \;
  			
  			\If{$ p_c < p_1$}{
  				\textbf{Set} $ Off \longleftarrow  Crossover(P_1,P_2)$  \;
  				
  			}
  			\textbf{Generate randomly} $p_2$   \;
  			
  			\If{$ p_m < p_2$}{
  				\textbf{Set} $ Off \longleftarrow  Mutation(Off)$  \;
  				
  			}
  			
  			\textbf{Repair} $Off$   \;
  			
  		}

  		\textbf{Set} $ P \longleftarrow  Off$  \;
  		
  	}
  	
  	\caption{GA procedure}
  		\label{GA}
  \end{algorithm}
  
  \subsubsection{Initial population}
  
  The initial population for a given size $P_{size}$ is randomly generated as follows:
  
  \begin{itemize}
  	\item For $i =1$ to $P_{size}$ do:
  	\item  Generate a random visiting order (see Table \ref{table:VOSMSS});
  	
  	\begin{table}[ !htb]
  		
  		\caption{Example of visiting order}
  		\centering
  		\centering
  		\begin{tabular}{ | c | c | c | c | c | c | c|c| }
  			\hline
  			Patients  & 5 (2)  & 1 (3)  & 3 (1) &1 (2)& 4 (3) & 2 (1) &  6 (2)  \\ 
  			\hline
  			Caregivers  &    & &  & &     &     &   \\   
  			\hline
  		\end{tabular}
  		
  		\label{table:VOSMSS}
  		
  	\end{table}
  	
  	\item For each patient, assign a qualified caregiver randomly selected (see Table \ref{table:AssignSMSS});
  	\begin{table}[ !htb]
  		
  		\caption{example of caregivers’ assignment to patients}
  		\centering
  		\begin{tabular}{ | c | c | c | c | c | c | c|c| }
  			\hline
  			Patients  & 5 (2)  & 1 (3)  & 3 (1) &1 (2)& 4 (3) & 2 (1) &  6 (2)  \\ 
  			\hline
  			Caregivers  & 1    & 2 &  2 & 1   & 1    & 2     &   1\\   
  			\hline
  		\end{tabular}

  		\label{table:AssignSMSS}
  		
  	\end{table}
  	
  	\item Calculate the transportation cost and the expected value of recourse for the stochastic model by equation \ref{fitSPR}. For the deterministic model, the transportation cost, penalized tardiness of services operations and caregivers' extra working times are computed by equation \ref{fitD}.

  \end{itemize}
  
  \subsection{Variable neighborhood search} \label{GVNSSMSS}
 
   The \ac{GVNS}  based heuristic and its parameters are further explained in the previous chapter (see sections \ref{VNSMTW} and \ref{VNSMSMTW}).

  \subsection{Numerical experiments }\label{NESMSS}
 
   \subsubsection{Tuning parameters}
  A series of tests are performed to find the best values of tuning parameters (see Table \ref{table:StochTP}). The stopping criterion for both \ac{GA} and \ac{GVNS}  based heuristics is set as the number of no improvement in the best solution found for a maximum number of iterations. The total requested services of all patients is denoted by $s$.
  
  \begin{table}[!hbt]
  	\caption{Tuning parameters}
  	
  	\centering
  	\renewcommand{\arraystretch}{1.2}
  	\begin{tabular}{ l|  l | l }
  		
  		\hline
  		\textbf{Algorithm }  & \textbf{Parameter} &  \textbf{Value} \\ 
  		\hline
  		
  		&$P_c$: crossover probability  & 0.8 \\ 
  		GA based heuristic  	&$P_s$: mutation probability   & 0.01 \\
  		(Deterministic model) &$P_{size} $ : population size & $s^2$ \\
  		& Stopping criterion   & $5s$ \\
  		&Tournament size      & $c+1$ \\
  		
  		\hline

  		& Shaking phase & $h=c+1$ \\
  		GVNS	based heuristic   & Stopping criterion   & $2s$ \\
  		& shaking phases order   &switch, inter-swap, intra-swap and shift \\
  		& local search methods order   & shift, switch, inter-swap and intra-swap \\

  		\hline
  		
  		& $\epsilon$ & $0.05$ \\
  		Monte Carlo simulation	& $MaxIterMCS$   & 100 \\
  		& $MaxIterGap$   & 10 \\

  		\hline
  		
  		&$P_c$: crossover probability  & 0.6 \\ 
  		GA based heuristic  	 &$P_s$: mutation probability  & 0.08 \\
  		(SPR model) &$P_{size} $ : population size & 100 \\
  		& Stopping criterion   & $50$ \\
  		&Tournament size      & $c+1$ \\
  		
  		\hline

  		Others	& $MaxIterSyn$    & $2c$\\
  		& $\beta$   & 100 \\

  		\hline
  		
  	\end{tabular}
  	\label{table:StochTP}
  	
  \end{table}
  \subsubsection{Test instances}\label{TISMSS}
  
  The test instances have been randomly generated. Deterministic parameters are generated as described in \ref{TIDMTW}. Travel time $T_{ij}$ and cost $c_{ij}$ are equal to the Euclidean distance $d_{ij}$ between patients’ locations truncated to an integer. In practice, cost and travel time between patients are quasi-proportional to  the distance. It is assumed that $T_{ij} =coef_1 \times d_{ij}$, $c_{ij} =coef_2 \times d_{ij}$ and $coef_1=coef_2=1$. Each patient requires single or double services, which are randomly drawn from $S=\{1,…,6\}$. 30\% of patients are considered requesting double services \cite{mankowska2014home}, 50\% of double services are supposed to be simultaneous  and the others 50\% are supposed to be received without synchronization. Two sets of instances are generated single services ($SS$) and multiple synchronized services ($MSS$), each one contains 3 subsets and are summarized in the Table \ref{table:StochTID}. '$Size$' is the number of instances in each subset. '$\#N$' is the number of patients. '$\#S$' is the total requested services by all patients. '$\#K$' is the number of caregivers available.

  \begin{table}[H]
  	\caption{ Tested instances details}
  	\centering
  	\def\arraystretch{1.3}
  	
  	\begin{tabular}{ | l|  l | l | l | l | l |l |}
  		
  		\hline
  		
  		Set & Subset & Size& $\# N$ &$\#S$  &$\#K$     \\   
  		
  		\hline
  		& A  &	7	& 10	&	10 & 3    \\   
  		
  		SS   & B &	7	& 25	&	25 & 5    \\   
  		
  		& C &	7	& 50	&	50 & 10   \\   
  		\hline
  		& D &	7	& 10	&	13 & 3   \\   
  		
  		MSS 	& E &	7	& 25	&	33	 & 5    \\   
  		
  		& F &	7	& 50	&	65 	& 10    \\

  		\hline

  	\end{tabular}

  	\label{table:StochTID}

  \end{table}
  
   Stochastic parameters (travel and service times) are  randomly generated the same way in \cite{shi2018modeling}. The travel time between each two patients follows a normal  distribution    $ \widetilde{T}_{ij} \sim \mathcal{N}(c_{ij}, (\frac{c_{ij}}{3})^2 ) $ and the service time also follows a normal  distribution $ \widetilde{t}_{is} \sim \mathcal{N}(t_{is},(\frac{t_{is}}{5})^2 ) $. $c_{ij}$ and $t_{is}$ are respectively the average of the traveling time between two patients $i$ and $j$ and the average time of a service operation $s$ at patient $i$. $\frac{c_{ij}}{3}$ and  $\frac{t_{is}}{5}$ are respectively the standard deviation values  of the traveling time between two patients $i$ and $j$ and service operation time $s$ for the patient $i$. $\alpha$ and $\gamma$ are set to 1.

  \subsubsection{Computational results}


  Instances are generated as described above  and solved within a time limit of 4 hours. CPLEX could solve the subsets $A$, $B$, $C$, $D$ and $E$ optimally. For the subset $F$, only instances $F1$ and $F6$ could be solved optimally. For other instances $F2$, $F3$, $F4$, $F5$ and $F7$ a feasible solution is found  with a gap respectively of 10.48\%, 2.66\%, 12.25\%, 10.84\%  and 1.26\%. The \ac{GVNS}  based heuristic could find the optimal solutions for the subsets $A$, $B$, $ D$ and $E$ at least once among the ten runs of each instance, except for the instances $E2$. The optimal solution is found also for the instances $C4$, $C5$ and  $C7$. The \ac{GA} based heuristic could find the optimal solutions for the subsets $A$, $B$ and $ D$ and for the instances $E1$, $E6$ and  $E7$ at least once among the ten runs of each instance (see Tables \ref{table:MIPSS} and \ref{table:MIPMSS}).
  \par
  
  The worst gaps  found by the \ac{GVNS}   for the subsets $A, B, C, D, E$ and $F$ are respectively  $0.99\%$ (instance $A6$), $0.77\%$ (instance $B6$), $3.53\%$ (instance $C4$), $0.00\%$ (subset $D$), $2.21\%$ (instance $E2$) and $17.71\%$ (instance $F5$). For the GA, are respectively  $2.60\%$ (instance $A6$), $2.82\%$ (instance $B2$), $14.98\%$ (instance $C4$), $1.04\%$ (instance $D6$), $5.71\%$ (instance $E3$) and $28.35\%$ (instance $F5$). Both \ac{GVNS}  and \ac{GA} based heuristics could find good quality solutions with a better efficiency for the \ac{GVNS}  (see Fig. \ref{DSS} and \ref{DMSS}).

  \par

  \begin{figure}[H]
  	\vspace{-7mm}
  	\begin{center}
  		\begin{tikzpicture}
  		\begin{axis}[
  		width=\textwidth, 
  		height=8cm,
  		grid=major, 
  		grid style={dashed,gray!30},
  		xlabel= Instances, 
  		ylabel= Objective function values,
  		xtick=data,
  		legend style={at={(0.08,0.6)},anchor=west},
  		,xticklabels={A1,A2,A3,A4,A5,A6,A7,B1,B2,B3,B4,B5,B6,B7,C1,C2,C3,C4,C5,C6,C7}
  		]
  		
  		\addplot+[] 
  		table[x=SS, y=Z1, col sep=comma] {src/benchmarking.txt}; 
  		
  		\addplot+[] 
  		table[x=SS, y=GVNS1, col sep=comma] {src/benchmarking.txt}; 
  		
  		\addplot+[] 
  		table[x=SS, y=GA1, col sep=comma] {src/benchmarking.txt};

  		\legend{CPLEX,GVNS average, \ac{GA} average}
  		\end{axis}
  		
  		\end{tikzpicture}
  		\vspace{-8mm}
  		\caption{Comparison of solutions found by CPLEX, \ac{GVNS}  and \ac{GA}  for the  instances of the set $SS$  with deterministic parameters.}
  		\label{DSS}

  		\begin{tikzpicture}
  		\begin{axis}[
  		width=\textwidth, 
  		height=8cm,
  		grid=major, 
  		grid style={dashed,gray!30},
  		xlabel= Instances, 
  		ylabel= Objective function values,
  		xtick=data,
  		legend style={at={(0.08,0.6)},anchor=west},
  		,xticklabels={D1,D2,D3,D4,D5,D6,D7,E1,E2,E3,E4,E5,E6,E7,F1,F2,F3,F4,F5,F6,F7}
  		]
  		
  		\addplot+[] 
  		table[x=MSS, y=Z2, col sep=comma] {src/benchmarking.txt}; 
  		
  		\addplot+[] 
  		table[x=MSS, y=GVNS2, col sep=comma] {src/benchmarking.txt}; 
  		
  		\addplot+[] 
  		table[x=MSS, y=GA2, col sep=comma] {src/benchmarking.txt};

  		\legend{CPLEX,GVNS average, \ac{GA} average}
  		\end{axis}
  	
  		\end{tikzpicture}
  		\vspace{-8mm}
  		\caption{Comparison of solutions found by CPLEX, \ac{GVNS}  and \ac{GA}  for the  instances of the set $MSS$  with deterministic parameters.}
  		\label{DMSS}
  	\end{center}

  \begin{center}
  	
  	\begin{tikzpicture}
  	\begin{axis}
  	[
  	width=14cm, 
  	height=7cm,
  	xlabel= Objective function values,
  	ytick={1,2,3},
  	yticklabels={CPLEX, \ac{GVNS}  Average, \ac{GA} Average},
  	grid=major, 
  	grid style={dashed,gray!30}, 
  	]
  	\addplot+[ 
  	boxplot prepared={
  		average=1185.79,
  		median=1157.50,
  		upper quartile=1567.00	,
  		lower quartile=777.00,
  		upper whisker=2011.00,
  		lower whisker=439.00
  	},
  	] coordinates{(0, 521)	(0, 715)	(0, 508)	(0, 817)	(0, 645)	(0, 439)	(0, 539)	(0, 1165)	(0, 993)	(0, 1131)	(0, 928)	(0, 1064)	(0, 1196)	(0, 1099)	(0, 1491)	(0, 1673)	(0, 1567)	(0, 1561)	(0, 1671)	(0, 1658)	(0, 1492)	(0, 769)	(0, 872)	(0, 709)	(0, 938)	(0, 777)	(0, 588)	(0, 609)	(0, 1317)	(0, 1361)	(0, 1338)	(0, 1150)	(0, 1246)	(0, 1251)	(0, 1145)	(0, 1754)	(0, 1864)	(0, 1731)	(0, 1962)	(0, 2011)	(0, 1808)	(0, 1730)
  	};
  	\addplot+[
  	boxplot prepared={
  		average=1210.31,
  		median=1163.30,
  		upper quartile=1618.10,
  		lower quartile=777.00,
  		upper whisker=2179.00,
  		lower whisker=443.40	
  	},
  	] coordinates {(0, 521)	(0, 715)	(0, 508)	(0, 817)	(0, 645)	(0, 443.4)	(0, 539)	(0, 1167.2)	(0, 996.6)	(0, 1131)	(0, 928)	(0, 1064)	(0, 1205.3)	(0, 1099)	(0, 1530.7)	(0, 1728.6)	(0, 1601.9)	(0, 1618.1)	(0, 1710.6)	(0, 1713.4)	(0, 1538.2)	(0, 769)	(0, 872)	(0, 709)	(0, 938)	(0, 777)	(0, 588)	(0, 609)	(0, 1336.1)	(0, 1391.7)	(0, 1357.5)	(0, 1159.4)	(0, 1268.6)	(0, 1259.3)	(0, 1147.6)	(0, 1903)	(0, 1933.6)	(0, 1822.7)	(0, 1955.3)	(0, 2179)	(0, 1853.1)	(0, 1783.1)
  	};
  	\addplot+[
  	boxplot prepared={
  		average= 1288.62,
  		median=1188.40,
  		upper quartile=1836.00,
  		lower quartile=777.00,
  		upper whisker=2502.50,
  		lower whisker=450.70
  	},
  	] 
  	coordinates {(0, 521.1)	(0, 718.9)	(0, 509.8)	(0, 817)	(0, 649.9)	(0, 450.7)	(0, 539)	(0, 1169.6)	(0, 1021.8)	(0, 1136.4)	(0, 947.4)	(0, 1064.4)	(0, 1215.9)	(0, 1099.6)	(0, 1635)	(0, 1901.2)	(0, 1735.4)	(0, 1836)	(0, 1901.3)	(0, 1837.4)	(0, 1646.8)	(0, 769)	(0, 872)	(0, 709.8)	(0, 938)	(0, 777)	(0, 594.2)	(0, 609)	(0, 1382.7)	(0, 1407.8)	(0, 1419.1)	(0, 1207.2)	(0, 1286.9)	(0, 1286.7)	(0, 1155.1)	(0, 2169.4)	(0, 2166.2)	(0, 2033.5)	(0, 2193.5)	(0, 2502.5)	(0, 2189)	(0, 2098.8)	};

  	\end{axis}
  	\end{tikzpicture}
  	
  	\vspace{-2mm}
  	\caption{Comparison of solutions found by CPLEX, \ac{GVNS}  and \ac{GA} with deterministic parameters.}
  	\label{Boxplot}
  \end{center}
  \end{figure}

 \begin{sidewaystable}

  		\caption{Results of solving the set $SS$ instances with deterministic parameters.}
  		\centering
  		\renewcommand{\arraystretch}{1.3}
  		\begin{tabular}{c|cccc | ccccc |ccccc}

  			\hline
  			\multicolumn{1}{c|}{} & \multicolumn{4}{c|}{\textbf{CPLEX} }&\multicolumn{5}{c|}{ \textbf{\ac{GVNS} (10 runs)}} &\multicolumn{5}{c}{ \textbf{ \ac{GA} (10 runs)}}\\ 
  			
  			\hline
  			
  			\textbf{ SS}    & \textbf{LB} & \textbf{Z} & \textbf{ Gap} &\textbf{CPU} & \textbf{Best} & \textbf{Worst}& \textbf{Average} &\textbf{ Gap} & \textbf{CPU}
  			& \textbf{Best} & \textbf{Worst}& \textbf{Average} &\textbf{ Gap} & \textbf{CPU}
  			\\ 
  			
  			\hline	
  			
  			A1	&	521	&	\textbf{521}	&	0.00\%	&	1.73	&	\textbf{521}	&	\textbf{521}	&	\textbf{521.00}	&	0.00\%	&	$<1$	&	\textbf{521}	&	522	&	521.10	&	0.02\%	&	$<1$	\\
  			A2	&	715	&	\textbf{715}	&	0.00\%	&	1.49	&	\textbf{715}	&	\textbf{715}	&	\textbf{715.00}	&	0.00\%	&	$<1$	&	\textbf{715}	&	754	&	718.90	&	0.54\%	&	$<1$	\\
  			A3	&	508	&	\textbf{508}	&	0.00\%	&	1.43	&	\textbf{508}	&	\textbf{508}	&	\textbf{508.00}	&	0.00\%	&	$<1$	&	\textbf{508}	&	526	&	509.80	&	0.35\%	&	$<1$	\\
  			A4	&	817	&	\textbf{817}	&	0.00\%	&	1.52	&	\textbf{817}	&	\textbf{817}	&	\textbf{817.00}	&	0.00\%	&	$<1$	&	\textbf{817}	&	\textbf{817}	&	\textbf{817.00}	&	0.00\%	&	$<1$	\\
  			A5	&	645	&	\textbf{645}	&	0.00\%	&	1.58	&	\textbf{645}	&	\textbf{645}	&	\textbf{645.00}	&	0.00\%	&	$<1$	&\textbf{	645}	&	660	&	649.90	&	0.75\%	&	$<1$	\\
  			A6	&	439	&	\textbf{439}	&	0.00\%	&	1.44	&	\textbf{439}	&	461	&	443.40	&	0.99\%	&	$<1$	&	\textbf{439}	&	461	&	450.70	&	2.60\%	&	$<1$	\\
  			A7	&	539	&	\textbf{539}	&	0.00\%	&	153.00	&	\textbf{539}	&	\textbf{539}	&	\textbf{539.00}	&	0.00\%	&	$<1$	&	\textbf{539}	&	\textbf{539}	&	\textbf{539.00}	&	0.00\%	&	$<1$	\\
  			
  			\hline	
  			
  			B1	&	1165	&	\textbf{1165}	&	0.00\%	&	2.66	&	\textbf{1165}	&	1178	&	1167.20	&	0.19\%	&	5.89	&	\textbf{1165}	&	1179	&	1169.60	&	0.39\%	&	14.96	\\
  			B2	&	993	&	\textbf{993}	&	0.00\%	&	3.58	&	\textbf{993}	&	1003	&	996.60	&	0.36\%	&	7.30	&	\textbf{993}	&	1070	&	1021.80	&	2.82\%	&	14.39	\\
  			B3	&	1131	&	\textbf{1131}	&	0.00\%	&	2.36	&	\textbf{1131}	&	\textbf{1131}	&	\textbf{1131.00}	&	0.00\%	&	10.22	&	\textbf{1131}	&	1142	&	1136.40	&	0.48\%	&	16.05	\\
  			B4	&	928	&	\textbf{928}	&	0.00\%	&	2.07	&	\textbf{928}	&	\textbf{928}	&	\textbf{928.00}	&	0.00\%	&	4.99	&	\textbf{928}	&	964	&	947.40	&	2.05\%	&	15.34	\\
  			B5	&	1064	&	\textbf{1064}	&	0.00\%	&	2.04	&	\textbf{1064}	&	\textbf{1064}	&	\textbf{1064.00}	&	0.00\%	&	5.68	&	\textbf{1064}	&	1070	&	1064.40	&	0.04\%	&	15.69	\\
  			B6	&	1196	&	\textbf{1196}	&	0.00\%	&	2.99	&	\textbf{1196}	&	1221	&	1205.30	&	0.77\%	&	6.36	&	\textbf{1196}	&	1243	&	1215.90	&	1.64\%	&	17.97	\\
  			B7	&	1099	&	\textbf{1099}	&	0.00\%	&	2.29	&	\textbf{1099}	&	\textbf{1099}	&	\textbf{1099.00}	&	0.00\%	&	4.39	&	\textbf{1099}	&	1102	&	1099.60	&	0.05\%	&	17.72	\\
  			
  			\hline	
  			
  			C1	&	1491	&	\textbf{1491}	&	0.00\%	&	40.70	&	1502	&	1559	&	1530.70	&	2.59\%	&	196.59	&	1570	&	1745	&	1635.00	&	8.81\%	&	251.65	\\
  			C2	&	1673	&	\textbf{1673}	&	0.00\%	&	206.39	&	1679	&	1777	&	1728.60	&	3.22\%	&	197.97	&	1751	&	1967	&	1901.20	&	12.00\%	&	254.08	\\
  			C3	&	1567	&	\textbf{1567}	&	0.00\%	&	673.93	&	1570	&	1627	&	1601.90	&	2.18\%	&	167.96	&	1650	&	1796	&	1735.40	&	9.70\%	&	275.72	\\
  			C4	&	1561	&	\textbf{1561}	&	0.00\%	&	27.19	&	\textbf{1561}	&	1702	&	1618.10	&	3.53\%	&	193.99	&	1621	&	1856	&	1836.00	&	14.98\%	&	291.48	\\
  			C5	&	1671	&	\textbf{1671}	&	0.00\%	&	99.44	&	\textbf{1671}	&	1784	&	1710.60	&	2.31\%	&	234.03	&	1817	&	2047	&	1901.30	&	12.11\%	&	338.54	\\
  			C6	&	1658	&	\textbf{1658}	&	0.00\%	&	23.85	&	1683	&	1739	&	1713.40	&	3.23\%	&	179.61	&	1738	&	1904	&	1837.40	&	9.76\%	&	322.35	\\
  			C7	&	1492	&	\textbf{1492}	&	0.00\%	&	32.76	&	\textbf{1492}	&	1595	&	1538.20	&	3.00\%	&	196.39	&	1601	&	1763	&	1646.80	&	9.40\%	&	320.45	\\
  			
  			\hline	
  			
  		\end{tabular}

  		\label{table:MIPSS}

  	\end{sidewaystable}

  	  	\begin{sidewaystable}

  		\caption{Results of solving the set $MSS$ instances with deterministic parameters.}
  		\centering
  		\renewcommand{\arraystretch}{1.3}
  		\begin{tabular}{c|cccc | ccccc |ccccc}

  			\hline
  			\multicolumn{1}{c|}{} & \multicolumn{4}{c|}{\textbf{CPLEX} }&\multicolumn{5}{c|}{ \textbf{\ac{GVNS} (10 runs)}} &\multicolumn{5}{c}{ \textbf{\ac{GA} (10 runs)}}\\ 
  			\hline

  			\textbf{ MSS}    & \textbf{LB} & \textbf{Z} & \textbf{ Gap} &\textbf{CPU} & \textbf{Best} & \textbf{Worst}& \textbf{Average} &\textbf{ Gap} & \textbf{CPU}
  			& \textbf{Best} & \textbf{Worst}& \textbf{Average} &\textbf{ Gap} & \textbf{CPU}	\\
  			\hline
  			
  			D1	&	769	&	\textbf{769}	&	0.00\%	&	1.64	&	\textbf{769}	&	\textbf{769}	&	\textbf{769.00}	&	0.00\%	&	$<1$	&	\textbf{769}	&	\textbf{769}	&	\textbf{769.00}	&	0.00\%	&	1.44	\\
  			D2	&	872	&	\textbf{872}	&	0.00\%	&	1.48	&	\textbf{872}	&	\textbf{872}	&	\textbf{872.00}	&	0.00\%	&	$<1$	&	\textbf{872}	&	\textbf{872}	&	\textbf{872.00}	&	0.00\%	&	1.22	\\
  			D3	&	709	&	\textbf{709}	&	0.00\%	&	1.78	&	\textbf{709}	&	\textbf{709}	&	\textbf{709.00}	&	0.00\%	&	$<1$	&	\textbf{709}	&	717	&	709.80	&	0.11\%	&	1.45	\\
  			D4	&	938	&	\textbf{938}	&	0.00\%	&	1.56	&	\textbf{938}	&	\textbf{938}	&	\textbf{938.00}	&	0.00\%	&	$<1$	&	\textbf{938}	&	\textbf{938}	&	\textbf{938.00}	&	0.00\%	&	1.25	\\
  			D5	&	777	&	\textbf{777}	&	0.00\%	&	1.55	&	\textbf{777}	&	\textbf{777}	&	\textbf{777.00}	&	0.00\%	&	$<1$	&	\textbf{777}	&	\textbf{777}	&	\textbf{777.00}	&	0.00\%	&	1.45	\\
  			D6	&	588	&	\textbf{588}	&	0.00\%	&	1.47	&	\textbf{588}	&	\textbf{588	}&	\textbf{588.00}	&	0.00\%	&	$<1$	&	\textbf{588}	&	606	&	594.20	&	1.04\%	&	1.39	\\
  			D7	&	609	&	\textbf{609}	&	0.00\%	&	1.48	&	\textbf{609}	&	\textbf{609}	&	\textbf{609.00}	&	0.00\%	&	$<1$	&	\textbf{609}	&	\textbf{609}	&	\textbf{609.00}	&	0.00\%	&	1.41	\\
  			
  			\hline

  			E1	&	1317	&	\textbf{1317}	&	0.00\%	&	18.52	&	\textbf{1317}	&	1362	&	1336.10	&	1.43\%	&	59.29	&	\textbf{1317}	&	1464	&	1382.70	&	4.75\%	&	60.88	\\
  			E2	&	1361	&	\textbf{1361}	&	0.00\%	&	8.07	&	1374	&	1401	&	1391.70	&	2.21\%	&	57.50	&	1387	&	1447	&	1407.80	&	3.32\%	&	58.41	\\
  			E3	&	1338	&	\textbf{1338}	&	0.00\%	&	7.15	&	\textbf{1338}	&	1383	&	1357.50	&	1.44\%	&	56.04	&	1353	&	1503	&	1419.10	&	5.71\%	&	63.79	\\
  			E4	&	1150	&	\textbf{1150}	&	0.00\%	&	5.57	&	\textbf{1150}	&	1219	&	1159.40	&	0.81\%	&	46.39	&	1171	&	1260	&	1207.20	&	4.74\%	&	55.50	\\
  			E5	&	1246	&	\textbf{1246}	&	0.00\%	&	3.59	&	\textbf{1246}	&	1294	&	1268.60	&	1.78\%	&	47.62	&	1247	&	1356	&	1286.90	&	3.18\%	&	64.12	\\
  			E6	&	1251	&	\textbf{1251}	&	0.00\%	&	11.25	&	\textbf{1251}	&	1289	&	1259.30	&	0.66\%	&	51.57	&	\textbf{1251}	&	1309	&	1286.70	&	2.77\%	&	67.60	\\
  			E7	&	1145	&	\textbf{1145}	&	0.00\%	&	5.14	&	\textbf{1145}	&	1157	&	1147.60	&	0.23\%	&	32.59	&	\textbf{1145}	&	1210	&	1155.10	&	0.87\%	&	46.88	\\
  			
  			\hline	
  			
  			F1	&	1754	&	\textbf{1754}	&	0.00\%	&	1971	&	1808	&	2002	&	1903.00	&	7.83\%	&	2766	&	2069	&	2251	&	2169.40	&	19.15\%	&	1138	\\
  			F2	&	1668.73	&	1864	&	10.48\%	&	14400	&	1832	&	2120	&	1933.60	&	13.70\%	&	1863	&	1989	&	2312	&	2166.20	&	22.97\%	&	1197	\\
  			F3	&	1684.98	&	1731	&	2.66\%	&	14400	&	1813	&	1839	&	1822.70	&	7.56\%	&	2052	&	1940	&	2141	&	2033.50	&	17.14\%	&	1291	\\
  			F4	&	1721.72	&	1962	&	12.25\%	&	14400	&	1908	&	2016	&	1955.30	&	11.95\%	&	2312	&	2105	&	2269	&	2193.50	&	21.51\%	&	969	\\
  			F5	&	1793.01	&	2011	&	10.84\%	&	14400	&	2040	&	2367	&	2179.00	&	17.71\%	&	2735	&	2296	&	2837	&	2502.50	&	28.35\%	&	1396	\\
  			F6	&	1808	&	\textbf{1808}	&	0.00\%	&	2207	&	1818	&	1907	&	1853.10	&	2.43\%	&	2096	&	2035	&	2309	&	2189.00	&	17.41\%	&	1487	\\
  			F7	&	1708.12	&	1730	&	1.26\%	&	14400	&	1738	&	1832	&	1783.10	&	4.21\%	&	1960	&	2012	&	2213	&	2098.80	&	18.61\%	&	1136	\\

  			\hline
  		\end{tabular}

  		\label{table:MIPMSS}

  \end{sidewaystable}

  \par
  
  In the SPR model, the transportation cost, the expected value of recourse and the objective function ($Z=cost+\widehat{\mathbb{E}(.)}$) are considered. Small instances (subset $A$ and $D$) are solved by both  \ac{GA}  and \ac{GVNS}  based heuristics. Solutions found by both heuristics are very close (see Table \ref{table:SPRSmallInstances}, the gap is computed as   (average$_{GVNS}$- average$_{GA}$)/average$_{GVNS}$).

  \begin{table*}[!htb]
  	
  	\caption{Results of solving  small instances (subsets $A$ and $D$) with stochastic parameters using \ac{GA} and \ac{GVNS}  based heuristics.}
  	\centering
  	\renewcommand{\arraystretch}{1.3}
  	
  	\begin{tabular*}{\textwidth}{@{\extracolsep{\fill}} c| cccc |cccc|c}
  		
  		\hline
  		\multicolumn{1}{c|}{} &\multicolumn{4}{c|}{ \textbf{\ac{GVNS} (10 runs)}} &\multicolumn{4}{c|}{ \textbf{\ac{GA} (10 runs)}}&\\ 
  		\hline

  		\textbf{}    & \textbf{Best} & \textbf{Worst}& \textbf{Average} & \textbf{CPU}
  		& \textbf{Best} & \textbf{Worst}& \textbf{Average} & \textbf{CPU} & \textbf{Gap}
  		\\ 
  		\hline
  		
  		A1	&	521.38	&	522.79	&	521.93	&	266.53	&	520.70	&	522.25	&	521.52	&	95.83	&	0.08\%	\\
  		A2	&	722.41	&	725.51	&	724.29	&	365.41	&	720.96	&	723.90	&	722.36	&	119.08	&	0.27\%	\\
  		A3	&	495.49	&	497.42	&	496.45	&	324.63	&	494.32	&	498.66	&	495.43	&	87.28	&	0.20\%	\\
  		A4	&	815.68	&	819.73	&	817.67	&	421.12	&	812.79	&	817.90	&	814.49	&	134.84	&	0.39\%	\\
  		A5	&	616.88	&	634.61	&	619.85	&	494.64	&	613.76	&	631.85	&	624.46	&	113.63	&	-0.74\%	\\
  		A6	&	430.95	&	434.63	&	432.94	&	477.53	&	428.35	&	452.80	&	432.83	&	127.54	&	0.03\%	\\
  		A7	&	537.15	&	538.80	&	537.96	&	408.06	&	534.82	&	536.95	&	535.92	&	122.84	&	0.38\%	\\
  		\hline
  		D1	&	778.68	&	783.76	&	780.38	&	1139.46	&	776.41	&	818.87	&	786.06	&	255.71	&	-0.73\%	\\
  		D2	&	872.23	&	872.98	&	872.71	&	621.22	&	872.06	&	872.99	&	872.52	&	184.24	&	0.02\%	\\
  		D3	&	688.40	&	698.00	&	693.05	&	918.01	&	683.96	&	697.95	&	694.09	&	192.73	&	-0.15\%	\\
  		D4	&	883.85	&	890.67	&	888.17	&	944.82	&	881.55	&	884.81	&	883.92	&	253.85	&	0.48\%	\\
  		D5	&	753.77	&	759.03	&	755.88	&	978.60	&	749.68	&	768.82	&	752.90	&	223.57	&	0.39\%	\\
  		D6	&	578.98	&	583.96	&	581.76	&	1221.90	&	576.52	&	588.86	&	580.31	&	273.40	&	0.25\%	\\
  		D7	&	609.43	&	610.66	&	609.94	&	1162.30	&	609.27	&	614.73	&	611.04	&	275.68	&	-0.18\%	\\

  		\hline
  	\end{tabular*}

  	\label{table:SPRSmallInstances}
  	
  \end{table*}

    Fig. \ref{Boxplot} shows that the minimum, first quartile and median found by CPLEX, \ac{GVNS}  average and \ac{GA} average are very close to each other, which indicate that the performances of these methods in the first 50\% of instances (subsets $A, B$ and $D$) are similar. According to Tables \ref{table:MIPSS} and \ref{table:MIPMSS}, optimal solutions are reached by both heuristics at least once in 10 runs. For the subset $E$, both heuristics still could reach solutions very close to those found by CPLEX. \ac{GVNS}  showed  a better efficiency and could reach 6 optimal solutions from 7 at least once in 10 runs. Fig.\ref{Boxplot} shows also that solutions reached by \ac{GA} are very dispersed in the second 50\% of instances (subsets $C, E$ and $F$). Fig. \ref{DSS} and \ref{DMSS} show that  \ac{GA} found solutions  for  subsets $C$ and $F$  a little bit further from those reached by CPLEX and \ac{GVNS}  compared to solutions of subset $E$. Therefore, we can conclude that this dispersion  is caused by objective function values of instances in subsets C and F. The solutions reached by \ac{GVNS}  and CPLEX have a very close distribution if we delete the outlier (the maximum).

  \par
  
Table \ref{table:SPRSmallInstances} shows clearly that \ac{GVNS}  is not suitable to be combined with the simulation to solve the SPR model. The worst CPU running time jumped from 494.64 seconds (instance $A5$ of the set $SS$) to 1221.90 seconds (instance $D6$ of the set $MSS$) while the worst CPU running time of \ac{GA} jumped from 134.184 seconds (instance $A4$ of the set $SS$) to 275.68 seconds (instance $D7$ of the set $MSS$). This is due to the complexity of neighborhoods sizes, which depends on the size $s$ of the problem (number of total requested services by all patients). Three from the four proposed neighborhood structures of \ac{GVNS}  have a complexity of $O(s^2)$. Exploring neighborhood solutions imply running the simulation for each new generated solution to estimate the expected value of recourse, which significantly increases the CPU running time. To illustrate that, we solved the first instance of each subset ($A1, B1, C1, D1, E1$ and $F1$) by \ac{GVNS}  for a single iteration. Fig. \ref{CPUTimeSPRGVNS} clearly shows that the CPU running time increases monotonically with the total number of services requested by all patients. The CPU running time jumped from 1.41 seconds (instance $A1$ with 10 services) to 2120.45 seconds (instance $F1$ with 65 services) to run \ac{GVNS}  for a single iteration.

  \begin{figure}[H]
  
  	\begin{center}
  		
  		\begin{tikzpicture}
  		\begin{axis}[
  		xlabel={Instances (\textit{Total requested services s})},
  		ylabel={ CPU running time (\textit{seconds})},
  		width=14cm,
  		height=9cm,
  		xticklabels={A1 \textit{(s=10)} , D1 \textit{(s=13)}, B1 \textit{(s=25)}, E1 \textit{(s=33)}, C1 \textit{(s=50)}, F1 \textit{(s=65)}},
  		xtick={1,2,3,4,5,6},
  		legend pos=north west,
  		grid=major, 
  		grid style={dashed,gray!30},
  		]
  		
  		\addplot[nodes near coords,
  		color=black,
  		mark=*,
  		]
  		table[meta=label] {
  			x   y       label
  			1	1.41	1.41 
  			2	3.54	3.54
  			3	60.04	60.04
  			4	102.39	102.39
  			5	1803.71	1803.71
  			6	2120.45	2120.45
  			
  		};
  		;

  		\end{axis}
  		\end{tikzpicture}
  		\caption{CPU running times elapsed to solve the first instance of each subset using the simulation embedded into \ac{GVNS}  for a single iteration.}
  		\label{CPUTimeSPRGVNS}
  	\end{center}
  \end{figure}
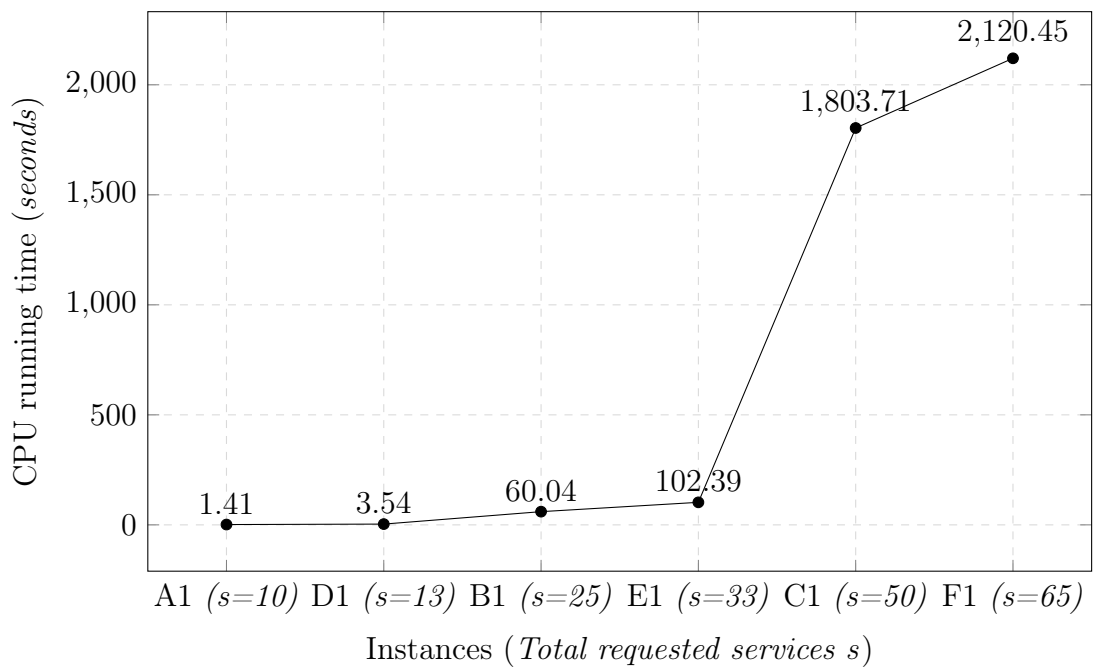

  	In the other hand, \ac{GA} parameters are independent of the problem size.  To solve medium (subsets $B$ and $E$) and large (subsets $C$ and $F$) instances  in reasonable computational time, we fixed the population size and the stopping criterion independently of the problem size (see Table \ref{table:StochTP}: SPR model). The SPR model is very complex in terms of CPU running time compared to the deterministic model, due to the expected value that must be estimated for each individual and in each iteration of the algorithm (see Table \ref{table:SPRLargeInstances1}). In other words, the expected recourse estimation procedure is executed each time the selection-crossover-mutation mechanism is applied.

  \begin{table*}[!htb]

  	\caption{Results of solving  medium and large instances  with stochastic parameters using \ac{GA} based heuristic.}
  	\centering
  	\renewcommand{\arraystretch}{1.3}
  	\begin{tabular*}{\textwidth}{@{\extracolsep{\fill}} c| cccc | c| cccc }
  		\hline
  		\multicolumn{10}{c}{ \textbf{\ac{GA} (10 runs)}}\\ 
  		\hline
  		
  		\textbf{SS}   &  \textbf{Best} & \textbf{Worst}& \textbf{Average} & \textbf{CPU}&	\textbf{MSS}   &  \textbf{Best} & \textbf{Worst}& \textbf{Average} & \textbf{CPU}	\\ 
  		\hline
  		
  		B1	&	1134.77	&	1342.94	&	1210.03	&	549.03	&	E1	&	1318.66	&	2030.43	&	1677.94	&	957.04	\\
  		B2	&	1002.53	&	1197.89	&	1085.04	&	748.76	&	E2	&	1371.74	&	2180.87	&	1673.77	&	862.95	\\
  		B3	&	1135.76	&	1250.39	&	1176.94	&	930.63	&	E3	&	1479.07	&	1767.76	&	1604.45	&	881.13	\\
  		B4	&	942.71	&	1081.69	&	1016.55	&	811.72	&	E4	&	1267.52	&	1481.07	&	1361.58	&	925.49	\\
  		B5	&	1076.11	&	1308.81	&	1167.70	&	841.30	&	E5	&	1325.51	&	2223.86	&	1636.49	&	836.41	\\
  		B6	&	1158.15	&	1289.61	&	1221.40	&	614.21	&	E6	&	1388.26	&	2595.32	&	1750.48	&	818.38	\\
  		B7	&	1116.84	&	1224.54	&	1182.24	&	1035.93	&	E7	&	1235.28	&	2688.13	&	1640.10	&	853.78	\\
  		\hline
  		C1	&	1692.52	&	2068.92	&	1856.73	&	4714.29	&	F1	&	2942.55	&	4547.42	&	3363.43	&	4756.74	\\
  		C2	&	1930.80	&	2244.35	&	2117.14	&	5442.66	&	F2	&	2465.93	&	3610.59	&	3114.22	&	5739.92	\\
  		C3	&	1775.89	&	1978.25	&	1881.40	&	6233.01	&	F3	&	2253.21	&	3138.75	&	2635.77	&	5478.67	\\
  		C4	&	1882.94	&	2368.98	&	2129.19	&	4268.39	&	F4	&	2388.92	&	3876.00	&	3148.09	&	4959.69	\\
  		C5	&	1950.67	&	2221.42	&	2063.72	&	6063.64	&	F5	&	2861.24	&	4806.11	&	3520.34	&	4753.81	\\
  		C6	&	1974.69	&	2135.67	&	2041.33	&	4233.30	&	F6	&	2522.03	&	3356.50	&	2908.03	&	5183.04	\\
  		C7	&	1712.97	&	2037.94	&	1826.58	&	6020.90	&	F7	&	2018.86	&	3162.31	&	2801.73	&	5535.88	\\

  		\hline
  	\end{tabular*}

  	\label{table:SPRLargeInstances1}
  	
  \end{table*}
  
  \par

  	In Table \ref{table:SPRvsD} we compared when deterministic optimal solutions of small instances (subsets $A$ and $D$) are considered as solutions for the stochastic model with solutions found by the SPR model.  The first stage (caregivers’ routes and their assignment to patients) is the optimal solution given by CPLEX. In the second stage, we run the simulation to estimate the expected value of recourse $ \widehat{\mathbb{E}(.)}$. The total cost $Z$ includes the transportation cost found by CPLEX and the cost of the recourse computed by the simulation, the costs of the two stages are computed independently. The \ac{GVNS}  (respectively GA) column  contains the average solution of 10 runs for solving the SPR model and the gaps are computed as \textit{(Z – Average)/Z}. In most cases, the values of the gaps are greater than 1, which shows that deterministic optimal solutions are not privileged in the stochastic case since the recourse is not optimized simultaneously with the transportation cost. In addition, the SPR model prioritizes solutions with higher difference between latest service times ($b_{i}$ and $e_{k}$) and caregivers' completion $C_{ik}$ to minimize the expected value of recourse. The more the margins $b_{i}$- $C_{ik}$ and $e_{k}$- $C_{ik}$ are higher the more the fluctuation of caregivers' arrival times remains robust and compatible with time windows.
  
  \begin{table}[!htb]
 
  	\caption{Deterministic solutions with stochastic parameters and solutions found by the SPR model for the small instances.}
  	\centering
  	\renewcommand{\arraystretch}{1.3}
  	
  	\begin{tabular*}{\textwidth}{@{\extracolsep{\fill}} c| c|cc |cc|cc}
  		
  		\hline
  		\multicolumn{1}{c|}{} & \textbf{CPLEX}& \multicolumn{2}{c|}{ \textbf{Simulation}}&	\multicolumn{2}{c|}{ \textbf{\ac{GVNS} (10 runs)}} &	\multicolumn{2}{c}{ \textbf{\ac{GA} (10 runs)}}\\ 
  		\hline

  		\textbf{}    & \textbf{Cost} & \textbf{$ \widehat{\mathbb{E}(.)}$}& \textbf{Z} & \textbf{Average} & \textbf{Gap}& \textbf{Average}  & \textbf{Gap}
  		\\ 
  		\hline
  		
  		A1	&	521	&	5.67	&	526.67	&	521.93	&	0.90\%	&	521.52	&	0.98\%	\\
  		A2	&	715	&	20.02	&	735.02	&	724.29	&	1.46\%	&	722.36	&	1.72\%	\\
  		A3	&	508	&	4.69	&	512.69	&	496.45	&	3.17\%	&	495.43	&	3.37\%	\\
  		A4	&	817	&	19.8	&	836.8	&	817.67	&	2.29\%	&	814.49	&	2.67\%	\\
  		A5	&	645	&	10.52	&	655.52	&	619.85	&	5.44\%	&	624.46	&	4.74\%	\\
  		A6	&	439	&	1.53	&	440.53	&	432.94	&	1.72\%	&	432.83	&	1.75\%	\\
  		A7	&	539	&	5.13	&	544.13	&	537.96	&	1.13\%	&	535.92	&	1.51\%	\\
  		
  		\hline
  		
  		D1	&	769	&	19.65	&	788.65	&	780.38	&	1.05\%	&	786.06	&	0.33\%	\\
  		D2	&	872	&	2.27	&	874.27	&	872.71	&	0.18\%	&	872.52	&	0.20\%	\\
  		D3	&	709	&	40.1	&	749.1	&	693.05	&	7.48\%	&	694.09	&	7.34\%	\\
  		D4	&	938	&	19.63	&	957.63	&	888.17	&	7.25\%	&	883.92	&	7.70\%	\\
  		D5	&	777	&	28.38	&	805.38	&	755.88	&	6.15\%	&	752.90	&	6.52\%	\\
  		D6	&	588	&	1.68	&	589.68	&	581.76	&	1.34\%	&	580.31	&	1.59\%	\\
  		D7	&	609	&	5.86	&	614.86	&	609.94	&	0.80\%	&	611.04	&	0.62\%	\\

  		\hline
  	\end{tabular*}

  	\label{table:SPRvsD}
  	
  \end{table}
  \par
  
    To sum up, even though the exact methods reach the optimal solution, their computational time increases monotonically with the problem size. The two proposed heuristics \ac{GVNS}  and \ac{GA} are able to solve large instances in a short computation time. \ac{GVNS}  showed the better efficiency in solving the instances for the deterministic model compared to GA. Despite this, \ac{GVNS}  was not suitable to be combined with simulation to solve the SPR model because three from the proposed neighborhood structures have a running time complexity of $O(s^2)$. To overcome this limitation, we embedded the simulation into \ac{GA} to solve the SPR model since its parameters do not depend on the problem size. Although the SPR model is very complex with a high run time, this model provides a robust scheduling in which uncertainties in term of traveling and caring times are taken into account. The computed expected value gives the average recourse that will be incurred for a given schedule without needing to know which patients have received the services with a tardiness and which caregivers have worked overtime.

  In this section, we introduce the recourse if a solution goes against constraints related to patients’ time windows and caregivers’ duty length.  We convert time windows from hard to soft using penalization for tardiness of services and remuneration for caregivers’ overtime. In the next section, we will investigate the case of travel and care times uncertainty with hard/fixed time windows. In others word, patients' availability periods must be respected, even with the presence of the uncertainty.


   \section{Stochastic home health care routing and scheduling problem with multiple hard time windows} \label{SectionSMTW}
  
  \subsection{Problem statement}\label{PSMTW}
  
 We make the following modifications to the problem defined in \ref{PSMSS} to introduce the new problem.  To ensure fairness among caregivers, a maximum number of patients not to exceed is fixed ($max_v$) instead of using the duty length for caregivers. Each patient requires a single service operation and allowed to be available in many periods $ [a_{il}, b_{il}]$. Parameters $a_{il}$ and $b_{il}$ are, respectively, the earliest and latest possible service times of the period $l \in L=\{1,2,...,p\}$, where $p$ is the number of patients' time windows. The decision maker could select any availability period $l$ to schedule the requested visit for each patient $i$.
  
  It is more likely that providing some services will not be compatible with patients’ preferred availability periods since travel $ \widetilde{T}_{ij}$ and service  $ \widetilde{t}_{is}$ times are supposed stochastic and patients’ time windows are assumed hard/fixed. Therefore, the recourse is defined as skipping a patient when carrying out the service will not be compatible with his availability periods. The goal is to define a daily planning that minimize the transportation cost and the expected value of recourse defined as minimizing the number of unvisited patients with respect to skills requirements, patients' times windows and the maximum of visits not to exceed by  each caregiver.

  \subsection{Mathematical formulation}\label{MFSMTW}
  The problem is formulated as a two-stage stochastic programming model recourse with multiple hard/fixed time windows. The first stage aims to define caregivers’ routes (patients’ visiting order) and patients’ assignments to caregivers. The second stage is to introduce the recourse caused by caregivers’arriving lately to patients, which is expressed as skipping patients' visits since their time windows are assumed hard/fixed. The recourse is to minimize the average number of unvisited patients. In the following, we present the notation of  sets, decision variables and parameters not defined  in the previous section \ref{MFSMSS} and those redefined in this section.
  
  
  \subsubsection{Deterministic parameters}
  \begin{itemize}

  	\item 	$Max_v$: maximum of patients that a caregiver could visit.
  	
  \end{itemize}
  
  \subsubsection{Stochastic parameters}
  
  \begin{itemize}
  	
  	
  	
  	\item $ \mathbb{E}(.)$: the expected value of the recourse of the second stage, which expresses the average number of uninvited patients.  
  \end{itemize}
  
  
  
  \subsubsection{Parameters for recourse model}
  
  \begin{itemize}
  	\item $v_i$:  equals  1 if the service operation requested by the patient $i$ will be skipped (will not be provided), 0 otherwise;
  	
  	\item $\alpha$: penalty cost for each unvisited patient.
  	
  \end{itemize}

  \subsubsection{Mathematical model}
  
  The stochastic programming recourse model formulation  proposed to solve this problem is adapted from \cite{bazirha2020scheduling} by using another recourse model  \cite{errico2016priori} (skip patients) instead of using a penalty cost for violating time windows \cite{shi2018modeling}. In addition, we add other constraints of the \ac{HHC} context such as multiple time windows and maximum number of patients to visit per caregiver. The model is defined as follows:

  \begin{flalign} 
  	\min  Z =\sum _{k=1}^c\sum _{i=0}^n \sum _{j=1}^{n+1} c_{ij}x_{ijk} +\mathbb{ E}[\min \sum _{i=1}^n  \alpha v_{i} ] && \nonumber
  \end{flalign}
  
  \quad s.t.
  
  \setlength{\belowdisplayskip}{0pt} \setlength{\belowdisplayshortskip}{0pt}
  \setlength{\abovedisplayskip}{0pt} \setlength{\abovedisplayshortskip}{0pt}
  
  \begin{flalign} \label{eqn:ptVisits1SMSS}
  	\sum _{i=0}^n \sum _{k=1}^c x_{ijk}=  1  &&  \forall  j\in N 
  \end{flalign}
  
  \begin{flalign}\label{eqn:ptVisits2SMSS}
  	\sum _{j=1}^{n+1} \sum _{k=1}^c x_{ijk}= 1 &&  \forall  i\in N 
  \end{flalign}

  \begin{flalign}\label{eqn:leftSMSS}
  	\sum _{i=0}^n  x_{i(n+1)k}= 1         &&  \forall  k\in K 
  \end{flalign}
  
  \begin{flalign}\label{eqn:retrunSMSS}
  	\sum _{j=1} ^{n+1}  x_{0jk}= 1         &&  \forall  k\in K
  \end{flalign}
  
  \begin{flalign}\label{eqn:fluxConservationSMSS}
  	\sum _{i=0} ^{n}  x_{imk}= \sum _{j=1}^{n+1}x_{mjk} &&  \forall  m\in N,  k\in K
  \end{flalign}
  
  
  \begin{flalign}\label{eqn:startingTime1SMSS}
  	\widetilde{S}_{ik}+\sum _{s=1}^{q}\widetilde{t}_{is}  y_{iks}  +\widetilde{T}_{ij}  \leq \widetilde{S}_{jk}+(1+v_i -x_{ijk})M &&  \forall
  	i\in N^0, j\in N^{n+1},  k\in K 
  \end{flalign}
  
  \begin{flalign}\label{eqn:startingTime2SMSS}
  	\widetilde{S}_{ik} +\widetilde{T}_{ij}  \leq \widetilde{S}_{jk}+(2-v_i-x_{ijk})M &&  \forall  i\in N^0, j\in N^{n+1},  k\in K
  \end{flalign}

  \begin{flalign}\label{eqn:YdefinitionSMSS}
  	\sum _{j=1} ^{n+1}  x_{ijk}= \sum_{s=1}^q y_{iks} &&  \forall  i\in N,  k\in K 
  \end{flalign}

  \begin{flalign}\label{eqn:skillsRequirementSMSS}
  	2 y_{iks} \leq \delta_{is}  +\Delta_{ks}  &&  \forall i\in N, s \in S, k\in K
  \end{flalign}

  \begin{flalign}\label{eqn:MaxToVisitSMSS}
  	\sum _{i=1} ^{n} \sum_{s=1}^q y_{iks}\leq Max_v &&  \forall  k\in K 
  \end{flalign}
  
  

  \begin{flalign}\label{eqn:twPatient1SMSS}
  	( \sum_{s=1}^q y_{iks} +z_{il} - v_i -2 ) M + a_{i} \leq 	\widetilde{S}_{ik}  &&  \forall   i\in N, l\in L, k\in K
  \end{flalign}
  
  \begin{flalign}\label{eqn:twPatient2SMSS}
  	\widetilde{S}_{ik} + \sum_{s=1}^q \widetilde{t}_{is} y_{iks}   \leq b_{i}  +(2 + v_i - z_{il} -\sum_{s=1}^q y_{iks})M &&  \forall i\in N, l\in L, k\in K 
  \end{flalign}
  
  \begin{flalign}\label{eqn:zDefinitionSMSS}
  	\sum_{l=1}^p z_{il} =1  &&  \forall  i\in N
  \end{flalign}

  \begin{flalign} \label{eqn:xDomSMSS}
  	x_{iik} =0  &&  \forall  i\in N,  k\in K
  \end{flalign}
  
  \begin{flalign} \label{eqn:sDomSMSS}
  	S_{ik} \geq 0  &&  \forall  i\in N,  k\in K 
  \end{flalign}

  \begin{flalign} \label{eqn:vDomSMSS}
  	v_i \geq 0 &&  \forall	  i\in N 
  \end{flalign}

  \begin{flalign} \label{eqn:xDom1SMSS}
  	x_{ijk} \in \{0,1\}  &&  \forall  i\in N, j\in N, k\in K
  \end{flalign}

  \begin{flalign} \label{eqn:yDomSMSS}
  	y_{iks} \in \{0,1\} \ &&  \forall   i\in N, k\in K, s\in S
  \end{flalign}
  
  \begin{flalign} \label{eqn:zDomSMSS}
  	z_{il} \in \{0,1\}  &&  \forall   i\in N, l\in L
  \end{flalign}

  The objective function is defined as minimizing caregivers' transportation cost and the expected value of the recourse caused by skipping patients if their availability periods are not respected. 
  
   Constraints (\ref{eqn:startingTime1SMSS}) and (\ref{eqn:startingTime2SMSS}) determine either service operation $s$ requested by  patient $i$ will be provided, or it will be skipped. Indeed, if ($v_i=0$, i.e. Constraints (\ref{eqn:startingTime1SMSS}) must be verified),  the service operation will be provided for  patient $i$ and the start time  for patient $j$ must respect completion time of providing the requested service operation for patient $i$. Otherwise ($ v_i =1$, i.e. Constraints (\ref{eqn:startingTime2SMSS}) must be verified), the service operation for patient $i$ will be skipped. Constraints (\ref{eqn:MaxToVisitSMSS}) guarantee that each caregiver does not exceed the maximum number of visits allowed. The rest of constraints are already explained in chapter \ref{ChapterDM}.

  \subsubsection{The expected recourse estimation procedure}
  
  The recourse model in stochastic programming depends on the nature of the problem and its constraints. Constraints containing stochastic parameters are more likely to be violated, therefore a recourse must be used to deal with the uncertainty.  In \cite{shi2018modeling}, the authors defined the recourse as a penalty cost for a tardiness of a service operation and a remuneration for caregivers’ extra working time. This recourse requires to be used with soft/flexible time windows.  In \cite{errico2016priori}, the authors defined the recourse when a route becomes infeasible as:  skipping the service at the current customer and skipping the visit at the next customer. This recourse is used with hard/fixed time windows. Since we suppose that patients’ time windows must be respected, we define the recourse as skipping providing a service operation for a patient when his availability periods could not be respected. The algorithm \ref{MCsimulation2} is used to estimate the excepted value of recourse. $sum_v$ is the total number of skipped visits. $T_j$ is the minimal tardiness of providing the service operation to patient $j$ considering all his availability periods. $p_j$ contains the selected period $l$ for the patient $j$. $V_k$ contains patients assigned to the caregiver $k$.

  \begin{algorithm}[!htb]
  	
  	\SetAlgoLined
  
  	\textbf{Initialization:} \;
  	\hspace{5mm}	 \textbf{- Set} $sum_{v}=0$ \;
  	\hspace{5mm}	 \textbf{- Set} $T_{j}=0$ \;
  	
  	\hspace{5mm}	 \textbf{- Set} $iter=0$ \;
  	
  	\While{( condition 1 or condition 2 is not reached )}{
  		\For{$k\leftarrow 1$ \KwTo $K$}  {
  			
  			\For  {$  j\in V_{k} $}{
  				
  				\textbf{Generate} randomly $\widetilde{T}_{ij} $ \; 
  				\textbf{Calculate} $\widetilde{A}_{jk}$ \;
  				\textbf{Generate} randomly $\widetilde{t}_{js}$   \;
  				
  				$T_{j} \leftarrow +\infty$ \;
  				\For{$l\leftarrow 1$ \KwTo $L$}  {
  					\textbf{Compute} the tardiness $T_l$ of providing the service $s$ using the period $l$
  					
  					\If{$ T_l <T_j $}{
  						
  						\textbf{Set} $ T_j \longleftarrow  T_l$  \;
  						\textbf{Set} $ p_j \longleftarrow  l$  \;
  					}
  				}

  				\textbf{Set} $ v_j \longleftarrow  0$  \;
  				\eIf{$ T_j  > 0 $}{
  					
  					\textbf{Set} $ v_j \longleftarrow  1$ \;
  				}  { 
  					\textbf{Compute} $\widetilde{S}_{jk}$  using the period $p_j$  \;
  				}
  				
  				$sum_v \longleftarrow sum_v + v_j$ \;
  			}

  		}
  		
  		$iter \longleftarrow iter + 1 $ \;
  		
  	}
  	\textbf{Set }  $\mathbb{ E}(.)\longleftarrow \frac {\alpha \times sum_v }{iter}$ \;
  	\caption{The expected recourse estimation procedure}
  		\label{MCsimulation2}
  \end{algorithm}

  \subsection{Genetic algorithm} \label{GASMTW}
   
   The same mutation operators defined  in \ref{mutationSMSS} are used. However, in the following, we define the new crossover operators used to solve the SPR model with multiple soft/flexible time windows.  As we mentioned  in \ref{crossoverSMSS}, The 3 crossover operators: 1-point crossover, 2-point crossover and uniform order crossover (UOX) were implemented and tested on some instances to choose the best one. For the SPR model with soft/flexible time windows, 2-point crossover operator gives better results.

  \subsubsection{Crossover operator}

  \begin{figure}[H]
  	
  	\begin{mdframed}[style=MyFrame,nobreak=true,align=center,userdefinedwidth=31em]
  		
  		\centering
  		
  		Parent 1: \hspace{3mm}
  		\begin{tabular}{ | c | c | c | c | c | c | c| }
  			\hline
  			5 (3)  &  3 (1) & 1 (2) & 2 (1) & 6 (2) &  4 (3)  \\ 
  			\hline
  			
  			2    &1 &  1 & 2    & 2        &   1\\ 
  			
  			\hline
  			
  		\end{tabular}
  		
  		\vspace{1mm}
  		Parent 2: \hspace{3mm}
  		\begin{tabular}{ | c | c | c | c | c | c | c| }
  			\hline
  			
  			4 (3) & 2 (1) & 3 (1) & 1 (2) &   6 (2) & 5 (3) \\ 
  			\hline
  			1 & 2 & 2    & 1   & 1  &   2\\   
  			
  			\hline
  			
  		\end{tabular}

  		\vspace{3mm}
  		
  		Offspring 1:
  		\begin{tabular}{ | c | c | c | c | c | c | c| }
  			\hline
  			4 (3) &  3 (1) &\color{blue} 1 (2) & \color{blue}2 (1) & \color{blue}6 (2) & 5 (3) \\ 
  			\hline
  			1   & \color{blue}1 &  \color{blue}1 & \color{blue} 2    & 1  &   2\\     
  			\hline
  		\end{tabular}
  		
  		\vspace{1mm}
  		Offspring 2:
  		\begin{tabular}{ | c | c | c | c | c | c | c| }
  			\hline
  			5 (3) & 2 (1)& \color{red}3 (1) & \color{red}1 (2) & \color{red}  6 (2) & 4 (3) \\ 
  			\hline
  			2 & \color{red}2 & \color{red}2    & \color{red}1   & 2  &   1\\    
  			\hline
  		\end{tabular}
  		
  	\end{mdframed}
  	\caption{Example of a parent and an offspring }
  	\label{fig:parentAndOffsrping}
  \end{figure}

  \vspace{-2mm}
  

  \begin{figure}[H]
  	
  	\begin{mdframed}[style=MyFrame,nobreak=true,align=center,userdefinedwidth=31em]

  		\centering
  		
  		Positions: $p_1 = 3$ and $p_2 = 5$

  		\vspace{3mm}
  		
  		Parent 1:\hspace{18 mm}
  		\begin{tabular}{ | c | c | c | c | c | c | c| }
  			\hline
  			5 (3)  &  3 (1) &\color{blue} 1 (2) & \color{blue}2 (1) & \color{blue}6 (2) &  4 (3)  \\ 
  			\hline
  		
  		\end{tabular}
  		\\
  			\vspace{1mm}
  		Parent 2:\hspace{18 mm}
  		\begin{tabular}{ | c | c | c | c | c | c | c| }
  			\hline
  			4 (3) & 2 (1) & \color{red}3 (1) & \color{red}1 (2) & \color{red}  6 (2) & 5 (3) \\ 
  			\hline

  		\end{tabular}
  		
  		\vspace{3mm}
  		
  		Step 2  (parent 1):\hspace{2mm}
  		\begin{tabular}{ | c | c | c | c | c | c | c| }
  			\hline
  			\hspace{8.7mm} & \hspace{8.7mm} &\color{blue} 1 (2) & \color{blue}2 (1) & \color{blue}6 (2) & \hspace{8.7mm} \\ 
  			\hline
  			
  		\end{tabular} 
  	  \\
  		\vspace{1mm}
  		Step 2   (parent 2):\hspace{2mm}
  		\begin{tabular}{ | c | c | c | c | c | c | c| }
  			\hline
  			\hspace{8.7mm} & \hspace{8.7mm} & \color{red}3 (1) & \color{red}1 (2) & \color{red}  6 (2) & \hspace{8.7mm} \\ 
  			\hline
  			
  		\end{tabular}
  		
  		\vspace{3mm}
  		
  		Offspring 1: \hspace{11 mm}
  		\begin{tabular}{ | c | c | c | c | c | c | c|c| }
  			\hline
  			4 (3) &  3 (1) &\color{blue} 1 (2) & \color{blue}2 (1) & \color{blue}6 (2) & 5 (3) \\ 
  			\hline
  		\end{tabular}
  		\\
  			\vspace{1mm}
  		Offspring 2: \hspace{11 mm}
  		\begin{tabular}{ | c | c | c | c | c | c | c|c| }
  			\hline
  			5 (3) & 2 (1)& \color{red}3 (1) & \color{red}1 (2) & \color{red}  6 (2) & 4 (3) \\ 
  			\hline
  			
  		\end{tabular}
  		
  	\end{mdframed}
  	\caption{Example of  patients' crossover  operator}
  	\label{fig:patientsCrossover}
  \end{figure}
  
  
    \vspace{-2mm}

  \begin{figure}[H]
  	
  	\begin{mdframed}[style=MyFrame,nobreak=true,align=center,userdefinedwidth=31em]
  		\begin{center}
  			
  		Positions: $p_1 = 2$ and $p_2 = 4$

  		\vspace{3mm}

  		Parent 1: \hspace{2mm}
  		\begin{tabular}{| c | c | c | c | c | c | c|}
  			
  			\hline
  			 \hspace{1.95mm} 2  \hspace{1.95mm}  &  \hspace{1.95mm} 1  \hspace{1.95mm} &  \hspace{1.95mm}  1  \hspace{1.95mm} &  \hspace{1.95mm} 2  \hspace{1.95mm}    &  \hspace{1.95mm} 2   \hspace{1.95mm} &  \hspace{1.95mm}  1  \hspace{1.95mm} \\ 
  			\hline
  		\end{tabular}
  	
  		\vspace{1mm}
  	    Parent 2: \hspace{2mm}
  		\begin{tabular}{ | c | c | c | c | c | c | c| }
  			
  			\hline
  			 \hspace{1.95mm}  1  \hspace{1.95mm}  &  \hspace{1.95mm} 2  \hspace{1.95mm} &  \hspace{1.95mm}  2  \hspace{1.95mm} &   \hspace{1.95mm} 1  \hspace{1.95mm} &   \hspace{1.95mm}  1  \hspace{1.95mm}  &  \hspace{1.95mm}   2  \hspace{1.95mm} \\   
  			\hline
  			
  		\end{tabular}

  		\vspace{3mm}
  		
  		Step 2 :\hspace{6mm}
  		\begin{tabular}{ | c | c | c | c | c | c | c| }
  			
  			\hline
  			\hspace{8.7mm}    & \color{blue}  \hspace{1.95mm} 1  \hspace{1.95mm}  &  \color{blue}  \hspace{1.95mm} 1  \hspace{1.95mm} & \color{blue}  \hspace{1.95mm} 2  \hspace{1.95mm}    & 	\hspace{8.7mm}   &   	\hspace{8.7mm} \\  
  			\hline
  		\end{tabular}
  	
  	  	\vspace{1mm}
  	
  	 	Step 2 :\hspace{6mm}
  		\begin{tabular}{ | c | c | c | c | c | c | c|c| }
  			
  			\hline
  				\hspace{8.7mm}  & \color{red}  \hspace{1.95mm} 2  \hspace{1.95mm}  & \color{red}  \hspace{1.95mm} 2  \hspace{1.95mm}    & \color{red}  \hspace{1.95mm} 1  \hspace{1.95mm}   & 	\hspace{8.7mm}   &  	\hspace{8.7mm} \\    
  			\hline
  		\end{tabular}

  		\vspace{3mm}
  		
  		offspring 1:
  		\begin{tabular}{ | c | c | c | c | c | c | c| }
  			
  			\hline
  			 \hspace{1.95mm} 1  \hspace{1.95mm}   & \color{blue}  \hspace{1.95mm} 1  \hspace{1.95mm} &  \color{blue}  \hspace{1.95mm} 1   \hspace{1.95mm} & \color{blue}  \hspace{1.95mm} 2  \hspace{1.95mm}  &  \hspace{1.95mm}  1  \hspace{1.95mm}  &  \hspace{1.95mm}   2  \hspace{1.95mm} \\    
  			\hline
  		\end{tabular}
  	
  		\vspace{1mm}
      	offspring 2:
  		\begin{tabular}{| c | c | c | c | c | c | c| }
  			
  			\hline
  			 \hspace{1.95mm} 2  \hspace{1.95mm} & \color{red}  \hspace{1.95mm} 2  \hspace{1.95mm} & \color{red}  \hspace{1.95mm} 2  \hspace{1.95mm}    & \color{red}  \hspace{1.95mm} 1  \hspace{1.95mm}   &  \hspace{1.95mm} 2  \hspace{1.95mm}  &  \hspace{1.95mm}   1  \hspace{1.95mm} \\    
  			\hline
  		\end{tabular}
  		
  	\end{center}
  	\end{mdframed}
  	\caption{Example of  caregivers'  crossover  operator}
  	\label{fig:caregiversCrossover}
  \end{figure}

   The crossover operation is the main genetic operator in \ac{GA} used to pass parents’ genes to their children. Many crossover operators have been proposed in the literature. In this study, we use the 2-point crossover operator to reproduce an offspring from two parents. It will be independently applied for each chromosome:
  \begin{itemize}
  	
  	\item Generate two random crossover points $p_1$ and $p_2$ in the parent; 
  	\item Copy the segment between points $p_1$ and $p_2$ from the first parent to the first offspring;
  	\item Copy the segment before  $p_1$ and the segment after $p_2$ from the second parent to the first offspring;
  	\item Repeat for the second offspring with the parent’s role reversed.
  	
  \end{itemize}

   The first two steps are the same for both chromosomes (patients and caregivers). However, the third step must be adapted to patients' chromosome. Indeed, the repetition of caregivers does not pose a problem when exchanging genes between parents  (see figures \ref{fig:parentAndOffsrping} and \ref{fig:caregiversCrossover}). Contrariwise, each patient must  appear only once in the solution, so the third step will be adapted by sorting  genes of the segment before  $p_1$ and the segment after $p_2$  in the same order as they appear at each opposite parent to avoid duplication or deletion of a patient from the solution (see figures \ref{fig:parentAndOffsrping} and \ref{fig:patientsCrossover}).

  
  
  \subsubsection{Fitness and selection }
 
  In this section, we use the same selection method (tournament selection) motivated and explained the previous section. 
  For each solution, three components are computed: caregivers’ transportation cost (see equation \ref{costFitness}); number of patients visited after exceeding the maximum number to visit per caregiver (see equation \ref{overVisits}); and the expected value for the SPR model (see equation \ref{expectedValue}). For the deterministic model, $F_2$ is equivalent to the number of unvisited patients. 
  Solutions are first compared according to the value of $F_1$, the individual with small value will be selected. If constraints \ref{eqn:MaxToVisitSMSS} are verified ($F_1=0$), solutions are compared according to the value of $F_2$. For the deterministic model this value must converge to 0. Contrariwise, for the SPR model, this value could be greater than 0. If solutions have the same values of $F_1$ and $F_2$, the solution with small transportation cost will be selected. We use this  lexicographic order to ensure convergence to feasible solutions ($F_1=0$) and to avoid using aggregation techniques for $F_2$ and $F_3$ since fixing weights is confusing and units are not the same. In addition, for any skipping a patient without providing the requested service operation, a wasting cost will  occurred. For example, if a caregiver will visit patients $p_1, p_2$ and $p_3$ and it happen that he will skip the patient $p_2$, the wasting cost occurred is $c_{p_1p_2} + c_{p_2p_3} - c_{p_1p_3}$ since he could be visit patient $p_3$ directly after patient $p_1$. Therefore, it interesting to minimize first the number of uninvited patients to ensure patients' satisfaction and avoid wasting costs.

  \begin{flalign} \label{overVisits}
  	F_{1} = \sum _{k=1}^c  max(\sum _{i=1} ^{n} \sum_{s=1}^q y_{iks} - Max_v,0)   &&
  \end{flalign}
  
  \begin{flalign} \label{expectedValue}
  	F_{2} = \mathbb{ E}[ \sum _{i=1}^n   v_{i} ] &&
  \end{flalign}
  \begin{flalign} \label{costFitness}
  	F_{3} =\sum _{k=1}^c\sum _{i=0}^n \sum _{j=1}^{n+1} c_{ij}x_{ijk}    &&
  \end{flalign}

  \subsubsection{Initial population}
  
  Given a population size $P_{size}$, the initial population for is randomly generated  as follows:
  
  \begin{itemize}
  	\item For $i =1$ to $P_{size}$ do:
  	\item Generate a random visiting order;
  	\item For each patient, assign a qualified caregiver selected randomly;
  	\item Compute components $F_1, F_2$ and $F_3$ (see equations \ref{overVisits}, \ref{expectedValue} and \ref{costFitness}).

  \end{itemize}

  \subsection{Numerical experiments} \label{NESMTW}

  \subsubsection{Test instances}

  The test instances have been randomly generated.   Deterministic and stochastic  parameters are generated respectively as described in \ref{TIDMTW} and  \ref{TISMSS}. Tuning parameters are fixed as follow: $P_c=0.4$, $P_m=0.08$, $P_{size}=n\times 20$, the stopping criterion for \ac{GA} based heuristic is fixed as the number of no improvement in the best solution found for $X$ iterations ($X=n \times 5$), and the tournament selection size is fixed to 2.
  
  Two sets of instances are generated. The first set is used with a single time window per patient and contains three subsets ( $A, B$ and $C$). The instance A1$\_$1 refers to the instance 1 of the category  A with single time window. Accordingly, the second set  has the same instances as the first set that are used with two availability periods  (see Table \ref{table:TIDSMTW}). The instance A1$\_$2 refers to the instance 1 of the category  A with double time windows. The same instances are used with single and double availability periods  to study their impact on solutions quality.

  \begin{table}[!htb]
  	\caption{ Tested instances details}
  	\centering
  	\def\arraystretch{1.25}
  	
  	\begin{tabular*}{\textwidth}{@{\extracolsep{\fill}} lllllll}
  		\hline
  		
  		Set & Subset & Size& Max$_v$ & N &L  & k    \\   
  		
  		\hline
  		& A$i$\_1  &		& 4& 10	&	1 & 3    \\   
  		
  		STW  & B$i$\_1 & $i \in \{1,2,...,7\}$	&8	& 25	&	1 & 5    \\   
  		
  		& C$i$\_1 &	&10	& 50	&	1 & 10   \\   
  		\hline
  		& A$i$\_2 &	&4	& 10	&	2 & 3   \\   
  		
  		MTW 	& B$i$\_2 &	$i \in \{1,2,...,7\}$&8	& 25	&	2	 & 5    \\   
  		
  		& C$i$\_2 &	&10	& 50	&	2 	& 10    \\

  		\hline

  	\end{tabular*}

  	\label{table:TIDSMTW}

  \end{table}

  \subsubsection{Computational results}

    Instances described above are solved within a time limit of 4 hours. 
  CPLEX could solve the subsets $A$ , $B$  and $C$  with single time window optimally except instances $C2$ and $C3$ for which a feasible solution is found with, respectively, a gap of $8.68\%$ and $6.85\%$. For multiple time windows,  the subsets $A$  are solved optimally and only instances $B1, B3, B4$ and $B6$ are solved optimally. For the instances $B2,B5$ and $B7$, a feasible  solution is found with, respectively, a gap of 10.78\%, 13.18\% and 14.04\%. Instances of the subset $C$ are hard to solve, CPLEX  is not able to resolve these instances within the time limit. This complexity is due to the multiple time windows, which is exponential. For $ n$ patients and $p$ availability periods for each one, we have  $p^n$ possibilities to select for each one a time window to receive care services. 
  
  The proposed \ac{GA} based heuristic could solve instances in  short computational CPU running times. For some instances, the optimal solutions are found, but for the others they are very close to those found by CPLEX (see figures  (\ref{STWD},\ref{MTWD}) and tables (\ref{STWD}, \ref{MTWD})). The complexity  due to the multiple time windows faced by CPLEX did not affected the \ac{GA} since we select for each patient the best time window independently of others patients, which avoid the exponential complexity. The worst CPU running time is on average 17.3 seconds and is elapsed to solve the instance $C5$ with double time windows. The worst gaps found for subsets with single time widows $A, B$ and $C$ equals respectively 1.13\% ($A6$), 6.27\% ($B2$) and 19.33\% ($C2$) (see Table \ref{STWD}). For instances with multiple time windows, the worst gaps found equals 1.16\% ($A1$), 21.01 \%($B2$)  and 47.44\% ($C2$). 
  
  The SPR model is solved by Mont Carlo simulation, which is used to estimate the expected value, embedded into th \ac{GA} based heuristic. For each instance, Caregivers' transportation cost and the average of unvisited patient is considered. The CPU time running time is significantly increasing for the SPR model compared the the deterministic version due to simulation that must be carried out for each new offspring to compute the expected value (see Table \ref{table:SPRM}). In addition, the expected value equals to zero for all instances except for instances $A3$, $A7$, $B2$ and $B6$ used with single time window, which show the robustness of the SPR model.
  
  A comparison is carried between instances with single and multiple time windows for both models. Figures \ref{STWvsMTWd} and \ref{STWvsMTWSPR} clearly show the advantage of adopting multiple time windows since the transportation cost is higher for instances with single availability period. Figures \ref{DvsSPRSTW} and \ref{DvsSPRMTW} show that caregivers' transportation cost is higher for the stochastic model because solutions with lower expected value are prioritized to avoid skipping patients without providing services operations and increase patients' satisfaction. Solutions with higher difference between patients' latest service times ($b_{il}$) and caregivers' completion $C_{ik}$ are prioritized to increase the chance of respecting  patients' time windows since travel and service times are supposed stochastic. The more the margin ($b_{il}$- $C_{ik}$) is higher the more the fluctuation of  caregivers' arrival times remains robust and compatible with patient's time windows. To illustrate that, we solved instance $A5$ using both models, the two solutions found are: for the deterministic model, caregivers' transportation cost is 662 (see tables \ref{table:ASTWD} and \ref{table:STWD}) and for the stochastic model is 744 (see tables \ref{table:ASTWSPR} and \ref{table:SPRM}). We computed starting and completion times for patients visited by the caregiver 1 for both models stochastic and deterministic. The minimum difference between the latest service time ($b_{i1}$) and the completion time $C_{i1}$ is 12 for the stochastic model and 4 for the deterministic model (see tables \ref{table:competionTimeD} and \ref{table:competionTimeSPR}).

  \begin{figure}[H]
  	\begin{center}
  		\begin{tikzpicture}
  			\begin{axis}[
  				width=\textwidth, 
  				height=7.5cm,
  				grid=major, 
  				grid style={dashed,gray!30},
  				xlabel= Instances, 
  				ylabel= Objective function values,
  				xtick=data,
  				legend style={at={(0.08,0.6)},anchor=west},
  				,xticklabels={A1,A2,A3,A4,A5,A6,A7,B1,B2,B3,B4,B5,B6,B7,C1,C2,C3,C4,C5,C6,C7}
  				]
  				
  				\addplot+[] 
  				table[x=STW, y=Z1, col sep=comma] {src/benchmarking.csv}; 
  				
  				\addplot+[] 
  				table[x=STW, y=Best1, col sep=comma] {src/benchmarking.csv}; 
  				
  				\addplot+[] 
  				table[x=STW, y=Worst1, col sep=comma] {src/benchmarking.csv}; 
  				
  				\addplot+[] 
  				table[x=STW, y=Average1, col sep=comma] {src/benchmarking.csv};

  				\legend{CPLEX,Best,Worst,Average}
  			\end{axis}
  			
  		\end{tikzpicture}
  		\vspace{-6mm}
  		\caption{CPLEX and \ac{GA} solutions comparison  for  instances of the set $STW$  with deterministic parameters.}
  		\label{STWD}

  		\begin{tikzpicture}
  			\begin{axis}[
  				width=\textwidth, 
  				height=7.5cm,
  				grid=major, 
  				grid style={dashed,gray!30},
  				xlabel= Instances, 
  				ylabel= Objective function values,
  				xtick=data,
  				legend style={at={(0.08,0.6)},anchor=west},
  				,xticklabels={A1,A2,A3,A4,A5,A6,A7,B1,B2,B3,B4,B5,B6,B7,C1,C2,C3,C4,C5,C6,C7}
  				]
  				
  				\addplot+[] 
  				table[x=STW, y=Z2, col sep=comma] {src/benchmarking.csv}; 
  				
  				\addplot+[] 
  				table[x=STW, y=Best2, col sep=comma] {src/benchmarking.csv}; 
  				
  				\addplot+[] 
  				table[x=STW, y=Worst2, col sep=comma] {src/benchmarking.csv}; 
  				
  				\addplot+[] 
  				table[x=STW, y=Average2, col sep=comma] {src/benchmarking.csv};

  				\legend{CPLEX,Best,Worst,Average}
  			\end{axis}
  			
  		\end{tikzpicture}
  		\vspace{-6mm}
  		\caption{CPLEX and \ac{GA} solutions comparison  for  instances of the set $MTW$  with deterministic parameters.}
  		\label{MTWD}

  		\begin{tikzpicture}
  			\begin{axis}[
  				width=\textwidth, 
  				height=7.5cm,
  				grid=major, 
  				grid style={dashed,gray!30},
  					xlabel= Instances, 
  				ylabel= Objective function values,
  				xtick=data,
  				legend style={at={(0.08,0.6)},anchor=west},
  				,xticklabels={A1,A2,A3,A4,A5,A6,A7,B1,B2,B3,B4,B5,B6,B7,C1,C2,C3,C4,C5,C6,C7}
  				]
  				
  				\addplot+[] 
  				table[x=STW, y=Z1, col sep=comma] {src/benchmarking.csv}; 
  				
  				\addplot+[] 
  				table[x=STW, y=Z2, col sep=comma] {src/benchmarking.csv};

  				\legend{STW,MTW}
  			\end{axis}
  			
  		\end{tikzpicture}
  		\vspace{-6mm}
  		\caption{Comparison of the best-found solutions for the deterministic model according to time windows}
  		\label{STWvsMTWd}
  		
  	\end{center}
  \end{figure}	
  
  \begin{figure}[H]
  	\begin{center}
  		
  		\begin{tikzpicture}
  			\begin{axis}[
  				width=\textwidth, 
  				height=7.5cm,
  				grid=major, 
  				grid style={dashed,gray!30},
  				xlabel= Instances, 
  				ylabel= Objective function values,
  				xtick=data,
  				legend style={at={(0.08,0.6)},anchor=west},
  				,xticklabels={A1,A2,A3,A4,A5,A6,A7,B1,B2,B3,B4,B5,B6,B7,C1,C2,C3,C4,C5,C6,C7}
  				]
  				
  				\addplot+[] 
  				table[x=STW, y=COST1, col sep=comma] {src/benchmarking.csv}; 
  				
  				\addplot+[] 
  				table[x=STW, y=COST2, col sep=comma] {src/benchmarking.csv};

  				\legend{STW,MTW}
  			\end{axis}
  			
  		\end{tikzpicture}
  		\vspace{-6mm}
  		\caption{Comparison of the best-found solutions for the SPR model according to time windows}
  		\label{STWvsMTWSPR}

  		\begin{tikzpicture}
  			\begin{axis}[
  				width=\textwidth, 
  				height=7.5cm,
  				grid=major, 
  				grid style={dashed,gray!30},
  				xlabel= Instances, 
  				ylabel= Objective function values,
  				xtick=data,
  				legend style={at={(0.08,0.6)},anchor=west},
  				,xticklabels={A1,A2,A3,A4,A5,A6,A7,B1,B2,B3,B4,B5,B6,B7,C1,C2,C3,C4,C5,C6,C7}
  				]
  				
  				\addplot+[] 
  				table[x=STW, y=Z1, col sep=comma] {src/benchmarking.csv}; 
  				
  				\addplot+[] 
  				table[x=STW, y=COST1, col sep=comma] {src/benchmarking.csv};

  				\legend{Deterministic model, SPR model}
  			\end{axis}
  			
  		\end{tikzpicture}
  		\vspace{-6mm}
  		\caption{Comparison of the best-found solutions  with single time window according to type of model }
  		\label{DvsSPRSTW}

  		\begin{tikzpicture}
  			\begin{axis}[
  				width=\textwidth, 
  				height=7.5cm,
  				grid=major, 
  				grid style={dashed,gray!30},
  				xlabel= Instances, 
  				ylabel= Objective function values,
  				xtick=data,
  				legend style={at={(0.08,0.6)},anchor=west},
  				,xticklabels={A1,A2,A3,A4,A5,A6,A7,B1,B2,B3,B4,B5,B6,B7,C1,C2,C3,C4,C5,C6,C7}
  				]
  				
  				\addplot+[] 
  				table[x=STW, y=Z2, col sep=comma] {src/benchmarking.csv}; 
  				
  				\addplot+[] 
  				table[x=STW, y=COST2, col sep=comma] {src/benchmarking.csv};

  				\legend{Deterministic model, SPR model}
  			\end{axis}
  			
  		\end{tikzpicture}
  		\vspace{-6mm}
  		\caption{Comparison of the best-found solutions  with multiple time windows according to type of model }
  		\label{DvsSPRMTW}

  	\end{center}
  \end{figure}

 \begin{table}[!htb]
	
	\caption{Numerical results of tested instances with single time window and deterministic parameters}
	\centering
	\renewcommand{\arraystretch}{1.25}
	
	\begin{tabular*}{\textwidth}{@{\extracolsep{\fill}} c | cccc|ccccc }
		\hline
		\multicolumn{1}{c|}{Instances} & \multicolumn{4}{c|}{\textbf{CPLEX} }&\multicolumn{5}{c}{ \textbf{\ac{GA} (10 runs)}} \\ 
		\hline

		\textbf{STW }     & \textbf{LB} & \textbf{Z}&\textbf{ Gap} & \textbf{CPU} & \textbf{Best} & \textbf{Worst}& \textbf{Average} &\textbf{ Gap} & \textbf{CPU }   \\  
		
		A1\_1	&	525.00	&	\textbf{525}	&	0.00\%	&	1.43	&	\textbf{525}	&	566	&	529.10	&	0.77\%	&	$<1$	\\
		A1\_2	&	754.00	&	\textbf{754}	&	0.00\%	&	1.28	&	\textbf{754}	&	754	&	754.00	&	0.00\%	&	$<1$	\\
		A1\_3	&	588.00	&	\textbf{588}	&	0.00\%	&	1.31	&	\textbf{588}	&	588	&	588.00	&	0.00\%	&	$<1$	\\
		A1\_4	&	817.00	&	\textbf{817}	&	0.00\%	&	1.43	&	\textbf{817}	&	817	&	817.00	&	0.00\%	&	$<1$	\\
		A1\_5	&	662.00	&	\textbf{662}	&	0.00\%	&	1.62	&	\textbf{662}	&	677	&	665.40	&	0.51\%	&	$<1$	\\
		A1\_6	&	439.00	&	\textbf{439}	&	0.00\%	&	1.37	&	\textbf{439}	&	464	&	444.00	&	1.13\%	&	$<1$	\\
		A1\_7	&	539.00	&	\textbf{539}	&	0.00\%	&	1.36	&	\textbf{539}	&	539	&	539.00	&	0.00\%	&	$<1$	\\
		
		\hline
		B1\_1	&	1165.00	&	\textbf{1165}	&	0.00\%	&	2.64	&	\textbf{1165}	&	1300	&	1232.50	&	5.48\%	&	8.34	\\
		B1\_2	&	993.00	&	\textbf{993}	&	0.00\%	&	2.40	&	\textbf{993}	&	1108	&	1059.40	&	6.27\%	&	8.71	\\
		B1\_3	&	1089.00	&	\textbf{1089}	&	0.00\%	&	3.58	&	\textbf{1089}	&	1162	&	1124.40	&	3.15\%	&	8.76	\\
		B1\_4	&	928.00	&	\textbf{928}	&	0.00\%	&	3.56	&	\textbf{928}	&	1001	&	965.00	&	3.83\%	&	9.12	\\
		B1\_5	&	1064.00	&	\textbf{1064}	&	0.00\%	&	2.27	&	\textbf{1064}	&	1246	&	1110.60	&	4.20\%	&	10.51	\\
		B1\_6	&	1196.00	&	\textbf{1196}	&	0.00\%	&	3.04	&	\textbf{1196}	&	1276	&	1226.40	&	2.48\%	&	10.84	\\
		B1\_7	&	1099.00	&	\textbf{1099}	&	0.00\%	&	1.80	&	\textbf{1099}	&	1172	&	1116.40	&	1.56\%	&	9.02	\\
		
		\hline
		
		C1\_1	&	1490.00	&	\textbf{1490}	&	0.00\%	&	411.76	&	1555	&	1719	&	1644.80	&	9.41\%	&	92.70	\\
		C1\_2	&	1529.64	&	1675	&	8.68\%	&	14400	&	1799	&	1996	&	1896.20	&	19.33\%	&	120.54	\\
		C1\_3	&	1425.26	&	1530	&	6.85\%	&	14400	&	1633	&	1853	&	1723.00	&	17.28\%	&	80.62	\\
		C1\_4	&	1603.00	&	\textbf{1603}	&	0.00\%	&	5389.00	&	1803	&	1998	&	1889.70	&	15.17\%	&	102.29	\\
		C1\_5	&	1482.00	&	\textbf{1482}	&	0.00\%	&	388.24	&	1614	&	1974	&	1752.70	&	15.44\%	&	141.29	\\
		C1\_6	&	1658.00	&	\textbf{1658}	&	0.00\%	&	172.16	&	1736	&	1999	&	1848.30	&	10.30\%	&	92.03	\\
		C1\_7	&	1492.00	&	\textbf{1492}	&	0.00\%	&	232.13	&	1608	&	1742	&	1675.00	&	10.93\%	&	91.38	\\

		\hline
	\end{tabular*}

	\label{table:STWD}
	
\end{table}

\begin{table}[!htb]
	
	\caption{Numerical results of tested instances with multiple time windows and  deterministic parameters}
	\centering
	\renewcommand{\arraystretch}{1.25}
	
	\begin{tabular*}{\textwidth}{@{\extracolsep{\fill}} c | cccc|ccccc}
		\hline
		\multicolumn{1}{c|}{Instances} & \multicolumn{4}{c|}{\textbf{CPLEX} }&\multicolumn{5}{c}{ \textbf{\ac{GA} (10 runs)}} \\ 
		\hline

		\textbf{MTW }     & \textbf{LB} & \textbf{Z}&\textbf{ Gap} & \textbf{CPU} & \textbf{Best} & \textbf{Worst}& \textbf{Average} &\textbf{ Gap} & \textbf{CPU }   \\  
		\hline
		
		A2\_1	&	512.00	&	\textbf{512}	&	0.00\%	&	1.76	&	\textbf{512}	&	524	&	518.00	&	1.16\%	&	$<1$	\\
		A2\_2	&	672.00	&	\textbf{672}	&	0.00\%	&	2.01	&	\textbf{672}	&	672	&	672.00	&	0.00\%	&	$<1$	\\
		A2\_3	&	588.00	&	\textbf{588}	&	0.00\%	&	1.74	&	\textbf{588}	&	593	&	588.60	&	0.10\%	&	$<1$	\\
		A2\_4	&	665.00	&	\textbf{665}	&	0.00\%	&	2.15	&	\textbf{665}	&	673	&	666.60	&	0.24\%	&	$<1$	\\
		A2\_5	&	593.00	&	\textbf{593}	&	0.00\%	&	2.86	&	\textbf{593}	&	598	&	594.90	&	0.32\%	&	$<1$	\\
		A2\_6	&	388.00	&	\textbf{388}	&	0.00\%	&	1.86	&	\textbf{388}	&	388	&	388.00	&	0.00\%	&	$<1$	\\
		A2\_7	&	507.00	&	\textbf{507}	&	0.00\%	&	1.99	&	\textbf{507}	&	507	&	507.00	&	0.00\%	&	$<1$	\\
		
		\hline
		
		B2\_1	&	855.00	&	\textbf{855}	&	0.00\%	&	403.83	&	935	&	1030	&	989.90	&	13.63\%	&	10.26	\\
		B2\_2	&	665.56	&	746	&	10.78\%	&	14400	&	759	&	881	&	842.60	&	21.01\%	&	10.85	\\
		B2\_3	&	905.00	&\textbf{	905}	&	0.00\%	&	4027.58	&	947	&	1034	&	985.50	&	8.17\%	&	10.59	\\
		B2\_4	&	739.00	&\textbf{	739}	&	0.00\%	&	512.62	&	743	&	852	&	782.50	&	5.56\%	&	9.60	\\
		B2\_5	&	726.67	&	837	&	13.18\%	&	14400	&	853	&	941	&	894.30	&	18.74\%	&	12.36	\\
		B2\_6	&	944.00	&	\textbf{944}	&	0.00\%	&	1065.68	&	950	&	1043	&	998.20	&	5.43\%	&	11.06	\\
		B2\_7	&	810.57	&	943	&	14.04\%	&	14400	&	948	&	1047	&	993.00	&	18.37\%	&	10.93	\\
		
		\hline
		
		C2\_1	&	800.04	&	-	&	-	&	14400	&	1330	&	1549	&	1406.80	&	43.13\%	&	125.68	\\
		C2\_2	&	834.54	&	-	&	-	&	14400	&	1467	&	1644	&	1587.80	&	47.44\%	&	162.22	\\
		C2\_3	&	999.36	&	-	&	-	&	14400	&	1437	&	1601	&	1519.60	&	34.24\%	&	124.98	\\
		C2\_4	&	860.14	&	-	&	-	&	14400	&	1474	&	1721	&	1622.20	&	46.98\%	&	138.26	\\
		C2\_5	&	812.12	&	-	&	-	&	14400	&	1367	&	1582	&	1453.30	&	44.12\%	&	173.21	\\
		C2\_6	&	842.53	&	-	&	-	&	14400	&	1417	&	1598	&	1505.03	&	44.02\%	&	138.54	\\
		C2\_7	&	842.50	&	-	&	-	&	14400	&	1361	&	1517	&	1436.30	&	41.34\%	&	143.47	\\

		\hline
	\end{tabular*}

	\label{table:MTWD}
	
\end{table}

\begin{table}[!htb]
	
	\caption{Numerical results of tested instances with stochastic parameters}
	\centering
	\renewcommand{\arraystretch}{1.25}
	
	\begin{tabular*}{\textwidth}{@{\extracolsep{\fill}} c | cccc|  cccc }
		\hline
		\multicolumn{1}{c|}{} & \multicolumn{4}{c|}{\textbf{single time window} }&\multicolumn{4}{c}{ \textbf{Multiple time windows}} \\ 
		\hline

		\textbf{Instances}    & \textbf{L}& \textbf{COST} & \textbf{E()}& \textbf{CPU} & \textbf{L}  & \textbf{COST} & \textbf{E()}& \textbf{CPU}   \\  
		\hline
		
		A1	&	1	&	525	&	0.00	&	2.32	&	2	&	524	&	0.00	&	2.29	\\
		A2	&	1	&	754	&	0.00	&	2.30	&	2	&	672	&	0.00	&	2.28	\\
		A3	&	1	&	609	&	\textbf{0.01}	&	3.73	&	2	&	594	&	0.00	&	2.40	\\
		A4	&	1	&	817	&	0.00	&	2.44	&	2	&	673	&	0.00	&	3.81	\\
		A5	&	1	&	744	&	0.00	&	3.15	&	2	&	593	&	0.00	&	2.36	\\
		A6	&	1	&	487	&	0.00	&	1.94	&	2	&	388	&	0.00	&	2.43	\\
		A7	&	1	&	586	&	\textbf{0.03}	&	2.51	&	2	&	507	&	0.00	&	2.48	\\
		\hline
		B1	&	1	&	1271	&	0.00	&	50.03	&	2	&	1149	&	0.00	&	68.24	\\
		B2	&	1	&	1274	&	\textbf{0.05}	&	51.01	&	2	&	942	&	0.00	&	54.80	\\
		B3	&	1	&	1261	&	0.00	&	94.36	&	2	&	1070	&	0.00	&	52.63	\\
		B4	&	1	&	1131	&	0.00	&	42.91	&	2	&	934	&	0.00	&	55.70	\\
		B5	&	1	&	1151	&	0.00	&	56.03	&	2	&	941	&	0.00	&	60.28	\\
		B6	&	1	&	1413	&	\textbf{0.16}	&	64.40	&	2	&	1104	&	0.00	&	57.36	\\
		B7	&	1	&	1285	&	0.00	&	39.99	&	2	&	1111	&	0.00	&	144.98	\\
		\hline
		C1	&	1	&	1696	&	0.00	&	751.00	&	2	&	1608	&	0.00	&	822.16	\\
		C2	&	1	&	2058	&	0.00	&	643.87	&	2	&	1774	&	0.00	&	666.43	\\
		C3	&	1	&	1954	&	0.00	&	1406.80	&	2	&	1563	&	0.00	&	596.02	\\
		C4	&	1	&	2030	&	0.00	&	590.64	&	2	&	1913	&	0.00	&	1083.88	\\
		C5	&	1	&	1917	&	0.00	&	1225.75	&	2	&	1546	&	0.00	&	961.67	\\
		C6	&	1	&	1772	&	0.00	&	906.39	&	2	&	1659	&	0.00	&	784.86	\\
		C7	&	1	&	1849	&	0.00	&	644.04	&	2	&	1623	&	0.00	&	652.07	\\

		\hline
	\end{tabular*}

	\label{table:SPRM}
	
\end{table}

\begin{table}[H]
	
	\caption{Solution of instance $A5$ with single time window and deterministic parameters}
	\centering
	\centering
	\begin{tabular}{  c | c | c | c | c | c |c|c | c | c |c }
		\hline
		Patients  & 5(6)&	2(6)&	9(5)&	8(5)&	10(2)&	6(6)&	3(5)&	1(2)&	7(5)&	4(5)	 \\ 
		\hline
		Caregivers &  1&3&	3&	1&	2&	3&	1&	2&	1&	3  \\   
		\hline
	\end{tabular}
	
	\label{table:ASTWD}
	
\end{table}

\begin{table}[H]
	
	\caption{Solution of instance $A5$ with single time window and stochastic parameters}
	\centering
	\centering
	\begin{tabular}{  c | c | c | c | c | c |c|c | c | c |c }
		\hline
		Patients  & 2(6)&	3(5)&	6(6)&	10(2)&	7(5)&	1(2)&	9(5)&	5(6)&	4(5)&	8(5)
		\\ 
		\hline
		Caregivers &  1&	3&	3&	2&	3&	2&	1&	1&	3&	1	 \\   
		\hline
	\end{tabular}
	
	\label{table:ASTWSPR}
	
\end{table}

\begin{table}[H]
	
	\caption{Caregiver 1 starting and completion times for assigned patients considering the solution  found by the deterministic model for instances $A5$ with single time window}
	\centering
	\centering
	\begin{tabular}{  l | c | c | c | c | c |c|c | c | c }
		\hline
		Patients  &\multicolumn{2}{c|}{ 5(6)}	& \multicolumn{2}{c|}{8(5)}	&	\multicolumn{2}{c|}{3(5)}	&	\multicolumn{2}{c|}{7(5)}\\	
		\hline
		Time windows $[a_{i1}, b_{i1}]$ & 200&	320&	275&	395&	222&	342&	425&	545
		\\   
		\hline
		Start and completion times $ [S_{i1}, C_{i1}]$& 200&	220&	277&	294&	318&	338&	425&	444
		\\
		\hline
		$b_{i1} - C_{i1} $ & &100	&&	101	&&	4	&&	101
		\\
		\hline
	\end{tabular}
	
	\label{table:competionTimeSPR}
	
\end{table}

\begin{table}[H]
	
	\caption{Caregiver 1 starting and completion times for assigned patients considering the solution  found by the SPR model for instances $A5$ with single time window}
	\centering
	\begin{tabular}{  l | c | c | c | c | c |c|c | c | c }
		\hline
		Patients  &  \multicolumn{2}{c|}{2(6)}	&	\multicolumn{2}{c|}{9(5)}	&	\multicolumn{2}{c|}{5(6)}&	\multicolumn{2}{c|}{8(5)}		\\ 
		\hline
		Time windows $[a_{i1}, b_{i1}]$ & 130&	250&	203&	323&	200&	320&	275&	395 \\   
		\hline
		Start and completion times $ [S_{i1}, C_{i1}]$& 130&	146&	230&	260&	288&	308&	365&	382\\
		\hline
		$b_{i1} - C_{i1} $ & &104&	&	63&  & 		12& & 13 \\
		\hline
	\end{tabular}
	
	\label{table:competionTimeD}
	
\end{table}

  \section{Conclusion}\label{SMConclusion}
  
The home health care companies aim to both minimize provided services cost and maximize patients' satisfaction. In the real world, travel and service times are not always deterministic.  Uncertainties may arise and affect the overall planning and service quality would be quite poor, which will cause patients’ dissatisfaction. These  companies seek to respect as much as possible patients’ time windows, which it becomes  challenging when dealing with stochastic travel or/and service times.  Two possible types of recourse can be used to deal with both time windows and the uncertainty of travel and service times. Either accepting providing services with a tardiness and time windows must be soft/flexible, or skipping a visit when a route becomes infeasible and time windows can be hard/fixed. In this chapter, the two types of recourse are used to deal stochastic service and travel times and time windows.

We proposed two-stage stochastic programming models with recourse. In the first model, we  suppose that patients can request multiple services with a possible synchronization if the services must simultaneously provided. Considering both the uncertainty of parameters and synchronization of multiple services increase the chance that time windows will not be respected. To overcome this issue, we  introduce the recourse if a solution goes against constraints related to patients’ time windows and caregivers’ duty length. The objective is to minimize the transportation cost and the expected value of recourse caused by  patients' delayed services and caregivers' extra working time.

In the second model, we assume that time windows are hard/fixed, which must be respected without any earliness or tardiness.  The recourse is defined as skipping patients when their time windows will be violated to guarantee that the services are provided without an earliness or tardiness. In addition, to increase the chance of providing the maximum of  services, we suppose also that patients are allowed to be available in many  periods. The objective is to minimize the transportation cost and the expected value of recourse, which expresses the average number of un-visited patients.  
  
The deterministic models are solved using CPLEX, the \ac{GA} and the  \ac{GVNS} based heuristics. \ac{GVNS}  and \ac{GA} are successfully tested using several instances randomly generated  from the literature. The tests prove the high performance of these two heuristics to deal with large instances in a little amount of time. \ac{GVNS}  and \ac{GA} are  able to reach optimal solutions for some instances  and yield near-optimal solutions for others. 

Monte Carlo simulation is  used to estimate the expected value of recourse, which is embedded into the \ac{GA} based heuristic to solve the proposed stochastic models. Computational results showed that \ac{GVNS}  is not suitable to be combined with the simulation to solve large instances. The complexity of the SPR model in terms of CPU running times is significant due to the expected value  that has to be estimated for each solution.

Future works could be addressed to propose a robust optimization approach to deal with the two introduced problems since this approach does not require knowing the distribution of stochastic parameters. Dynamic programming would be necessary to avoid revising the overall planning in the event that a requested service is canceled or a new risky service should be provided.

    \chapter{Multi-objective home health care routing and scheduling problem} \label{ChapterMOMTW}
    \vspace{-0.6cm}

   \section{Introduction}
  
     In this chapter, we study the  \ac{HHCRSP} considering multi-objectives that must be simultaneously optimized.  Most studies transform multi-objective problems into a mono-objective case using aggregation techniques. Assigning weights to conflicting  objectives is a confusing task and requires the decision maker' experience and knowledge of the problem. In addition, a sensitivity analysis  of weights should be conducted, which makes the process of solving multi-objective problems more complicated. The lexicographical order approach is based on a priori knowledge of some decision-maker preferences to establish the lexicographic order, which is not a simple task. 
  
  \ac{MOEAs} have been very popular in solving multi-objective problems (MOP). \ac{MOEAs} could obtain multiple Pareto front solutions in a single simulation run as they are a population based. Several \ac{MOEAs} have been proposed in the literature, which can be broadly grouped under three categories:

  \begin{itemize}
  	
  	\item  Domination-based: in this group, the selection strategy is based on Pareto dominance in order to  ensure convergence. Furthermore, an explicit diversity preservation scheme  must be used to maintain the diversity of solutions belonging to the same group. The most known \ac{MOEAs} of this group  are \ac{NSGA-II} \cite{deb2002fast} and strength pareto evolutionary algorithm (SPEA2) \cite{zitzler2001spea2}; 
  	
  	\item Indicator-based: this group  measures the fitness of solutions  with a performance indicator such as hypervolume indicator. $S$ metric selection evolutionary multiobjective optimization algorithm (SMS-EMOA) \cite{emmerich2005emo} and simple indicator based evolutionary algorithm (SIBEA) \cite{zitzler2007hypervolume} are famous representative of \ac{MOEAs} of this group;
  	
  	\item Decomposition-based: The methods of this group are based on the scalarization techniques to convert a  multi-objective problem  into a mono-objective optimization sub-problems, which are solved simultaneously. \ac{MOEA/D} \cite{zhang2007moea} and multi-objective genetic local search algorithm (MOGLS) \cite{ishibuchi1998multi} are the most popular algorithms of this category. 
  	
  \end{itemize}

  	Few studies \cite{braekers2016bi,decerle2019memetic,fathollahi2020bi} have dealt with the multi-objectives case in the home health care context by using methods based on the concept of Pareto dominance to approximate the Pareto front using methods of domination-based category. To our knowledge, methods from indicator-based and decomposition-based are not applied in the home health care context. In this chapter, we adopt the Pareto and decomposition based approaches to deal with the multi-objective HHCRSP.  According to \cite{zhang2007moea}, it is very time-consuming, if not impossible, to obtain the entire Pareto front for most multi-objective problems. In addition, the decision maker may not be interested in having  all Pareto optimal  solutions. \ac{MOEA/D} algorithm solves  multi-objective problems  faster since ranking solutions based on Pareto dominance is not used. This study deals  with three objectives so \ac{NSGA-II} algorithm  remains efficient for this kind of optimization problems. According to \cite{trivedi2016survey}, efficiently combining dominance and decomposition based approaches can result in high performance many objective optimizers. 
  	
  	Furthermore, most existing studies only consider a single availability period  per patient and no work, as far as we know, has been proposed for  multi-objective problems  with multiple time windows  in the \ac{HHC} context. This study considers multiple time windows for patients that aim to ensure their availability and give the decision maker more flexibility to schedule patients' visits, which is not possible with a single time window per patient. In \cite{bazirha2019daily}, instances with single and multiple availability periods were compared according to the sum of  earliness and tardiness of services operations as well as caregivers' waiting times. Computational results showed  that instances with multiple availability periods  are better optimized than those with a single period.  Based on this, it is worth trying multiple time windows when optimizing multi-objective problems to improve the approximate Pareto set, which is proven on tested instances.

     This chapter is  divided as follows. The definition of the problem studied is described in section \ref{MultiPMTW}. In section \ref{MOEA}, MOEA/D, \ac{NSGA-II} and the hybrid algorithm used to solve the problem are presented. The test instances, the experimental settings and the performance of the algorithms are measured by analyzing the experimental results in section \ref{MultiNE}. This chapter ends with a conclusion in section \ref{MultiC}.

    \section{Multi-objective home health care routing and scheduling problem with multiple  time windows}
  
  \subsection{Problem statement}\label{MultiPMTW}
  
  
  In this chapter, we aim to minimize caregivers' traveling times ($f_1$), their waiting times ($f_2$) and to balance their workload ($f_3$). The goal is to simultaneously optimize these objectives using \ac{MOEAs} while ensuring that the selected period per patient, caregivers' time windows and skills requirements are respected. The decision maker is supposed to be involved a posteriori in order to select a final preferred solution from the approximate set of Pareto's optimal solutions.

  \begin{flalign} 
  \min   f_1= \sum _{i=0}^n \sum _{j=1}^{n+1} T_{ij}\sum_{k=1}^cx_{ijk}   && \nonumber
  \end{flalign}
  
  \begin{flalign} 
  \min    f_2= \sum_{i=1}^n \sum_{k=1}^c(S_{ik}-A_{ik})  && \nonumber
  \end{flalign}
  
  \begin{flalign} 
  \min   f_3= \sum_{k=1}^c D_k  && \nonumber
  \end{flalign}

  \subsection{Solution approaches:  multi-objective evolutionary algorithms }\label{MOEA}
  
   The encoding/decoding of a solution is described in chapter \ref{ChapterDM}. Crossover and mutation operators are described in chapter \ref{ChapterSM} more precisely in section \ref{SectionSMS}.

   \subsubsection{Dominance}
  
  The notion of Pareto dominance is introduced by Pareto \cite{pareto1964cours} to determine  the set of non-dominated solutions. Considering $m$ as the number of objective functions, Pareto dominance is defined as follows:

  \begin{definition}{Dominance:}
  	A solution $x$ dominates ($\preceq$) a solution $y$ if and only if the following conditions are verified :
  	\begin{itemize}
  		\item $ \forall  i \in \{1, ..., m\}: f_{i}(x)  \leq f_{i}(y)$;
  		\item $  \exists i \in \{1, ..., m\}: f_{i}(x)  < f_{i}(y)$.
  	\end{itemize}
  	
  \end{definition}
  
  Evolutionary algorithms are designed to solve unconstrained problems. \cite{deb2002fast} proposed constrained-domination to transform a constrained problem to an unconstrained one and is defined as follows:
  
  \begin{definition}{Constrained-domination:} 
  	A solution $x$ is said to constrained-dominate a solution  $y$, if one of the following conditions is verified:
  	\begin{enumerate}
  		\item Solution $x$ is feasible and solution $y$ is not;
  		\item Solutions $x$ and $y$ are both infeasible and the  solution $x$ has a smaller overall constraint violation;
  		\item Solutions $x$ and $y$ are feasible and $x$ dominates $y$.
  	\end{enumerate}
  	
  	The solution generated by crossover and mutation operators may sometimes not be able to satisfy the constraints of patients' availability periods, caregivers' duty length or skills requirements. The infeasible assignments could be remedied by randomly replacing unqualified caregivers with skilled caregivers.  The constrained-domination is used to deal with  caregivers' duty length and patients' time windows that cannot be met. We ultimately use a penalty's cost related to the tardiness of services operations and  caregivers' overtime. 
  	
  \end{definition}
  
  \subsubsection{Initial population}\label{IP}
  
  The initial population for a given size $P_{size}$ is randomly generated as follows:
  
  \begin{itemize}
  	\item For $i =1$ to $P_{size}$ do:
  	\item   Generate a random visiting order (see Table \ref{table:3});
  	
  	\begin{table}[ htb]
  		
  		\caption{Example of visiting order}
  		\centering
  		\centering
  		\begin{tabular}{ | c | c | c | c | c | c | c|c| }
  			\hline
  			Patients  & 5 (2)  & 1 (3)  & 3 (1) &1 (2)& 4 (3) & 2 (1) &  6 (2)  \\ 
  			\hline
  			Caregivers  &    & &  & &     &     &   \\   
  			\hline
  		\end{tabular}
  		
  		\label{table:3}
  		
  	\end{table}
  	
  	\item For each patient, assign a qualified caregiver randomly selected (see Table \ref{table:4});
  	\begin{table}[!htb]
  		
  		\caption{example of caregivers’ assignment to patients}
  		\centering
  		\begin{tabular}{ | c | c | c | c | c | c | c|c| }
  			\hline
  			Patients  & 5 (2)  & 1 (3)  & 3 (1) &1 (2)& 4 (3) & 2 (1) &  6 (2)  \\ 
  			\hline
  			Caregivers  & 1    & 2 &  2 & 1   & 1    & 2     &   1\\   
  			\hline
  		\end{tabular}

  		\label{table:4}
  		
  	\end{table}
  	
  	\item Calculate the objective function values $f_1, f_2$ and $f_3$.

  \end{itemize}

  \subsubsection{Decomposition based: MOEA/D}
  The \ac{MOEA/D} is proposed by \cite{zhang2007moea}. It decomposes a multi-objective optimization problem into a number  of sub-problems  based on the distances between their aggregation coefficient vectors and optimizes them simultaneously. The optimal solutions to two neighboring sub-problems should be very similar as their coefficient vectors are very close. According to \cite{zhang2007moea} it is very time-consuming, if not impossible, to obtain the complete Pareto front for most multi-objective problems. Furthermore, the decision maker may not be interested in having  all Pareto optimal  solutions. Three scalarization techniques are used with \ac{MOEA/D} in \cite{zhang2007moea}: Weighted Sum Approach, Tchebycheff Approach and Boundary Intersection (BI) Approach. Tchebycheff Approach is adopted to solve our problem and is defined as follows:
  \par
  \begin{equation}
  \text{minimize  } \qquad  g^{te}(x| \lambda, z^{*})= \max_{1 \leq i \leq m}\{ \lambda_i| f_i(x) - z^*| \}
  \end{equation}
  
  \begin{equation}
  \text{subject to }\qquad x \in \Omega
  \end{equation}
  
  \par
  where $z^*=(z^*_1,...,z^*_m)^T $ is the reference point with $ z^*_i=min \{ f_i(x) | x \in \Omega \} $.
  
  The  \ac{MOEA/D} framework is described by Algorithm \ref{MOEA/D}. At each generation $t$, \ac{MOEA/D} with the Tchebycheff approach maintains :
  \begin{itemize}
  	\item  a population of N points $x^{1},......., x^{N} \in \Omega$, where $x^{i}$   is the current solution of the $i$ th sub-problem; 
  	\item a set of F-values  $FV^{1},......., FV^{N},$ where  $FV^{i}$ is the F-value at $x^{i}$, i.e., $FV^{i}=F(x^{i})$; 
  	\item A vector  $z=(z_{1},.....,z_{m})^{T}$, where $z_{i}$ is the best value found so far for the objective $f_{i}$; 
  	\item  an external population (EP), which is used to store non-dominated solutions found during the search.
  \end{itemize}

  \begin{algorithm}[!htb]
  	\SetAlgoLined
 
  	\textbf{Input:}\\
  	a stopping criterion \;
  	N the number of the sub-problems considered in \ac{MOEA/D} \;
  	a uniform spread of N weight vectors: $\lambda^{1},......., \lambda^{N}$ \;
  	T: the number of the weight vectors in the neighborhood of each weight vector\;  
  	\textbf{ Output: } EP \;	
  	
  	\textbf{Initialization:} \;
  	\textbf{Set} \textbf{EP}=$\emptyset$ \;
  	\For{$i\leftarrow 1$ \KwTo $N$}{
  		\textbf{set} $B(i)=\{i_{1},.....,i_{T}\}$, where  $\lambda^{i_{1}},......., \lambda^{i_{T}}$ are the T closest	weight vectors to $\lambda^{i}$ \;
  	}
  	\textbf{Generate} an initial population $x^{1},......., x^{N}$\;
  	\textbf{Set} $FV^{i}=F(x^{i})$\;
  	\textbf{Initialize}  $z=(z_{1},.....,z_{m})^{T}$ by a problem-specific method \;
  	
  	\While{(the stopping condition is not reached)}{
  		\For{$i\leftarrow 1$ \KwTo $N$}  {
  			\textbf{Select} Randomly two indexes $k$ and $l$ from $B(i)$ \;
  			$ y \longleftarrow$ \textbf{Crossover}($x^k$,$x^l$) \;
  			$ y' \longleftarrow$ \textbf{Mutation}($y$) \;
  			\textbf{Repair} ($y'$) \;
  			
  			\For{$j\leftarrow 1$ \KwTo $m$}  {

  				\If{$z_{j}< f_{j}(y')$ }{
  					\textbf{set} $ z_j \longleftarrow f_{j}(y')$ \;
  				}
  			}
  			\ForEach{ $j \in B(i)$}{
  				
  				\If{ ($ g^{te}(y'| \lambda^j, z) < g^{te}(x^j| \lambda^j, z))$ }
  				{
  					\textbf{set} $ x^j \longleftarrow y'$ \;
  					\textbf{set} $ FV^j \longleftarrow F(y')$ \;
  				}
  			}
  			\textbf{Remove} from \textbf{EP} all the vectors dominated by $F(y')$ \;
  			\textbf{Add} $F(y')$ to \textbf{EP} if no vectors in \textbf{EP} dominate $F(y')$ \;
  		}	
  	}	
  	\caption{ \ac{MOEA/D} Framework}
  	 	\label{MOEA/D}
  \end{algorithm}

  \subsubsection{Pareto based: NSGA-II}
  \ac{NSGA-II} \cite{deb2002fast}, is an elitist evolutionary multi-objective algorithm  that extends the \ac{GA}  to deal with multi-objective problems. \ac{NSGA-II} is based on the idea of the non dominance to rank solutions into different fronts, and crowding distances to select the most spreedest solutions from the same front. The  \ac{NSGA-II} framework is described by Algorithm \ref{NSGA-II}.
  
  \begin{algorithm}[!htb]
  	\SetAlgoLined
  
  	\textbf{Input:}\\
  	a stopping criterion \;
  	population size: N \;
  	
  	\textbf{ Output: } population $P_t $ \;	
  	
  	\textbf{Initialization:} \;
  	\textbf{Set} \textbf{$Q_0$}=$\emptyset$ \;	
  	\textbf{Generate} an initial population $P_0$ of size $N$\;
  	\textbf{Ranking} $P_0$\;
  	\textbf{Crowding} $P_0$\;
  	\textbf{Set}  $ t \longleftarrow 0$\;
  	
  	\While{(the stopping condition is not reached)}{
  		\For{$i\leftarrow 1$ \KwTo $N$}  {
  			
  			\textbf{Set} $k  \longleftarrow$ \textbf{TournamentSelection}($P_t$, 2) \;
  			\textbf{Set} $l  \longleftarrow$ \textbf{TournamentSelection}($P_t$, 2) \;
  			
  			\textbf{Set} 	$ y_1 \longleftarrow$ \textbf{Crossover}($x^k$,$x^l$) \;

  			\textbf{Mutation} ($y_1$) \;
  			\textbf{Repair} ($y_1$) \;
  			
  			\textbf{Set} $ Q_t \longleftarrow  Q_t \cup \{y_1 \}$   \;	
  		}	
  		\textbf{Set} $ R_t \longleftarrow  P_t \cup Q_t$   \;	
  		\textbf{Ranking} $R_t = F_1 \cup F_2,....,\cup F_{Last}$\;

  		\textbf{Set} $ i \longleftarrow  1$\;
  		\textbf{Set}  $ sum \longleftarrow  0$\;
  		\While{( $ sum < N$ )}
  		{
  			\textbf{Crowding} $F_i$\;
  			\textbf{Set} 	$ sum \longleftarrow  sum + \mid F_i \mid $\;	
  			\textbf{Set} 	$ i \longleftarrow  i+ 1$\;	
  		}
  		
  		\While{( $ \mid P_{t+1}  \mid \cup \mid F_i \mid < N$ )}
  		{
  			\textbf{Set} $ P_{t+1} \longleftarrow  P_{t+1} \cup F_i$   \;
  			\textbf{Set} $ i \longleftarrow  i+1  $\;	
  		}
  		\textbf{Select} best ($N - \mid  P_{t+1} \mid$) solutions from $F_i$ in terms of Crowding distance\;
  		\textbf{Set}  $ t \longleftarrow  t+1  $\;	
  	}	
  	\caption{ \ac{NSGA-II} Framework}
  		\label{NSGA-II}
  \end{algorithm}

  \subsubsection{ Hyprid NSGA-II with MOEA/D}
  The following algorithm is the hybridization of two  algorithms \ac{NSGA-II} (Algorithm \ref{NSGA-II}) and \ac{MOEA/D} (Algorithm \ref{MOEA/D}), which are successively performed  to improve the efficiency of the overall search. According to \cite{trivedi2016survey}, efficiently combining dominance and decomposition based approaches can lead to better performance during the resolution phase. In addition, the hybrid algorithm can deal with multiple objectives (more than three) and avoid issues known in \ac{NSGA-II} since obtained solutions by \ac{NSGA-II} are further enhanced using MOEA/D. First, an initial population is generated (\textbf{line: 6}) for \ac{NSGA-II} to return non-dominated solutions, which are then stored in $N_D$ (\textbf{line: 7}). Another initial population $ P_1 $ of size $ N $ will be constructed from non-dominated solutions $N_D$ and possibly other randomly generated solutions. We first set $ P_1 $ to the set of non-dominated solutions of \ac{NSGA-II} (\textbf{line: 8}) and we then complete it, if the size N is not yet reached, by  $(N-|Nd|)$ solutions randomly generated (\textbf{lines: 9-12}) as described in subsection \ref{IP} by setting $P_{size}=N-|Nd|$. Afterwards, \ac{MOEA/D} algorithm is performed to enhance this new initial population, and store solutions in the archive $EP$.
  
  
  \begin{algorithm}
  	\SetAlgoLined
  	
  	\textbf{Input:} \\
  	NSGA-$II_{input}$  \;
  	MOEA/$D_{input}$\;
  	$N_D$ : non-dominated solutions returned by NSGA-II\;
  	\textbf{Output: } EP \;	
  	
  	\textbf{Generate} an initial population $P_0$ of size $N$\;
  	\textbf{Set} $N_D \longleftarrow$ \textbf{NSGA-II($P_0,N$)}\;

  	\textbf{Set}  $ P_1 \longleftarrow   N_D $\;
  	\While{( $ \mid P_1 \mid \leq N$ )}
  	{
  		\textbf{Generate} randomly a solution $X$ \;	
  		\textbf{Set} $ P_1 \longleftarrow  P_1  \cup \{X\}$\;
  	}
  	\textbf{Set}  $ EP \longleftarrow $  \textbf{MOEA/D($P_1, \lambda, N, T$)} \;

  	\caption{ Hybrid \ac{NSGA-II} with MOEA/D}
  	
  	\label{Hybrid}
  \end{algorithm}

  \subsection{Numerical experiments }\label{MultiNE}
  
  \subsubsection{Test instances}

 The test instances  described  and generated in \ref{NESMTW} are used to test the algorithms.

  \subsubsection{Computational results}

  \begin{table}[!h]
  	\caption{ Tuning parameters}
  	\centering
  	\def\arraystretch{1.3}
  	
  	\begin{tabular}{ lll}
  		
  		\hline
  		
  		\textbf{Algorithm} &\textbf{ Parameter }& \textbf{Value}  \\   
  		\hline

  		& Mutation probability  &   0.08 \\  
  		\textbf{NSGA-II} & Crossover probability  &   0.8 \\  
  		& Population size   & $ 10 \times n$ \\  
  		
  		\hline
  		
  		& Mutation probability  &   0.08 \\  
  		\textbf{MOEA/D} & Crossover probability  &   1 \\  
  		& Population size   & $ 10 \times n$ \\  
  		& T  &   n \\  
  		& Archive size   & $ 10 \times n$ \\  
  		\hline
  		
  		\textbf{Hybrid } &\multicolumn{2}{l}{The same parameters of \ac{NSGA-II} and \ac{MOEA/D} are used }    \\     
  		\hline

  	\end{tabular}

  	\label{table:TM}

  \end{table}
  

  
    Instances are solved with a maximum number iterations of $5\times10^4$. Tuning parameters are set as described as in Table \ref{table:TM}. At each iteration, the number of non-dominated solutions maintained by \ac{NSGA-II} is always equal to the population size $N$. Indeed, Parent and offspring populations are combined in a mating pool $R_t$ then the best $N$ solutions are selected using Pareto dominance and the crowding distance. MOEA /D could use an unlimited archive to store non-dominated solutions, but here we have set its size to that of the population. \ac{NSGA-II} and \ac{MOEA/D} algorithms improve the initial population and try to move it towards the Pareto optimal solutions. The hybrid algorithm first finds non-dominated solutions by the mechanism of \ac{NSGA-II} and then attempts to further improve them using the procedure of MOEA/D.
  
  In order to perform the comparison between the three algorithms, we run them with a single simulation, to solve each instance by considering single then double time windows and using the same initial population. Indicators performance used to compare performance of these algorithms are:
  \begin{itemize}
  	\item \textit{hyper-volume indicator \cite{zitzler1999multiobjective}}: it measures  the hyper-volume covered by a set of solutions in the objective space. This metric has two issues. On the one hand, it needs a reference point and, on the other hand, the complexity increases monotonically with the number of objective functions. The algorithm with higher value performs better;
  	
  	\item\textit{ coverage of two sets \cite{zitzler1999multiobjective}:} this indicator computes the proportion of a set of solutions $B$ dominated by another set  of solutions $A$. The value 0 means that solutions in the set $B$ are not dominated by any solutions from the set $A$, the value 1 indicates that each solution in $B$ is dominated  at least by one solution of $A$. it is computed as follows:
  	\begin{equation}
  	Cov(A,B) = \frac{\mid \{ \; b \; \in \;B:\; \exists \; a \; \in \; A, \; a \;\prec \;b \} \mid}{\mid B \mid}
  	\end{equation}

  	\item \textit{Pareto solution size}: this indicator measures the cardinality of the Pareto front and represents the number of trade-off solutions from whose the decision maker will select a final preferred solution;
  	
  	\item\textit{ CPU time}: the most known indicator for comparing the complexity of algorithms, best algorithms have shortest CPU times. 
  \end{itemize}

  Before we compute the hyper-volume indicator, we consider (1,1,1) as a reference point and we normalize objective function values using equation (\ref{eqn:fnorm}). $f_i^{min}$ and $f_i^{max}$ are respectively the minimum and the maximum objective function values found by the three algorithms considering each instance with a single and double availability periods.
  
  \begin{equation} \label{eqn:fnorm}
  f_i^{normalized} = \frac{f_i - f_i^{min}}{f_i^{max} - f_i^{min}}
  \end{equation}  
  
  \begin{table}[!htb]
  	\caption{Comparison between instances with one and two availability periods according to coverage set indicator}
  	\centering
  	\def\arraystretch{1.3}
  	
  	\begin{tabular}{l|ll|ll|ll  }
  		\hline
  		
  		\multirow{3}{*}{\textbf{Instances}} & \multicolumn{2}{l|}{\textbf{NSGA-II}}	 & \multicolumn{2}{c|}{\textbf{MOEA/D}}	& \multicolumn{2}{c}{ \textbf{Hybrid} } \\ 
  		
  		\cline{2-3} \cline{4-5} \cline{6-7}

  		& 1 period& 2 periods& 1 period& 2 periods& 1 period& 2 periods \\	
  		& 2 periods& 1 period& 2 periods& 1 period& 2 periods& 1 period\\
  		
  		\hline
  		
  		A1	&	0	&	\textbf{1}	&	0	&	\textbf{0.947}	&	0	&	\textbf{1}	\\
  		A2	&	0	&	\textbf{0.929}	&	0	&	\textbf{0.929}	&	0	&	\textbf{0.929}	\\
  		A3	&	0.093	&	\textbf{0.533}	&	0.088	&	\textbf{0.533}	&	0.088	&	\textbf{0.533}	\\
  		A4	&	0	&	\textbf{1}	&	0	&	\textbf{1}	&	0	&	\textbf{1}	\\
  		A5	&	0	&	\textbf{1}	&	0	&	\textbf{1}	&	0	&	\textbf{1}	\\
  		A6	&	0	&	\textbf{1}	&	0	&	\textbf{1}	&	0	&	\textbf{1}	\\
  		A7	&	0	&	\textbf{1}	&	0	&	\textbf{1}	&	0	&	\textbf{1}	\\
  		
  		\hline
  		
  		B1	&	0	&	\textbf{1}	&	0	&	\textbf{1}	&	0	&	\textbf{1}	\\
  		B2	&	0	&	\textbf{1}	&	0	&	\textbf{1}	&	0	&	\textbf{1}	\\
  		B3	&	0	&	\textbf{1}	&	0	&	\textbf{1}	&	0	&	\textbf{1}	\\
  		B4	&	0	&	\textbf{1}	&	0	&	\textbf{1}	&	0	&	\textbf{1}	\\
  		B5	&	0	&	\textbf{1}	&	0	&	\textbf{1}	&	0	&	\textbf{1}	\\
  		B6	&	0	&	\textbf{1}	&	0	&	\textbf{1}	&	0	&	\textbf{1}	\\
  		B7	&	0	&	\textbf{1}	&	0	&	\textbf{1}	&	0	&	\textbf{1}	\\
  		
  		\hline
  		
  		C1	&	0	&	\textbf{1}	&	0	&	\textbf{1}	&	0	&	\textbf{1}	\\
  		C2	&	0	&	\textbf{1}	&	0	&	\textbf{1}	&	0	&	\textbf{1}	\\
  		C3	&	0	&	\textbf{1}	&	0	&	\textbf{1}	&	0	&	\textbf{1}	\\
  		C4	&	0	&	\textbf{0.958}	&	0	&	\textbf{0.808}	&	0	&	\textbf{0.960}	\\
  		C5	&	0	&	\textbf{1}	&	0.045	&	\textbf{0.244}	&	0	&	\textbf{1}	\\
  		C6	&	0	&	\textbf{1}	&	0	&	\textbf{1}	&	0	&	\textbf{1}	\\
  		C7	&	0	&	\textbf{1}	&	0	&	\textbf{1}	&	0	&	\textbf{1}	\\

  		\hline

  	\end{tabular}

  	\label{table:CoverSetAccTW}

  \end{table}

  \begin{table}[p]
  	\caption{Comparison of algorithms  according  to hypervolume indicator, front size and CPU time.}
  	\centering
  	\def\arraystretch{1.0}
  	\hspace*{-2.2cm}\setlength\tabcolsep{1.5pt}
  	\begin{tabular}{ l|lll|lll|lll  }
  		
  		\hline
  		
  		\multirow{2}{*}{\textbf{}}& \multicolumn{3}{c|}{\textbf{Hypervolume}}  & \multicolumn{3}{c|}{\textbf{Front size}} & \multicolumn{3}{c}{\textbf{CPU time}}\\	
  		
  		\cline{2-4} \cline{5-7} \cline{8-10}
  		
  		\rot{\textbf{Instances}}& \rot{\textbf{MOEA/D}} &\rot{\textbf{NSGA-II}}  &\rot{\textbf{Hybrid}} &\rot{\textbf{MOEA/D}} &\rot{\textbf{NSGA-II}}  &\rot{\textbf{Hybrid}}& \rot{\textbf{MOEA/D}} &\rot{\textbf{NSGA-II}}  &\rot{\textbf{Hybrid}} \\
  		
  		\hline

  		A1-1	&	0.247077	&	0.247077	&	0.247077	&	19	&	19	&	19	&	\textbf{11.47}	&	66.62	&	78.32	\\
  		A2-1	&	0.37748	    &	0.37748		&	0.37748		&	14	&	14	&	14	&	\textbf{11.54}	&	65.30	&	76.45	\\
  		A3-1	&	0.244444	&	0.244444	&	0.244444	&	15	&	15	&	15	&	\textbf{11.60}	&	66.25	&	77.44	\\
  		A4-1	&	0.027405	&	0.027405	&	0.027405	&	12	&	12	&	12	&	\textbf{11.69}	&	70.49	&	81.27	\\
  		A5-1	&	0.097665	&	0.097665	&	0.097665	&	20	&	20	&	20	&	\textbf{11.21}	&	61.63	&	72.57	\\
  		A6-1	&	0.121518	&	0.121518	&	0.121518	&	\textbf{23}	&	22	&	22	&	\textbf{11.73}	&	59.02	&	71.21	\\
  		A7-1	&	0.22591		&	0.22591	&	0.22591	&	18	&	18	&	18	&	\textbf{11.41}	&	70.23	&	82.12	\\
  		\hline
  		A1-2	&	0.818409	&	0.818243	&	\textbf{0.820858}	&	70	&	63	&	\textbf{78}	&	\textbf{14.97}	&	57.70	&	74.69	\\
  		A2-2	&	0.864433	&	\textbf{0.864814}	&	\textbf{0.864814}	&	\textbf{22}	&	16	&	18	&	\textbf{12.33}	&	58.40	&	70.28	\\
  		A3-2	&	\textbf{0.459722}	&	0.459718	&	\textbf{0.459722}	&	\textbf{57}	&	54	&	\textbf{57}	&	\textbf{15.07}	&	58.74	&	72.86	\\
  		A4-2	&	0.669616	&	\textbf{0.695665}	&	\textbf{0.695665}	&	18	&	\textbf{33}	&	\textbf{33}	&	\textbf{12.69}	&	56.43	&	68.81	\\
  		A5-2	&	0.771339	&	\textbf{0.842371}	&	\textbf{0.842371}	&	12	&	\textbf{15}	&\textbf{15}	&	\textbf{11.10}	&	58.24	&	69.63	\\
  		A6-2	&	0.712418	&	0.734699	&	\textbf{0.735039}	&	19	&	30	&	\textbf{31}	&	\textbf{11.63}	&	56.78	&	69.27	\\
  		A7-2	&	0.683382	&	\textbf{0.688122}	&	\textbf{0.688122}	&	12	&	12	&	\textbf{16}	&	\textbf{11.39}	&	61.92	&	72.26	\\
  		
  		\hline
  		
  		B1-1	&	0.248781	&	0.277003	&	\textbf{0.278062}	&	\textbf{250}	&	217	&	\textbf{250}	&	\textbf{88.78}	&	260.42	&	380.87	\\
  		B2-1	&	0.385349	&	0.42151	&	\textbf{0.423894}	&	\textbf{250}	&	212	&	\textbf{250}	&	\textbf{85.79}	&	251.51	&	373.38	\\
  		B3-1	&	0.285949	&	0.325361	&	\textbf{0.326721}	&	182	&	216	&	\textbf{250}	&	\textbf{89.88}	&	261.69	&	384.91	\\
  		B4-1	&	0.316459	&	0.387703	&	\textbf{0.390148}	&	141	&	213	&	\textbf{250}	&\textbf{89.61}	&	260.92	&	379.11	\\
  		B5-1	&	0.398129	&	0.419883	&	\textbf{0.427596}	&	\textbf{250}	&	224	&	\textbf{250}	&	\textbf{90.80}	&	255.05	&	378.25	\\
  		B6-1	&	0.392038	&	0.414666	&	\textbf{0.41899}	&	\textbf{250}	&	219	&	\textbf{250}	&	\textbf{104.76}	&	245.72	&	365.50	\\
  		B7-1	&	0.209099	&	0.237091	&	\textbf{0.240017}	&	173	&	215	&	\textbf{250}	&\textbf{97.36}	&	269.63	&	386.60	\\
  		\hline
  		B1-2	&	0.764303	&	0.879726	&	\textbf{0.884279}	&	201	&	74	&	\textbf{250}	&	\textbf{87.18}	&	267.97	&	388.79	\\
  		B2-2	&	0.77076	&	0.868432	&	\textbf{0.875051}	&	94	&	193	&	\textbf{250}	&	\textbf{82.40}	&	250.35	&	358.56	\\
  		B3-2	&	0.755139	&	0.961991	&	\textbf{0.967665}	&	72	&	132	&	\textbf{145}	&	\textbf{82.11}	&	246.97	&	343.27	\\
  		B4-2	&	0.78797	&	0.934304	&	0.\textbf{937031}	&	72	&	209	&	\textbf{250}	&\textbf{87.93}	&	248.31	&	373.53	\\
  		B5-2	&	0.881285	&	0.945488	&	\textbf{0.948753}	&	82	&	145	&	\textbf{152}	&\textbf{86.7}8	&	265.72	&	363.97	\\
  		B6-2	&	0.763252	&	0.905613	&	\textbf{0.908228}	&	70	&	175	&	\textbf{190}	&	\textbf{86.55}	&	261.62	&	368.07	\\
  		B7-2	&	0.637109	&	0.837473	&	\textbf{0.846069}	&	50	&	171	&	\textbf{249}	&	\textbf{79.24}	&	251.93	&	366.44	\\
  		
  		\hline
  		C1-1	&	0.307671	&	0.504296	&	\textbf{0.514147}	&	\textbf{500}	&	456	&	\textbf{500}	&\textbf{471.16}	&	1012.95	&	1561.80	\\
  		C2-1	&	0.273967	&	0.386123	&	\textbf{0.399}	&	\textbf{500}	&	473	&	\textbf{500}	&\textbf{472.01}	&	982.63	&	1506.07	\\
  		C3-1	&	0.278814	&	0.566633	&	\textbf{0.572585}	&	64	&	446	&	\textbf{500}	&\textbf{364.22}	&	934.72	&	1469.67	\\
  		C4-1	&	0.279834	&	0.399067	&	\textbf{0.432489}	&	\textbf{500}	&	471	&	\textbf{500}	&\textbf{422.63}	&	992.47	&	1539.73	\\
  		C5-1	&	0.392415	&	0.564469	&	\textbf{0.564884}	&	\textbf{500}	&	468	&	\textbf{500}	&\textbf{443.17}	&	971.72	&	1538.01	\\
  		C6-1	&	0.319871	&	0.505397	&	\textbf{0.506527}	&	382	&	464	&	\textbf{500}	&	\textbf{429.46}	&	963.15	&	1543.40	\\
  		C7-1	&	0.230476	&	0.494416	&	\textbf{0.500799}	&	71	&	467	&	\textbf{500}	&	\textbf{421.11}	&	952.47	&	1515.98	\\
  		\hline
  		C1-2	&	0.630014	&	0.925657	&	\textbf{0.927298}	&	82	&	410	&	\textbf{500}	&	\textbf{387.89}	&	963.36	&	1467.46	\\
  		C2-2	&	0.636132	&	0.952033	&	\textbf{0.957938}	&	103	&	410	&	\textbf{500}	&	\textbf{398.10}	&	951.38	&	1519.20	\\
  		C3-2	&	0.47998	&	0.967442	&	\textbf{0.969775}	&	40	&	354	&	\textbf{500}	&	\textbf{376.67}	&	929.88	&	1419.33	\\
  		C4-2	&	0.577446	&	0.895814	&	\textbf{0.917401}	&	73	&	412	&	\textbf{500}	&	\textbf{385.00}	&	990.48	&	1492.62	\\
  		C5-2	&	0.374383	&	0.946422	&	\textbf{0.950055}	&	89	&	417	&	\textbf{500}	&\textbf{396.70}	&	1031.74	&	1592.74	\\
  		C6-2	&	0.555871	&	0.958149	&	\textbf{0.960983}	&	92	&	366	&	\textbf{500}	&	\textbf{383.10}	&	950.24	&	1477.62	\\
  		C7-2	&	0.530525	&	0.925249	&	\textbf{0.925615}	&	76	&	412	&	\textbf{500	}&	\textbf{387.90}	&	986.33	&	1536.91	\\

  		\hline

  	\end{tabular}

  	\label{table:HvFsCPU}

  \end{table}


  \begin{table}[p]
  	\caption{ Comparison of algorithms according  to coverage set indicator.}
  	\centering
  	\def\arraystretch{1.0}
  	
  	\begin{tabular}{ l|ll|ll|ll   }
  		\hline

  		\multirow{2}{*}{\textbf{Instances}}& \textbf{NSGA-II} & \textbf{MOEA/D}  & \textbf{Hybrid}  & \textbf{NSGA-II} & \textbf{Hybrid} & \textbf{MOEA/D}  \\
  		
  		\cline{2-2} \cline{3-3} \cline{4-4}	 \cline{5-5} \cline{6-6} \cline{7-7}
  		
  		& \textbf{MOEA/D}  & \textbf{NSGA-II} & \textbf{NSGA-II} & \textbf{Hybrid}  & \textbf{MOEA/D} & \textbf{Hybrid} \\
  		
  		\hline
  		
  		A1-1	&	0	&	0	&	0	&	0	&	0	&	0	\\
  		A2-1	&	0	&	0	&	0	&	0	&	0	&	0	\\
  		A3-1	&	0	&	0	&	0	&	0	&	0	&	0	\\
  		A4-1	&	0	&	0	&	0	&	0	&	0	&	0	\\
  		A5-1	&	0	&	0	&	0	&	0	&	0	&	0	\\
  		A6-1	&\textbf{0.043}	&	0	&	0	&	0	&\textbf{0.043}	&	0	\\
  		A7-1	&	0	&	0	&	0	&	0	&	0	&	0	\\
  		\hline
  		A1-2	&	0.114	&	\textbf{0.127}	&\textbf{0.127}	&	0	&\textbf{0.114}	&	0	\\
  		A2-2	&	\textbf{0.227}	&	0.063	&	\textbf{0.0625}	&	0	&\textbf{0.227}	&	0	\\
  		A3-2	&	0	&	0	&	0	&	0	&	0	&	0	\\
  		A4-2	&	\textbf{0.333}	&	0	&	\textbf{0.333}	&	0	&	0	&	0	\\
  		A5-2	&	\textbf{0.417}	&	0	&	0	&	0	&	\textbf{0.417}	&	0	\\
  		A6-2	&	0	&	0	&	0	&	0	&	0	&	0	\\
  		A7-2	&	0	&	\textbf{0.083}	&	\textbf{0.083}	&	0	&	0	&	0	\\
  		
  		\hline
  		
  		B1-1	&	\textbf{0.408}	&	0.097	&	\textbf{0.101}	&	0	&	\textbf{0.508}	&	0.200	\\
  		B2-1	&	\textbf{0.524}	&	0.024	&	\textbf{0.094}	&	0	&	\textbf{0.704}	&	0	\\
  		B3-1	&	\textbf{0.462}	&	0.106	&	\textbf{0.222}	&	0	&	\textbf{0.489}	&	0.040	\\
  		B4-1	&	\textbf{0.702}	&	0.056	&	\textbf{0.235}	&	0	&	\textbf{0.844}	&	0.008	\\
  		B5-1	&	0.060	&	\textbf{0.205}	&	\textbf{0.205}	&	0	&	0.092	&	\textbf{0.184}	\\
  		B6-1	&	0.140	&	\textbf{0.201}	&	\textbf{0.269}	&	0	&	\textbf{0.248}	&	0.032	\\
  		B7-1	&	0.035	&	\textbf{0.098}	&	\textbf{0.060}	&	0	&	0.040	&	\textbf{0.084}	\\
  		\hline
  		B1-2	&	\textbf{0.568}	&	0.015	&	\textbf{0.065}	&	0	&	\textbf{0.568}	&	0.016	\\
  		B2-2	&	\textbf{0.543}	&	0.057	&	\textbf{0.161}	&	0	&	\textbf{0.543}	&	0.048	\\
  		B3-2	&	\textbf{1}	&	0	&	\textbf{0.091}	&	0	&	\textbf{1}	&	0	\\
  		B4-2	&	\textbf{0.583}	&	0.029	&	\textbf{0.249}	&	0	&	\textbf{0.597}	&	0.016	\\
  		B5-2	&	\textbf{0.866}	&	0.021	&	\textbf{0.014}	&	0	&	\textbf{0.866}	&	0.020	\\
  		B6-2	&	\textbf{1}	&	0	&	\textbf{0.126}	&	0	&	\textbf{1}	&	0	\\
  		B7-2	&	\textbf{1}	&	0	&	\textbf{0.211}	&	0	&	\textbf{1}	&	0	\\
  		
  		\hline
  		
  		C1-1	&	\textbf{0.152}	&	0.037	&	\textbf{0.061}	&	0	&	\textbf{0.150}	&	0.036	\\
  		C2-1	&	\textbf{0.144}	&	0.061	&	\textbf{0.106}	&	0	&	\textbf{0.144}	&	0.058	\\
  		C3-1	&	\textbf{0.484}	&	0    	&	\textbf{0.085}	&	0	&	\textbf{0.484}	&	0	\\
  		C4-1	&	\textbf{0.258}	&	0.053	&	\textbf{0.244}	&	0	&	\textbf{0.258}	&	0.050	\\
  		C5-1	&	0.028	&	\textbf{0.060}	&	\textbf{0.006}	&	0	&	0.028	&	\textbf{0.056}	\\
  		C6-1	&	0.058	&	\textbf{0.078}	&	\textbf{0.013}	&	0	&	0.058	&	\textbf{0.072}	\\
  		C7-1	&	\textbf{0.592}	&	0.002	&	\textbf{0.105}	&	0	&	\textbf{0.592}	&	0.002	\\
  		\hline
  		C1-2	&	\textbf{1}	&	0	&	\textbf{0.078}	&	0	&	\textbf{1}	&	0	\\
  		C2-2	&	\textbf{1}	&	0	&	\textbf{0.320}	&	0	&	\textbf{1}	&	0	\\
  		C3-2	&	\textbf{1}	&	0	&	\textbf{0.232}	&	0	&	\textbf{1}	&	0	\\
  		C4-2	&	\textbf{1}	&	0	&	\textbf{0.330}	&	0	&	\textbf{1}	&	0	\\
  		C5-2	&	\textbf{1}	&	0	&	\textbf{0.091}	&	0	&	\textbf{1}	&	0	\\
  		C6-2	&	\textbf{1}	&	0	&	\textbf{0.194}	&	0	&	\textbf{1}	&	0	\\
  		C7-2	&	\textbf{1}	&	0	&	\textbf{0.056}	&	0	&	\textbf{1}	&	0	\\

  		\hline
  		
  	\end{tabular}

  	\label{table:AlgorithmsCoverageSet}

  \end{table}

  \begin{figure}[p]
  	\begin{center}
  		\begin{tikzpicture}
  		\begin{axis}[
  		width=\linewidth, 
  		height=7.5cm,
  		grid=major, 
  		grid style={dashed,gray!30},
  			xlabel= Instances, 
  		ylabel= Hypervolume,
  		xtick=data,
  		legend style={at={(0.4,0.11)},anchor=west},
  		,xticklabels={A1,A2,A3,A4,A5,A6,A7,B1,B2,B3,B4,B5,B6,B7,C1,C2,C3,C4,C5,C6,C7}
  		]
  		
  		\addplot+[] 
  		table[x=I, y=M1, col sep=comma] {src/HVTW1.csv}; 
  		
  		\addplot+[] 
  		table[x=I, y=M2, col sep=comma] {src/HVTW1.csv};

  		\legend{MOEA/D: 1 TW ,MOEA/D: 2 TW}
  		\end{axis}
  		
  		\end{tikzpicture}
  		\vspace{-6mm}
  		\caption{Hypervolume indicator for \ac{MOEA/D} algorithm according to time windows (TW)}
  		\label{MOEA/D/TW}

  		\begin{tikzpicture}
  		\begin{axis}[
  		width=\textwidth, 
  		height=7.5cm,
  		grid=major, 
  		grid style={dashed,gray!30},
  			xlabel= Instances, 
  		ylabel= Hypervolume,
  		xtick=data,
  		legend style={at={(0.4,0.11)},anchor=west},
  		,xticklabels={A1,A2,A3,A4,A5,A6,A7,B1,B2,B3,B4,B5,B6,B7,C1,C2,C3,C4,C5,C6,C7}
  		]
  		
  		\addplot+[] 
  		table[x=I, y=N1, col sep=comma] {src/HVTW1.csv}; 
  		
  		\addplot+[] 
  		table[x=I, y=N2, col sep=comma] {src/HVTW1.csv};

  		\legend{NSGA-II: 1 TW ,NSGA-II: 2 TW}
  		\end{axis}
  		
  		\end{tikzpicture}
  		\vspace{-6mm}
  		\caption{Hypervolume indicator for \ac{NSGA-II} algorithm according to time windows (TW)}
  		\label{NSGA-II/TW}

  		\begin{tikzpicture}
  		\begin{axis}[
  		width=\textwidth, 
  		height=7.5cm,
  		grid=major, 
  		grid style={dashed,gray!30},
  			xlabel= Instances, 
  		ylabel= Hypervolume,
  		xtick=data,
  		legend style={at={(0.4,0.11)},anchor=west},
  		,xticklabels={A1,A2,A3,A4,A5,A6,A7,B1,B2,B3,B4,B5,B6,B7,C1,C2,C3,C4,C5,C6,C7}
  		]
  		
  		\addplot+[] 
  		table[x=I, y=H1, col sep=comma] {src/HVTW1.csv}; 
  		
  		\addplot+[] 
  		table[x=I, y=H2, col sep=comma] {src/HVTW1.csv};

  		\legend{Hybrid: 1 TW ,Hybrid: 2 TW}
  		\end{axis}
  		
  		\end{tikzpicture}
  		\vspace{-6mm}
  		\caption{Hypervolume indicator for Hybrid \ac{NSGA-II} with \ac{MOEA/D} algorithm according to time windows (TW)}
  		\label{Hybrid/TW}
  	\end{center}
  \end{figure}
  

  \begin{figure}[p]
  	\begin{center}
  		\begin{tikzpicture}
  		\begin{axis}[
  		width=\textwidth, 
  		height=9cm,
  		grid=major, 
  		grid style={dashed,gray!30},
  		xlabel= Instances, 
  		ylabel= Hypervolume,
  		xtick=data,
  		legend style={at={(0.4,0.11)},anchor=west},
  		,xticklabels={A1,A2,A3,A4,A5,A6,A7,B1,B2,B3,B4,B5,B6,B7,C1,C2,C3,C4,C5,C6,C7}
  		]
  		
  		\addplot+[] 
  		table[x=I, y=N1, col sep=comma] {src/HVTW1.csv}; 
  		
  		\addplot+[] 
  		table[x=I, y=H1, col sep=comma] {src/HVTW1.csv}; 
  		
  		\addplot+[] 
  		table[x=I, y=M1, col sep=comma] {src/HVTW1.csv};

  		\legend{NSGA-II,Hybrid,MOEA/D}
  		\end{axis}
  		
  		\end{tikzpicture}
  		\caption{ Comparison of algorithms for instances with one availability period  according  to hypervolume indicator }
  		\label{AlgTw1}
  		
  		\vspace{10mm}
  		
  		\begin{tikzpicture}
  		\begin{axis}[
  		width=\textwidth, 
  		height=9cm,
  		grid=major, 
  		grid style={dashed,gray!30},
  		xlabel= Instances, 
  		ylabel= Hypervolume,
  		xtick=data,
  		legend style={at={(0.4,0.11)},anchor=west},
  		,xticklabels={A1,A2,A3,A4,A5,A6,A7,B1,B2,B3,B4,B5,B6,B7,C1,C2,C3,C4,C5,C6,C7}
  		]
  		
  		\addplot+[] 
  		table[x=I, y=N2, col sep=comma] {src/HVTW1.csv}; 
  		
  		\addplot+[] 
  		table[x=I, y=H2, col sep=comma] {src/HVTW1.csv}; 
  		
  		\addplot+[] 
  		table[x=I, y=M2, col sep=comma] {src/HVTW1.csv};

  		\legend{NSGA-II,Hybrid,MOEA/D}
  		\end{axis}
  		
  		\end{tikzpicture}
  		\caption{ Comparison of algorithms for instances with two availability periods  according  to hypervolume indicator}
  		\label{AlgTw2}
  	\end{center}
  \end{figure}

  According to CPU time indicator, \ac{MOEA/D} solved instances faster. In the case where the decision maker is not interested to obtain the whole Pareto front, \ac{MOEA/D} may be the best choice to solve the \ac{HHCRSP} (see Table \ref{table:HvFsCPU}). According to the hypervolume indicator,  \ac{NSGA-II} and the hybrid algorithms outperformed  \ac{MOEA/D} algorithm. Taking into consideration both hypervolume and Pareto front size indicators showed that  the hybrid algorithm  outperformed  \ac{MOEA/D} and \ac{NSGA-II} algorithms.
  
  
  The comparison according to the coverage metric is also carried out (see Table \ref{table:AlgorithmsCoverageSet}). This indicator did not differentiate algorithms performances for small instances A (single or double time windows) and B (single time window) because the algorithms generated very few non-dominated solutions compared to medium and large instances. However, For instances $B$ (double time windows) and $C$ (single or double time windows), \ac{NSGA-II}  performs better than MOEA/D. The percentage of solutions found by \ac{MOEA/D} and dominated by solutions of \ac{NSGA-II} is remarkable. The hybrid algorithm outperformed both \ac{NSGA-II} and \ac{MOEA/D} algorithms. 
  In most of the 42 instances,  except instances $B5$-1, $B7$-1, $C5$-1, and $C6$-1, the percentage of \ac{MOEA/D} solutions dominated by those of the hybrid algorithm is greater and significant compared to the proportion of solutions of the hybrid algorithm that are dominated by those of MOEA/D.
  
  For all instances, \textit{Cov(NSGA-II, Hybrid)=0} because the set of solutions found by \ac{NSGA-II} are improved by MOEA/D. In the worst case, \ac{MOEA/D} will not be able to improve solutions found by \ac{NSGA-II} , which is the case for instances of set A with single time window  and instances A3-2, A5-2 and A6-2. 
  Accordingly, \textit{Cov(Hybrid, NSGA-II)$>0$  } for all instances, except for previously mentioned  instances, which indicate \ac{MOEA/D} could  further improve the quality of solutions found by NSGA-II. Indeed, \textit{Cov(Hybrid, NSGA-II)$>0$ } means that hybrid could  reach better solutions, which dominate some solutions found by NSGA-II. The higher  the value of \textit{Cov(Hybrid, NSGA-II) }, the higher the proportion of dominated solutions reached by NSGA-II.
  
  In order to study the impact of the number of availability periods on the quality of solutions, the coverage set metric is computed to compare the approximated Pareto front found for each instance considering single and double time windows. Table \ref{table:CoverSetAccTW} clearly shows the advantage of instances with double availability periods. A significant improvement is performed to the approximated Pareto front when double time windows for patient are considered. For all instances, except $A1, A2, A3, C4$ and $C4$, each solution from the approximated Pareto front for the instances with single time window is dominated at least by one solution from the approximated Pareto front for the instances with double time windows\textit{ (Cov(2 periods, 1 period)=1)}. Accordingly, solutions from the approximated Pareto front for instances with double time windows, except for previously mentioned  instances, are not dominated by any solution from the set of instances with single time window ($Cov(1 \;period,\; 2 \;periods)=0$).
  
  To sum up, the \ac{MOEA/D} solves instances faster and is suitable if the decision maker is not interested in having all Pareto optimal solutions. The \ac{NSGA-II} algorithm is very efficient and could found solutions close to the hybrid algorithm in terms of quality considering hyper-volume indicator (see figures \ref{AlgTw1} and \ref{AlgTw2}). However,  \ac{NSGA-II}  is not suitable for problems with more than three objectives. The hybrid algorithm  outperformed  \ac{NSGA-II} and \ac{MOEA/D} algorithms, which is justified by performance indicators  used (hyper-volume, coverage set and Pareto front size: see tables \ref{table:HvFsCPU} and \ref{table:AlgorithmsCoverageSet}). Using multiple time windows improved solutions quality (see figures \ref{MOEA/D/TW}, \ref{NSGA-II/TW} and \ref{Hybrid/TW}) and these solutions dominate solutions of instances with single  availability period (see Table \ref{table:CoverSetAccTW}).

  \section{Conclusion}\label{MultiC}

  Home health care companies often deal with conflicting  objectives, such as minimizing the cost while increasing patients' satisfactions. Most studies transform multi-objective problems  by  into a mono-objective case using scalarazition techniques. The issue is how to  assign weights to objective functions, which may be a confusing task for the decision maker. Consequently, he needs to be involved  a posteriori to select a final preferred solution from  the non-dominated solutions set. In addition, most studies consider only single time window per patient in which he will be available to be visited by the assigned caregiver. 
  
  This chapter deals with the \ac{HHCRSP} with multi-objectives and multiple availability periods of patients.  The goal is to minimize the travel and waiting times as well as to balance caregivers’ workload. Two approaches are adopted to tackle the multi-objective case: Pareto based, where \ac{NSGA-II} is implemented  and decomposition based, where \ac{MOEA/D} is implemented. An efficient hybrid algorithm is proposed by performing \ac{NSGA-II} and then \ac{MOEA/D} successively, which would improve  solutions quality found by NSGA-II. 
  
  Performance  indicators  proved that the hybrid algorithm outperforms \ac{MOEA/D} and \ac{NSGA-II} algorithms in terms of solutions quality while \ac{MOEA/D} solves instances faster.  \ac{MOEA/D} is more suitable if the decision maker is not interested in having all Pareto optimal solutions. Besides, the benefit of using multiple availability periods for patients is highlighted by the significant improvement in the quality  of the solutions after carrying out a comparison with single time window per patient.
  
  Future works could be addressed to adapt these approaches to deal with multiple services operations and even considering the interdependence between them since patients need several care activities per day. In addition, it would be interesting to compare the performance of the hybrid algorithm with some indicators-based approaches.

\chapter*{General conclusion and perspectives}  
\addcontentsline{toc}{chapter}{General conclusion and perspectives}
      
    \markboth{General conclusion and perspectives}{} 
     
       \ac{HHC} aims to assist patients at home and to help them to live with greater independence, avoiding hospitalization or admission to care institutions. Providing health care services  is expected to generate advantages such as, a decrease in hospital admissions, a decrease in hospitalization duration and the ability for patients to remain in their homes and receive care and assistance.  \ac{HHC} companies face two major issues when planning schedules, namely assigning qualified caregivers to requested services and defining routes for caregivers to visit patients.
    
    Establishing plannings should be done taking into consideration  constraints such as patients’ availability periods and caregivers’ qualification and optimizing one or more criteria. In addition, the planning is affected by the type of the model (deterministic, stochastic, multi-objectives). The deterministic models are not robust since the uncertainties that may arise are not considered. In multi-objective models, the decision maker has more choices between solutions as the objective functions are simultaneously optimized. From this perspective, the research work in this thesis focuses on modeling and proposing solution methods to solve the home health care routing and scheduling problem, which extend the \ac{VRP} with time windows. 
     
     In chapter \ref{ChapterLR}, we analyzed the works carried out in the literature, which cover a very large variety concerning the constraints considered, the optimized objectives and the methods used to solve the problem. The first limitation in previous works is the no-consideration of multiple time windows for patients as well as most of these studies that consider the synchronization of the services are limited only to double. The second issue is related to uncertainty of parameters, which
is only dealt with in a few papers. The predefined schedule should be adapted for any change in practical situations. Otherwise, there will probably be delays in the services for patients who have not yet been visited, which will cause their dissatisfaction. The third limitation relates to the solution approaches used to solve the \ac{HHCRSP} multi-objective optimization model. Indeed, in the most cases, this model is solved using aggregation techniques, but assigning weights to conflicting goals is a difficult task.

In chapter \ref{ChapterDM}, we have proposed a deterministic version of the \ac{HHCRSP} that can deal with an arbitrary number of time windows and services per patient as well as taking into account of the  synchronization of simultaneous services and skill requirement. We proposed a new mathematical model as well as a \ac{GVNS}  based heuristic, which is combined with a strategy to solve the model with large instances. \ac{GVNS}  generates solutions (caregivers' routes and their assignment) while the strategy selects, for each patient, a time window and ensures the synchronization of simultaneous services.

In chapter \ref{ChapterSM}, we have considered the uncertainty of travel and service times. We proposed two SPR models to cope with the uncertainty of parameters. In the first one, the recourse is defined as a penalization of tardiness of services operation and a remuneration for caregivers’ overtime. In the second one, we assumed that patients’ times are hard/fixed and must be respected. Therefore, we used another recourse defined as skipping a patient if his time windows will not be respected. We showed through computational results the adequacy of the \ac{GA} to embed the simulation to solve the SPR model.

In chapter \ref{ChapterMOMTW},  we used algorithms designed to approach the Pareto front rather than using the aggregation techniques.  The decision maker will be involved a posteriori to select a solution from not dominated solutions, which avoid assigning weights to objective functions. Caregivers’ travel and waiting times as well the balance of workload are three criteria that are simultaneously optimized.   Approaches based on Pareto and decomposition, with multi-objective evolutionary algorithms are used to solve approximate Pareto front. Computational results and performance measures inferred that the \ac{MOEA/D} algorithm solved instances faster while the hybrid algorithm found solutions that better approximate the Pareto front. 

Following this works, several research perspectives are available to continue improving and/or extending our models. 

Our proposed models can be extended to cover other important features such as,  considering lunch break for caregivers. Thus, allowing caregivers to take breaks between the earliest and the latest of the break period where they can take lunch. Considering the continuity of cares is also an important feature, which allow  patients to be visited by the same caregiver and during approximately the same time. Our contributions are based on daily planning, considering medium and long horizon plannings  will lead to a better optimization of scheduling and routing of caregivers as well as to avoid defining daily caregivers plannings.

In the stochastic models, we considered that the uncertainty of  the parameters  related to travel and service times follows the normal law. However, if the distributions of parameters exposed to the uncertainty  are unknown, it would be better to use the robust optimization. The uncertainty in this model is not stochastic, but rather deterministic and set-based. In addition, we considered the uncertainty of of travel and service times, but the demand could also be exposed to the uncertainties and future works might address this issue. 

We used a simulation-based optimization  to the solve the SPR model  since computing the expected real value by an explicit mathematical formula is very complex. The computational results showed that the time resolution of some instances is relatively long because the simulation requires many iterations to find a good estimation. As future work, it would be interesting to compute  the expected real value by an explicit mathematical formula or to approximate it in order to avoid using the simulation. Otherwise, one might think of machine learning methods in the hope of reducing the number of iterations in the simulation method to estimate the expected value of each solution.

One of the main problems is the lack of generic benchmark instances that can be used to compare the performance of new methods and models. In most cases, especially when new features are considered, researchers randomly generate new instances to test their models and methods, as the instances in the literature are designed for classical problems such as TSP, CVRP, \ac{VRPTW}...etc. This problem can be tackled by generating new generic reference instances that contain the features studied in the literature.

     \cleardoublepage

\let\clearpage\relax
  
  

\begin{spacing}{1}
	\setlength{\bibsep}{0.9pt}
	
	\nocite{*}   
		
	\bibliographystyle{abbrv}
		\markboth{Bibliography}{}
	
	\bibliography{document}
	
   
 \newpage	 

\end{spacing}

\end{document}